\documentclass[final]{elsarticle}

\usepackage{lineno,hyperref}

\modulolinenumbers[5]

\usepackage[letterpaper,margin = 0.7 in]{geometry}

\usepackage{graphicx}

\usepackage{bm}

\graphicspath{{figure/}}

\usepackage{amsmath}

\usepackage{mathrsfs}

\usepackage{amsthm}

\usepackage{amssymb}

\usepackage{relsize}

\DeclareMathOperator*{\argmin}{arg\,min}

\usepackage{subcaption}

\usepackage{setspace}

\usepackage{tikz}

\usepackage{tabularx}

\usetikzlibrary{shapes.geometric, arrows}
\tikzstyle{startstop} = [rectangle, rounded corners, minimum width=1.5cm, minimum height=1cm,text centered, draw=black, fill=white]
\tikzstyle{io} = [trapezium, trapezium left angle=80, trapezium right 
angle=100, minimum width=1.0cm, minimum height=0.8cm, text centered, draw=black, fill=white]
\tikzstyle{process} = [rectangle, minimum width=1.5cm, minimum height=1cm, text centered, draw=black, fill=white]
\tikzstyle{decision} = [diamond, minimum width=1.5cm, minimum height=1cm, text centered, draw=black, fill=white]
\tikzstyle{arrow} = [thick,->,>=stealth]
\usetikzlibrary{automata,positioning}

\usepackage{framed,multirow}

\usepackage{amssymb}

\usepackage{latexsym}

\usepackage{url}

\usepackage{xcolor}

\definecolor{newcolor}{rgb}{.8,.349,.1}

\usepackage[utf8]{inputenc}
\usepackage{times}

\usepackage[ruled, vlined, linesnumbered, resetcount]{algorithm2e}

\theoremstyle{plain}

\theoremstyle{definition}

\usepackage[textsize=tiny]{todonotes}

\makeatletter
\newcommand{\algorithmfootnote}[2][\footnotesize]{%
	\let\old@algocf@finish\@algocf@finish% Store algorithm finish macro
	\def\@algocf@finish{\old@algocf@finish% Update finish macro to insert "footnote"
		\leavevmode\rlap{\begin{minipage}{\linewidth}
				#1#2
		\end{minipage}}%
	}%
}
\makeatother

\journal{Computer Methods in Applied Mechanics and Engineering}

\begin{document}
	
	\begin{frontmatter}

		\title{Projection pursuit adaptation on polynomial chaos expansions}%

		\author[1]{Xiaoshu Zeng\corref{cor1}}
		\ead{xiaoshuz@usc.edu}
		
		\author[1]{Roger Ghanem}
		\ead{ghanem@usc.edu}

		\cortext[cor1]{Corresponding author}

		%\fntext[fn1]{This is author footnote for second author.}  
		
		\address[1]{University of Southern California, 210 KAP Hall, Los Angeles, CA90089, United States}

		\begin{abstract}
			%%%
			The present work addresses the issue of accurate stochastic approximations in high-dimensional parametric space using tools from uncertainty quantification (UQ). The basis adaptation method and its accelerated algorithm in polynomial chaos expansions (PCE) were recently proposed to construct low-dimensional approximations adapted to specific quantities of interest (QoI). The present paper addresses one difficulty with these adaptations, namely their reliance on quadrature point sampling, which limits the reusability of potentially expensive samples. Projection pursuit (PP) is a statistical tool to find the ``interesting'' projections in high-dimensional data and thus bypass the curse-of-dimensionality. In the present work, we combine the fundamental ideas of basis adaptation and projection pursuit regression (PPR) to propose a novel method to simultaneously learn the optimal low-dimensional spaces and PCE representation from given data. While this projection pursuit adaptation (PPA) can be entirely data-driven, the constructed approximation exhibits mean-square convergence to the solution of an underlying governing equation and is thus subject to the same physics constraints. The proposed approach is demonstrated on a borehole problem and a structural dynamics problem, demonstrating the versatility of the method and its ability to discover low-dimensional manifolds with high accuracy with limited data. In addition, the method can learn surrogate models for different quantities of interest while reusing the same data set.
			%%%%
		\end{abstract}
		
		\begin{keyword}
			Surrogate modeling\sep
			polynomial chaos expansion\sep
			projection pursuit\sep
			high-dimensional models\sep
			dimension reduction\sep
			data-driven
		\end{keyword}

	\end{frontmatter}

%%%%%%%%%%%%%%%%%%%%%%%%%%%%%%%%%%%%%%%%%%%%%%%%%%%%%%%%%%%%%%%

%%%%%%%%%%                  New Section                %%%%%%%%

%%%%%%%%%%%%%%%%%%%%%%%%%%%%%%%%%%%%%%%%%%%%%%%%%%%%%%%%%%%%%%%

%\tableofcontents

\section{Introduction}
% Things to cover in this section
% 1. Curse-of-dimensionality in UQ and ML
% 2. Basis adaptation in PCE
% 3. Projection pursuit regression
% 4. Interpretability after combing these two methods 
% 5. Uncertainty of the surrogate model? Least squres lead to probabilisitic model of PCE coefficients

Recent developments in sensing have enabled highly resolved characterization of physical systems, justifying the drive for computational models of increasingly higher resolution and constantly pushing computational resources to the limit of hardware capacity. This computational bottleneck challenge is further exacerbated when uncertainty is acknowledged, which, at least conceptually, necessitates the simultaneous contemplation of a number of plausible realities. Allowing uncertainty to infiltrate computational models, however, can serve, with the proper perspective, to alleviate the computational burden, allowing consistency between numerical accuracy, statistical confidence, and experimental paucity. Response prediction can thus be achieved with the help of statistical surrogate models, which we achieve in the present paper by augmenting machine learning (ML) paradigms with an uncertainty quantification (UQ) dimension. Typically, these surrogate models are trained on replicas obtained by evaluating the physical model over a set of input parameters. Each of these evaluations will typically involve an expensive high-fidelity model, which is typically described by a large number of parameters. This combination of high resolution and high parameterization continues to challenge the accuracy of surrogate modeling. Similarly, in a big data set, the data of various control variables are generated such that comprehensive information on the quantity of interest (QoI) is considered in the predictive model. Usually, the number of data required for prediction with good accuracy increases with the number of parameters. Surrogate modeling in high-dimensional parameter space is thus ubiquitous, remaining at the forefront of predictive science. 

UQ is a tool to characterize the uncertainty in specific responses or outputs, given that the input parameters of the system are uncertain \cite{ghanem1991, ghanem1999propagation, sudret2007uncertainty, najm2009uncertainty, soize2017uncertainty,chen2017, zeng2017}. Generally, the uncertainties in the input parameters are represented by random variables that are characterized by their probability distributions. Usually, the supports of the parameters are first provided. This information is sometimes complemented by marginal distributions for each parameter. Only on rare occasions is joint information on the input parameters available at the outset. Such information may be weaved into the statistical model following data-driven updating. In all cases, the output is computed over a set of input parameters that, in principle, spans the support of their joint probability measure in parameter space. In a standard UQ setting, the probabilistic model of the output can be quantified by propagating the uncertainties from the input space through the computational model. Surrogate models, on the other hand, aim to construct a mapping from input parameters to output variables. The ensuing representation can be used as a surrogate to efficiently map the support of the input parameters into the support and distribution of the quantities of interest. In the literature, it is not uncommon to use UQ models as an efficient surrogate for further analyses that require repeated evaluation of the computational model, for example, in \cite{marzouk2007stochastic, marzouk2009dimensionality, marzouk2009stochastic, sudret2017surrogate}. The quality of these surrogates is clearly paramount. Polynomial chaos representations \cite{ghanem1991} can be shown to be faithful mathematical representations that converge in the mean-square sense to the true solution of the underlying governing stochastic equations.

For high-dimensional UQ problems, the standard Monte Carlo (MC) sampling has the attractive feature that its convergence is independent of the number of random parameters, thus can avoid the curse-of-dimensionality. However, it has slow convergence in general, and the rate is dependent on the complexity of the solution being sought. Thus the development of UQ methods such as polynomial chaos expansion (PCE) \cite{ghanem1991, ghanem1999}, generalized polynomial chaos expansion (gPCE) \cite{xiu2002, xiu2002-2}, Gaussian processes \cite{mackay1998introduction, seeger2004gaussian, bilionis2012multi}, probability density evolution methods \cite{li2009, chen2016, zeng2017, chen2017}, manifold methods \cite{soize2021probabilistic, zhang2021stochastic, giovanis2020data}, Wiener path integral technique\cite{kougioumtzoglou2012analytical, psaros2019wiener}, and others. PCE expands the QoI on a set of multivariate Hermite polynomial bases where the approximation has a mean-squared sense convergence with respect to the maximum order or degree of the polynomial bases. Its solid mathematical foundation, accuracy, simplicity in application, and convergence property gained much attention from science and engineering applications \cite{ghanem1991,ghanem1999, ghanem1999propagation, najm2009uncertainty}. However, the number of PCE terms in the expansion grows factorially with the parameters' dimension and the maximum polynomial order, leading to the curse-of-dimensionality. In a more common non-intrusive implementation, the PCE coefficients associated with each term are computed as a multi-dimensional integral, where the dimension is equal to the parameters' dimension. One can typically use numerical quadrature \cite{babuvska2007, maitre2010} to avoid excessive computation arising from the model evaluation to compute the PCE coefficients. Sparse quadrature \cite{ smolyak1963, gerstner1998, novak1999, maitre2010} can reduce the cost without sacrificing accuracy. Another family of methods, for example, least angle regression \cite{blatman2010, blatman2011}, and compressive sensing \cite{doostan2011, sargsyan2014, hampton2015compressive}, propose a sparse representation by expanding QoI on a selection of significant polynomial bases. These methods attenuate the computational burden but still face challenges if the dimension is excessively large or the computational model is expensive. Active subspace \cite{constantine2014, constantine2015} is a dimension reduction technique that finds the important directions with the greatest variability by an eigendecomposition of the covariance matrix of the gradient. In practice, the gradient of the QoI usually takes much work to obtain. The basis adaptation method \cite{tipireddy2014} in the framework of PCE is another dimension reduction technique that finds directions adapted to low-dimensional or scalar QoI. The idea is based on rotating the Gaussian input variables with a proper rotation matrix so that the QoI can be represented in a low-dimensional space. The construction of the rotation matrix requires the first-order PCE coefficients, which can be obtained by a first-order pilot PCE. Basis adaptation has been applied to practical engineering problems \cite{thimmisetty2017, ghauch2019} that reduce the dimension to only a few. This method can have slow convergence, resulting in relatively high reduced dimensions. Thus, the authors developed the accelerated basis adaptation method \cite{zeng2021accelerated}. The basis adaptation method requires the construction of low-dimensional adapted PCEs, where the PCE coefficients are usually computed via sparse quadrature to reduce computational cost. In addition, the method is adapted to scalar or low-dimensional QoIs.

In ML, dimension reduction techniques, for example, feature selection \cite{cai2018feature, murphy2012machine}, matrix factorization \cite{howley2005effect, reddy2020analysis}, manifold learning \cite{mcinnes2018umap, lin2008riemannian},  and autoencoder methods \cite{han2018autoencoder, baldi2012autoencoders}, have been able to increase the performance when the dimension of the data is high. Among these, the principal component analysis (PCA) is arguably the most widely used method, useful for discovering components that have ranked importance and can capture the data's salient structure. Projection pursuit (PP) \cite{friedman1974projection, huber1985projection, friedman1987exploratory} is another dimension reduction technique used to find the most ``interesting'' projections in high-dimensional data that can overcome the curse-of-dimensionality. The low-dimensional interesting projections are found by maximizing a so-called projection index stage-wisely. Different definitions of projection index lead to specific methods, for example, PCA and discriminant analysis. 
PP has received attention in ML in the recent two decades. For example, in \cite{lee2005projection}, new projection indices were proposed to provide low-dimensional projections for efficient supervised classification. For the classification of complex data, a projection pursuit constructive neural network was proposed to discover the simplest models in \cite{grochowski2008projection}. Based on the Jensen-Shannon divergence, a supervised projection pursuit was proposed in \cite{barcaru2019supervised} for better classification and visualization and proved superior to PCA. In \cite{grear2021molecular}, PP was used to develop a novel recurrent neural network for discriminant analysis. These works all belong to the scope of the classification. The combination of statistics and ML methods to improve the prediction performance is not unusual, for example, the Bayesian neural networks \cite{olivier2021bayesian, yang2021b, leibig2017leveraging}.
One of the extended ideas of PP is the projection pursuit regression (PPR) \cite{friedman1981projection}. In PPR, the response surface is approximated by a sum of univariate smooth functions of linear combinations, represented by projections of the input parameters (or predictors). The projections in PPR are found in a forward stage-wise strategy until new projections cannot improve the approximation significantly. Since the number of projections can be taken arbitrarily large, with a proper choice of smooth functions, the PPR model is able to approximate any continuous function in the space spanned by the parameters \cite{hastie2009elements}. Such a PPR model is thus a universal approximator. PPR has been applied to many practical applications, for example, in \cite{qianjian2011application, ferraty2013functional, durocher2015nonlinear,  cui2017assessment}.

The basis adaptation in PCE and the PPR have many shared features. They find the important or interesting directions in the input parameter space adapted to the output. The rotation matrix in adaptation corresponds to the projections in PPR. In contrast, the PCE of the adapted variables in adaptation corresponds to the smooth functions of the projected variables. However, these two methods are different in several aspects. The rotation matrix in adaptation is constructed based on a first-order PCE, and the directions are mutually orthogonal. At the same time, the projections in PPR are found one at a time, and no restrictions are applied among the projections. In addition, the functional links input parameters to output response in adaptation are Hermite polynomial bases, while in PPR, one must find proper smooth functions suited for the application. Also, the PCE coefficients in adaptation are usually found by sparse numerical quadrature (requiring the evaluation of the model on specific parameters) to reduce the computational cost. In contrast, the projections and smooth functions in PPR are discovered based on given data. Then, the advantages of these two methods can be summarized as follows: (1) the adaptation method uses Hermite polynomial bases that can represent any variable in the parameters space, and it has mean-squared convergence to the maximum polynomial order; (2) the PPR model can find projections and smooth functions simultaneously in a data-driven manner without querying the model.  

In this paper, we integrate the basis adaptation method in PCE and the PPR model to take advantage of both methods. Specifically, we use PCEs as smooth functions in the PPR model. In a naive integration, the PCEs work as univariate functions on each projected variable. To take advantage of the completeness and convergence properties of the multivariate Hermite polynomial bases, we will maintain the orthogonality of the projections and use a multivariate PCE on all the projected variables in the approximation. Switching from an additive of many univariate models to a multivariate model, in general, will be less cost-efficient. However, the purpose is to find a low-dimensional representation of the output; thus, the number of projections considered is supposed to be small. In such circumstances, seeking a more accurate multivariate representation is reasonable. We refer to the integrated multivariate model as the projection pursuit adaptation (PPA) model. Note that the projections construct the rotation matrix in the PPA model; thus, the multivariate PCE on the projected variables can be transformed back to the original physical space for interpretation purposes. The proposed model will be compared with the basis adaptation and the PPR methods. The paper is organized as follows: Section \ref{sec:adaptation} introduces the basis adaptation method, including the accelerated algorithm, in the framework of PCE; Section \ref{sec:ppa} first introduces the detailed procedure of the PPR, then derives the integrated PPA model with a detailed procedure also given; the applications on a borehole model and a space structure are shown in Section \ref{sec:application}; conclusions are finally presented in Section \ref{sec:conclusions}.

\section{Basis adaptation on polynomial chaos expansions}
\label{sec:adaptation}
\subsection{Polynomial chaos expansions}
Consider a probability space $(\Omega, \mathcal{F}, \mathbb{P})$ with $\Omega$ the sample space, $\mathcal{F}$ the $\sigma$-algebra on $\Omega$, and $\mathbb{P}$ the probability measure. Let $\bm{\xi} = (\xi_1, \ldots, \xi_d) \in \mathbb{R}^d$ be a set of independent standard Gaussian variables on $\Omega$. Then, the closed linear span (this allows for the case where $d\rightarrow\infty$) of $\bm{\xi}$ defines a Gaussian Hilbert space, denoted as $\mathcal{H}$. Denoting the space of second order random variables on $(\Omega, \mathcal{F}, \mathbb{P})$ by $L^2(\Omega)$, an inner product on $L^2(\Omega)$ can be introduced as 
\begin{equation}
	\langle u, v \rangle_{\mathcal{H}} = \mathbb{E} [u v],\quad \forall u,v\in L^2(\Omega),
\end{equation}
where $\mathbb{E}[\cdot]$ is the mathematical expectation with respect to the probability measure $\mathbb{P}$. Let $\bm{\alpha} = (\alpha_1, \ldots, \alpha_d) \in \mathcal{J}_d := (\mathbb{N}_0)^d$ be a $d$-dimensional multi-index with its length defined as $|\bm\alpha| = \sum_{i=1}^{d} \alpha_i$,  $h_{\bm{\alpha}}$ be multivariate Hermite polynomials, and $\psi_{\bm{\alpha}} = h_{\bm{\alpha}}/ \lVert h_{\bm{\alpha}}\lVert_2$ be the normalized Hermite polynomials, then $\{ \psi_{\bm{\alpha}},\, \bm{\alpha} \in \mathcal{J}_d \}$ is a complete orthonormal basis on $L^2(\Omega)$. We also let $\mathcal{F}(\mathcal{H})$ denote the $\sigma$-algebra generated by $\mathcal{H}$. Then we can represent any random variable $Y \in L^2 (\mathcal{F}(\mathcal{H}))$ as the following polynomial chaos expansion (PCE) \cite{ghanem1991}
\begin{equation}
	\label{eq:PCE}
	Y(\bm{\xi}) = \sum_{\bm{\alpha} \in \mathcal{J}_d} Y_{\bm{\alpha}} \psi_{\bm{\alpha}} (\bm{\xi}) \,.
\end{equation}
In Eq~\eqref{eq:PCE}, $Y_{\bm{\alpha}}$ is the PCE coefficient associated with the basis $\psi_{\bm{\alpha}}$.

In a computational context, a truncated PCE is usually employed where the maximum order of the polynomials is specified. Suppose the maximum order of the PCE is $p$, then the order $p$ PCE of $Y$ can be written as,
\begin{equation}
	\label{eq:PCEp}
	\widehat{Y}_p(\bm{\xi}) = \sum_{\bm{\alpha} \in \mathcal{J}_{d,p}} Y_{\bm{\alpha}} \psi_{\bm{\alpha}} (\bm{\xi})\,,
\end{equation}
where $\mathcal{J}_{d,p} \subset \mathcal{J}_d$ is the subset of indices of maximum length $p$ (i.e., $|\bm{\alpha}| \leq p$). The above expression, $\widehat{Y}_p(\bm{\xi})$ converges to $Y$ in the mean-squared sense \cite{janson1997} as $p \rightarrow \infty$. With abuse of notation, we will henceforth use $Y(\bm\xi)$ to denote $\widehat{Y}_p(\bm\xi)$.
The total number of PCE terms in Eq~\eqref{eq:PCEp} is
\begin{equation}
	P = {d+p \choose p} = \frac{(d+p)!}{d! p!}\,,
\end{equation}
which grows factorially with the dimension of the parameter space $d$ and the maximum order $p$. In a non-intrusive implementation of the PCE \cite{ghanem1991, babuvska2007, maitre2010}, the PCE coefficients $Y_{\bm{\alpha}}$ are computed by projecting $Y(\bm{\xi})$ on the space spanned by the Hermite polynomial basis $\{\psi_{\bm{\alpha}}(\bm{\xi}),\, \bm{\alpha} \in \mathcal{J} _{d, p} \}$ (referred to as the original space) as
\begin{equation}
	\label{eq:project}
	Y_{\bm{\alpha}} = \langle Y, \psi_{\bm{\alpha}}\rangle\,, \quad \bm{\alpha} \in \mathcal{J}_{d,p}\,.
\end{equation}
The projection can be calculated by, for example, the sampling-based Monte Carlo Simulation (MCS). In this setting, a set of samples $(\Xi, \mathcal{Y}) = \left\{ \left( \bm{\xi}^{(i)}, Y^{(i)} \right) \right\}_{i=1}^{N}$ is evaluated. These samples are then used to estimate the PCE coefficients $Y_{\bm{\alpha}}$ by minimizing the squared $l^2$ norm of the residual
\begin{equation}
	\left\lVert Y - \widetilde{Y} \right\lVert_2^2 = \sum_{i=1}^{N} \begin{bmatrix}
		Y(\bm{\xi}^{(i)}) - \sum_{\bm{\alpha} \in \mathcal{J}_{d,p}} Y_{\bm{\alpha}} \psi_{\bm{\alpha}} (\bm{\xi}^{(i)})
	\end{bmatrix}^2\,.
\end{equation}
The optimization problem can be solved by Ordinary Least Squares (OLS). Other methods to choose the sampling set $(\Xi, \mathcal{Y})$ could result in a weighted least squares problem. To obtain a stable and accurate estimation of the PCE coefficients, the number of required samples is typically of order $N \sim \mathcal{O}(P \ln P)$, and in some cases $N \sim \mathcal{O}(P^2 \ln P)$ \cite{cohen2017optimal}. With the factorial increased $P$ with respect to $d$ and $p$, the required computation effort becomes prohibitive when the dimension of parameters is high. 

Eq~\eqref{eq:project} can also be approximated using numerical quadrature  \cite{maitre2010}, where the physical model must be evaluated at specific quadrature points. Usually, the required samples could reduce a lot compared to MC samples. The sparse quadrature rules \cite{gerstner1998, novak1999, smolyak1963} can reduce the required quadrature points by adopting a sparse sampling algorithm. Other methods propose sparse PCE representations where only the dominant PCE terms are considered in the expansion \cite{blatman2010, doostan2011, sargsyan2014, schwab2011}. These methods partially attenuate the computation burden. However, the challenge of the curse-of-dimensionality remains as the underlying physical problem is mixed with high-dimensionality and a complex computational model.

\subsection{Basis adaptation}
\subsubsection{Rotation}
Basis adaptation is a dimension reduction technique for PCE in which the Gaussian input variables are rotated by a proper rotation matrix $\bm{A}$ such that the QoI can be represented on a low-dimensional space.

Consider a rotation matrix $\bm{A} \in \mathbb{R}^{d\times d}$ such that $\bm{AA}^T = \bm{I}$, and let $\bm{\eta}$ be the rotated variables defined as
\begin{equation}
	\bm{\eta} = \bm{A\xi}\,, \qquad \bm{\eta} = \begin{bmatrix}
		\bm{\eta}_r \\ \bm{\eta}_{\neg r}
	\end{bmatrix}\,,
\end{equation}
where $\bm{\eta}_r$ are the first $r$ components, and $\bm{\eta}_{\neg r}$ are the remaining $(d-r)$ components of $\bm{\eta}$. Then, instead of expanding $Y$ with respect to the original variable $\bm{\xi}$, we can equivalently represent $Y$ by the rotated variables $\bm{\eta}$ as
\begin{equation}
	\label{eq:PCEadapt}
	Y^{\bm{A}}(\bm{\eta}) = \sum_{\bm{\beta} \in \mathcal{J}_{d,p}} Y^{\bm{A}}_{\bm{\beta}} \psi_{\bm{\beta}}(\bm{\eta}) \,.
\end{equation}
Same as in the original space, the PCE coefficients $Y_{\bm{\beta}}^{\bm{A}}$ can be computed by projecting $\bm{Y}^{\bm{A}}(\bm\eta)$ on the space spanned by $\{ \psi_{\bm{\beta}}(\bm\eta),\, \bm{\beta} \in \mathcal{J}_{d,p} \}$ (referred to as the rotated space) as
\begin{equation}
	Y_{\bm{\beta}}^{\bm{A}} = \left\langle Y^{\bm{A}}, \psi_{\bm{\beta}}\right\rangle\,, \qquad \bm{\beta} \in \mathcal{J}_{d,p}\,.
\end{equation}
Define
\begin{equation}
	\psi_{\bm{\alpha}}^{\bm{A}}(\bm{\xi}) =\psi_{\bm{\alpha}}(\bm{A\xi})= \psi_{\bm{\alpha}}(\bm{\eta})\,.
\end{equation}
Then, by the equivalence of Eq~\eqref{eq:PCEp} and Eq~\eqref{eq:PCEadapt}, we have
\begin{equation}
	\label{eq:coeff_transfer}
	Y_{\bm{\alpha}}^{\bm A} = \sum_{\bm{\beta} \in \mathcal{J}_{d,p}} Y_{\bm\beta} \left\langle \psi_{\bm \beta}, \psi_{\bm \alpha}^{\bm A} \right\rangle,  \qquad \quad Y_{\bm{\alpha}} = \sum_{\bm{\beta}\in\mathcal{J}_{d, p}}Y_{\bm{\beta}}^{\bm A} \left\langle \psi_{\bm \beta}^{\bm A},\psi_{\bm \alpha} \right\rangle, 
	\ \qquad \bm{\alpha}\in \mathcal{J}_{d,p}\,.
\end{equation}
The above equation provides the foundation to transfer a PCE model from one representation to another and can be used to compare expansions with respect to different Gaussian variables.

\subsubsection{Classical Gaussian adaptation}
The rotation changes the variables the QoI expands on without changing the quality of the overall PCE, and without affecting the requisite computational burden. The idea of basis adaptation is to find a proper rotation matrix such that the probabilistic information of the QoI is concentrated on a low-dimensional space to achieve dimension reduction. Therefore, the construction of the rotation matrix becomes crucial.

The rotation matrix is constructed purely based on a first-order PCE expansion on the original space in classical Gaussian adaptation. If the first-order pilot PCE that requires a limited number of model evaluations is constructed, then the rotation matrix is created such that
\begin{equation}
	\label{eq:gauss_adapt}
	\eta_1 = \sum_{i=1}^{d}Y_{\bm e_i} \xi_i\,,
\end{equation}
where $\bm e_i$ is a $d$-dimensional multi-index with one at the $i$th location and zeros elsewhere, then $\{Y_{\bm{e}_i}\}_{i=1}^{d}$ is the set of first-order PCE coefficients or Gaussian coefficients in the original space. Eq~\eqref{eq:gauss_adapt} defines the first row of the rotation matrix $\bm{A}$, with which the first rotated variable captures the complete Gaussian components of the QoI. Next, the other rows of $\bm{A}$ are defined based on the sensitivity of the original variables. The sensitivities of the original variables $\bm{\xi}$ are quantified by the absolute values of their Gaussian coefficients. As a result, the second row of $\bm{A}$ is defined such that $\eta_2$ is the most sensitive variable of $\bm{\xi}$; the third row of $\bm{A}$ is defined such that $\eta_3$ is the second most sensitive variable of $\bm{\xi}$; and so on. Finally, a Gram-Schmidt procedure is applied to $\bm{A}$ to make it a rotation matrix. The complete procedure to construct the rotation matrix for the classical basis adaptation is summarized in Algorithm \ref{algo1}.
\begin{algorithm}[htb]
	\caption{Construction of rotation matrix by classical Gaussian adaptation}\label{algo1}
	
	Construct a pilot first-order PCE in the original space to find $\{Y_{\bm e_i}\} _{i=1}^d$;
	
	Construct the first row of the rotation matrix $\bm{A} \in \mathbb{R} ^{d \times d}$ such that
	\begin{equation*}
		\eta_1 = \sum_{i=1}^{d} Y_{\bm e_i} \xi_i \,.
	\end{equation*}

	Rank the first order coefficients by absolute value in descending order and record their indices in the original coordinates as $\{ \kappa_j \}_{j=1}^{d}$;
	
	For $j = 2, \ldots, d$ construct the $j^{\text{th}}$ row of $\bm{A}$ such that
	\begin{equation}
		\eta_j  =  \xi_{\kappa_{j-1}}.
	\end{equation}
	
	Perform Gram-Schmidt procedure on $\bm{A}$ to make it a rotation matrix.
\end{algorithm}

The Gaussian components are typically the most important information in many applications. Therefore, the first rotated variable $\eta_1$ conveys the key information of the QoI. Since the other rotated variables are initially constructed as the most sensitive variables in the original space in descending order, the rotated variables have a descending order of importance to the QoI, motivating the seeking of a reduced representation of the QoI by the first several rotated variables. The rows of the rotation matrix $\bm{A}$ are sometimes referred to as the adaptation or adapted directions.

Suppose the first $r$ rotated variables are adequate to represent the QoI, $Y$, then we can write
\begin{equation}
	\label{eq:adapted_PCE}
	Y^{\bm{A}_r}(\bm{\eta}_r) = \sum_{\bm{\beta} \in \mathcal{J}_{r,p}} Y^{\bm{A}_r}_{\bm{\beta}} \psi_{\bm{\beta}}(\bm{\eta}_r) \,,
\end{equation}
where $\bm{A}_r$ denotes the sub-matrix of the first $r$ rows of $\bm{A}$; $\mathcal{J}_{r,p}$ is a set of $r$-dimensional multi-indices. We will refer to Eq~\eqref{eq:adapted_PCE} as the $r$-dimensional adapted PCE, or simply $r$-dimensional adaptation. The PCE coefficient $Y_{\bm{\beta}}^{\bm{A}_r}$ can be calculated by projecting $Y^{\bm{A}_r}(\bm{\eta}_r)$ on the space spanned by the Hermite polynomial basis $\{ \psi_{\bm{\beta}}(\bm{\eta}_r),\, \bm{\beta} \in \mathcal{J}_{r,p} \}$ (referred to as the reduced rotated space) as
\begin{equation}
	\label{eq:project_rotate}
	Y_{\bm{\beta}}^{\bm{A}_r} = \left\langle Y^{\bm{A}_r}, \psi_{\bm{\beta}}\right\rangle\,, \qquad \bm{\beta} \in \mathcal{J}_{r,p}\,.
\end{equation}
Regardless of the methods used to compute Eq~\eqref{eq:project_rotate}, the physical model $Y$ must be evaluated on the original variables $\bm{\xi}$ instead of $\bm{\eta}_r$. Thus, the following approximation is employed
\begin{equation}
	Y^{\bm{A}_r}(\bm{\eta}_r) = Y(\tilde{\bm{\xi}}) = Y\left( \bm{A}^T \begin{bmatrix} \bm{\eta}_r & \bm{0} \end{bmatrix}^T \right) \,.
\end{equation}
That is, when projecting $\bm{\eta}_r$ back to the original coordinates, we assume that $\bm{\eta}_{\neg r} =\bm{0}$. The assumption introduces negligible errors when the first $r$ components of $\bm{\eta}$ are adequate to capture the probabilistic information of the QoI. 

The reduced dimension $r$ remains to be determined. A systematic way starts with $r=1$, where we construct an $r$-dimensional adaptation. Then, $r$ is sequentially increased. The procedure is stopped when the difference between two adaptations with successive dimensions is within a pre-specified tolerance. The measure of the difference between two adaptations, such as $r$ and $(r+1)$ dimensional adaptations, can be quantified by the distance of either their PCE coefficients or the associated Kernel Density Estimate (KDE) of the probability density functions (PDF) of the QoI. 
In the former, the $r$-dimensional adaptation is projected to the $(r+1)$-dimensional space, leading to a representation in the $(r+1)$-dimensional space (same idea as presented in Eq~\eqref{eq:coeff_transfer}). The projected PCE coefficients can then be compared to the ones of the $(r+1)$-dimensional adaptation.
In the latter, many MC samples are first generated based on the constructed $r$- and $(r+1)$-dimensional adaptations independently. These MC samples are then used to construct two KDE estimations of the PDF of the QoI, which are then compared to find a quantified difference.

In many applications, the classical basis adaptation can identify low-dimensional spaces adapted to the QoI, for example, \cite{thimmisetty2017, ghauch2019}. In \cite{zeng2021accelerated}, the authors proposed accelerated basis adaptation methods that can further enhance the accuracy of the reduced dimension $r$. The accelerated algorithm reduces the dimension when the classical adaptation still requires middle-to-high dimensions.

\subsubsection{Accelerated adaptation}
In \cite{zeng2021accelerated}, two acceleration methods of basis adaptation are proposed. The first method corrects the low-order information of the adapted PCE by utilizing the first-order pilot PCE constructed in the first step of adaptation, see Algorithm \ref{algo1}. The second method updates the rotation matrix by accounting for the information from higher dimensions and high-order adaptations while the dimension is sequentially increased. We will introduce the ideas behind these two methods briefly. For a detailed algorithm, please refer to \cite{zeng2021accelerated}.

Suppose the first-order pilot PCE in the original space is constructed as
\begin{equation}
	\label{eq:pilot}
	Y(\bm{\xi}) = \sum_{\bm{\alpha} \in \mathcal{J}_{d,1}} Y_{\bm{\alpha}} \psi_{\bm \alpha}(\bm{\xi}) = Y_{\bm{e}_0} + Y_{\bm{e}_1}\xi_1 + \cdots  + Y_{\bm{e}_d}\xi_d\,.
\end{equation}
Denoting the vector of Gaussian coefficients as $\bm{Y} ^{\bm{e}} _{\text{coeff}} = (Y_{\bm{e}_1}, \ldots, Y_{\bm{e}_d})$, then the $r$-dimensional adaptation of order up to one can be expanded as
\begin{equation}
	\label{eq:1d_adapt}
	{Y}^{\bm{A}_r}(\bm{\eta}_r) = \sum_{\bm{\beta}\in\mathcal{J}_{r,1}} Y_{\bm{\beta}}^{\bm{A}_r} \psi_{\bm{\beta}}(\bm{\eta}_r) = Y_{\bm{e}_0}^{\bm{A}_r} + Y_{\bm{e}_1}^{\bm{A}_r}\eta_1 + \cdots  + Y_{\bm{e}_r}^{\bm{A}_r}\eta_r\,.
\end{equation}
To measure the difference in the zero and first-order expansions, the PCE \eqref{eq:1d_adapt} can be projected to the original space to obtain the associated PCE coefficients as
\begin{equation}
	\widetilde{Y}_{\bm{e}_0}  = Y_{\bm{e}_0}^{\bm A_r}\,, \quad
	\widetilde{Y}_{\bm{e}_1}  = \sum_{i=1}^{r} Y_{\bm{e}_i}^{\bm A_r} A_{i1} \,, \quad
	\dots\,, \quad
	\widetilde{Y}_{\bm{e}_d}  = \sum_{i=1}^{r} Y_{\bm{e}_i}^{\bm A_r} A_{id} \,,
\end{equation}
where $A_{ij}$ denotes the $(i, j)$ component of the rotation matrix $\bm{A}$. When $r=1$, the projected coefficients can be computed explicitly as
\begin{equation}
	\widetilde{Y}_{\bm{e}_0}  = Y_{\bm{e}_0}^{\bm A_r}\,, \quad
	\widetilde{Y}_{\bm{e}_1}  = \frac{Y_{\bm{e}_1}^{\bm A_r} Y_{\bm{e}_1}}{\lVert \bm{Y^e}_{\text{coeff}} \lVert_2} \,, \quad
	\dots\,, \quad
	\widetilde{Y}_{\bm{e}_d}  = \frac{Y_{\bm{e}_1}^{\bm A_r} Y_{\bm{e}_d}}{\lVert \bm{Y^e}_{\text{coeff}} \lVert_2} \,.
\end{equation}
Defining ${\widetilde{\bm{Y}}}^{\bm{e}}_{\text{coeff}} = ( \widetilde{Y}_{\bm{e}_1}, \ldots,  \widetilde{Y}_{\bm{e}_d} )$, then the errors of the Gaussian coefficients are quantified as
\begin{equation}
	\bm{\epsilon}_{1} = \bm{Y^e}_{\text{coeff}} - \widetilde{\bm{Y}}^{\bm{e}}_{\text{coeff}} = \left( 1 - \frac{Y_{\bm{e}_1}^{\bm{A}_r}}{ \lVert \bm{Y^e}_{\text{coeff}} \lVert_2 } \right) \bm{Y^e}_{\text{coeff}}\,.
\end{equation}
We see that the arithmetic errors are proportional to the exact Gaussian coefficients $\bm{Y^e}_{\text{coeff}}$, implying that more significant errors appear in the terms with greater Gaussian coefficients. However, greater Gaussian coefficients imply greater sensitivities in the variable, and greater errors in the terms with greater sensitivity introduce more significant errors in the PCE. As dimension $r$ of the adaptation increases, the errors in Gaussian coefficients decrease and will recover to the exact values when $r=d$. In order to eliminate the errors in the zero and first-order information of the adapted PCE, the authors proposed to correct the zero and first-order PCE coefficients on the adapted space by the exact zero and first-order PCE coefficients on the original space. The idea is accomplished by projecting Eq~\eqref{eq:pilot} to the $r$-dimensional adapted space and obtaining the projected coefficients as
\begin{equation}
	\widetilde{Y}^{\bm A_r}_{\bm{e}_0}  = Y_{\bm{e}_0}\,, \quad
	\widetilde{Y}^{\bm A_r}_{\bm{e}_1}  = \sum_{i=1}^{d} Y_{\bm{e}_i} A_{1i} \,, \quad
	\dots\,, \quad
	\widetilde{Y}^{\bm A_r}_{\bm{e}_r}  = \sum_{i=1}^{d} Y_{\bm{e}_i} A_{ri} \,.
\end{equation}
These coefficients are the exact zero and first-order PCE coefficients in the $r$-dimensional adapted space and can be used to correct the coefficients in Eq~\eqref{eq:1d_adapt}. Then, the zero and first-order information are guaranteed to be accurate.

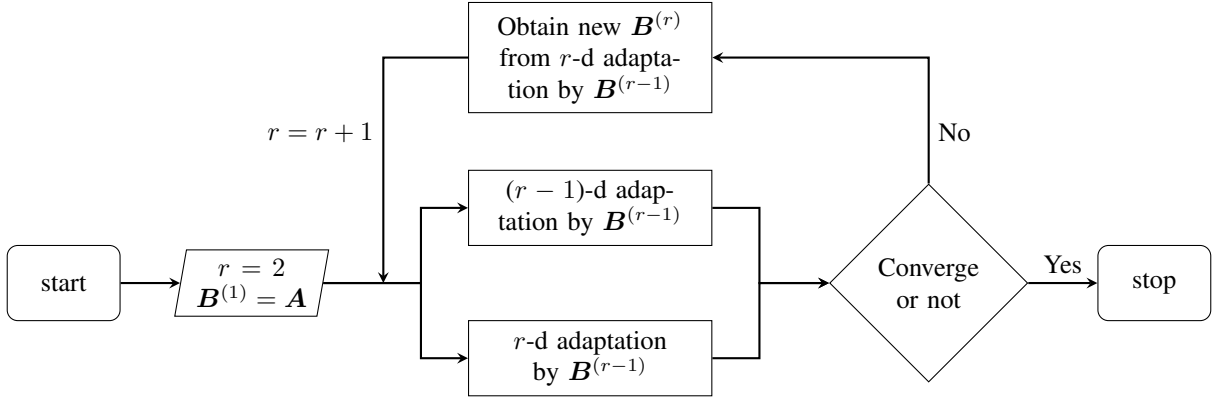
\begin{figure}[htb]
	\centering
	\begin{tikzpicture}
		[%
		>=stealth,
		node distance=3.0cm,
		on grid,
		auto
		]
		\node (start) at (-7, 0) [startstop] {start};
		\node (in) at (-4.5, 0) [io, text width=1.5cm] {$r=2$  $\bm{B}^{(1)}=\bm A$};
		\coordinate (aux1) at (-2.25,0);
		\coordinate (aux0) at (-2.75,0);
		\node (1d) at (0, 1)[process,  text width=3cm] 
		{($r-1$)-d adaptation by $ \bm{B}^{(r-1)}$};
		\node (2d) at (0, -1) [process, text width=3cm] {$r$-d adaptation by $\bm{B}^{(r-1)}$} ;
		\coordinate (aux2) at (2.25,0);
		\node (dec) at (4.5,0) [decision, text width=1.5cm] {Converge or not};
		\node (stop) at (7.5,0) [startstop] {stop};
		\node (new_A1) at (0,3) [process,  text width=3cm, yshift=0cm] {Obtain new $\bm{B}^{(r)}$ from $r$-d adaptation by $\bm{B}^{(r-1)}$};
		\draw [arrow] (start) -- (in);
		\draw [arrow] (in) -- (aux1) |- (1d);
		\draw [arrow] (in) -- (aux1) |- (2d);
		\draw [arrow] (2d) -| (aux2) -- (dec);
		\draw [arrow] (1d) -| (aux2) -- (dec);
		\draw [arrow] (dec) -- node[anchor=south]{Yes}(stop);
		\draw [arrow] (dec) |- node[anchor=west, yshift=-1.0cm]{No} (new_A1);
		\draw [arrow] (new_A1) -| node[anchor=east, yshift=-1.0cm]{$r=r+1$}(aux0);
	\end{tikzpicture}
	\caption{Workflow of the sequentially optimized adaptation method \cite{zeng2021accelerated}.} \label{workflow}
\end{figure}

For the second method, we noticed that with the increased $r$ in the adaptation, the QoI could be better represented by the adapted PCE. In addition, the $r$-dimensional adapted PCE retains the information of the first $r$ rotated variables up to order $p$. However, the rotation matrix in the classical adaptation only considers the first-order information and remains constant even when more probabilistic information becomes available. The second method in \cite{zeng2021accelerated} proposed to update the rotation matrix $\bm{A}$ as $r$ is sequentially increased. Specifically, suppose the $r$-dimensional adaptation from the rotation matrix $\bm{A}$ is available. In that case, the best one-dimensional adaptation close to the $r$-dimensional adaptation is discovered from an updated rotation matrix $\bm{B}$. The one-dimensional adaptation by $\bm{B}$ can be written as
\begin{equation}
	{Y}^{\bm{B}_1} (\zeta_1) = \sum_{ \bm{\gamma} \in \mathcal{J}_{1,p}} Y_{\bm{\gamma}}^{\bm{B}_1} \psi_{ \bm{\gamma}} (\zeta_1) = \sum_{ \bm{\gamma} \in \mathcal{J}_{1,p}}  Y_{\bm{\gamma}}^{\bm{B}_1} \psi_{ \bm{\gamma}}^{\bm B_1} (\bm\xi) \,,
\end{equation}
where $\bm{\zeta} = \bm{B}\bm{\xi}$ is the new adapted variable by $\bm{B}$. To compute the PCE coefficients of the new one-dimensional adaptation, one would usually require to evaluate the physical model $Y$ on some parameters. To avoid the additional computation, the evaluation is instead performed on the $r$-dimensional adapted PCE by $\bm{A}$, the best available model. For any choice of $\bm{B}$, one can easily construct the one-dimensional adapted PCE of the QoI. Then, an optimization procedure can be employed to find the best $\bm{B}$ that minimizes the difference between the one-dimensional adaptation by $\bm{B}$ and the $r$-dimensional adaptation by $\bm{A}$. This method can be combined with the first method that always corrects the zero and first-order PCE coefficients in the adapted PCEs. The combined method is named the sequentially optimized adaptation (SOA) with its workflow in Figure \ref{workflow}. Please refer to \cite{zeng2021accelerated} for details.

The basis adaptation is a great tool for dimension reduction in UQ. However, when it comes to computing the adapted PCE coefficients, for example, by Eq~\eqref{eq:project_rotate}, we typically use sparse quadrature rules to reduce the number of evaluations while having good accuracy at the same time. In this setting, the model must be evaluated at specific quadrature points. We will next introduce a novel data-driven method that can simultaneously find the optimal rotation matrix and the PCE coefficients.

\section{Projection pursuit adaptation}
\label{sec:ppa}
\subsection{Projection pursuit regression}
Projection pursuit (PP) intends to find the most ``interesting'' projections in high-dimensional data. These projections are discovered such that they reveal the most information embedded in the data. Therefore, the projections of PP and the adaptation directions by the basis adaptation method have the same purpose but are discovered with different algorithms. The advantages and disadvantages of these two methods are worth exploring. In the sense of adapting input to output, the projection pursuit regression (PPR) was introduced \cite{friedman1981projection}, where the projections and the functions from input space to output space are built simultaneously. 

In a general supervised learning setting, assume we have a data set $\{(\bm{x}^i, y^i)\}_{i=1}^{N}$ with $\bm{x}^i \in \mathbb{R}^{d}$ and $y^i \in \mathbb{R}$. Let $\bm{c}_{j}$, $j=1,\ldots, r$ be $d$-dimensional known vectors. Then the PPR model seeks to approximate $y^i$ as follows
\begin{equation}
	y^i \approx \sum_{j=1}^{r}f_j(\bm{c}_{j}^T \bm{x}^i)\,,
\end{equation}
where $\{f_j : \mathbb{R} \mapsto \mathbb{R}\}$ is a set of univariate smooth functions.  This model aims to find $(f_j, \bm{c}_{j})$ such that the approximation fits the data in some optimal sense. Specifically, the following optimization problem is solved to fit such a model with the given data
\begin{equation}
	\argmin_{(f_j, \bm{c}_{j})} \sum_{i=1}^{N} \left[y^i - \sum_{j=1}^{r} f_j (\bm{c}_{j}^T \bm{x}^i)\right]^2 \,.
\end{equation}
In this optimization problem, the number of required projections, $r$, is determined by a stage-wise greedy algorithm where we add a pair of $(f_j, \bm{c}_{j})$ at each stage and stop the procedure if additional information cannot improve the model fitting significantly. For example, suppose we are in stage $j$ and want to find the optimal $(f_j, \bm{c}_{j})$. Let the previously found parameters, $\{ (f_l, \bm{c}_l) \}_{l=1}^{j-1}$ be fixed, then the residual of the $i$th data point is obtained as
\begin{equation}
	\nu^i = y^i - \sum_{l=1}^{j-1} f_l(\bm{c}_l^T \bm{x}^i)\,.
\end{equation}
Then, the optimal $(f_j, \bm{c}_{j})$ can be found by solving the following optimization problem
\begin{equation}
	\label{eq:opt_ppr_1}
	\argmin_{(f_j, \bm{c}_{j})} \sum_{i=1}^{N} \left[\nu^i -  f_j (\bm{c}_{j}^T \bm{x}^i)\right]^2 \,.
\end{equation}
When the total residual of the data $\sum_{i=1}^{N} \nu^i$ is converged, we stop the procedure. 

An alternative to the greedy procedure is to iteratively update the pairs $(f_k, \bm{c}_{k})$ for $k=1, \ldots, j$ through backfitting. While this approach could result in significantly fewer terms in the model, the computation cost is increased, and the improvement in predictions is insignificant \cite{hastie2009elements}. In addition, the same performance can usually be achieved by adding more terms to the model. Therefore, we will focus on the stage-wise greedy procedure.

The optimization problem \eqref{eq:opt_ppr_1} remains to be solved in each stage. The solution can be obtained by an alternating optimization procedure, where the optimization switches between $f_j$ and $\bm{c}_j$ until converge. Specifically, if the vector $\bm{c}_{j}$ is given, then the argument of $f_j$ is a one-dimensional variable. The problem becomes finding the optimal $f_j$ such that the residual is minimized. This can be done by the regression method. On the other hand, assuming $f_j$ is given, one needs to find the optimal $\bm{c}_{j}$ to minimize the residual. This can be done by a Gauss-Newton search, a quasi-Newton method that keeps only the zero and first derivatives in a Taylor expansion. Let $\bm{c}_{j,\text{old}}$ be the current estimate of $\bm{c}_{j}$, then we can write the following Taylor expansion
\begin{equation}
	\label{eq:taylor}
	f_j(\bm{c}_{j}^T \bm{x}^i) \approx f_j(\bm{c}_{j,\text{old}}^T \bm{x}^i) + \dot{f}_j(\bm{c}_{j,\text{old}}^T \bm{x}^i) (\bm{c}_{j} - \bm{c}_{j,\text{old}})^T \bm{x}^i \,.
\end{equation}
Substituting Eq~\eqref{eq:taylor} into Eq~\eqref{eq:opt_ppr_1} (with $f_j$ given) leads to
\begin{equation}
	\label{eq:opt_ols}
	\argmin_{\bm{c}_{j}} \sum_{i=1}^{N} \left[\nu^i -  f_j (\bm{c}_{j}^T \bm{x}^i)\right]^2 
	\approx
	\argmin_{\bm{c}_{j}}
	\sum_{i=1}^{N} \dot{f}_j(\bm{c}_{j,\text{old}}^T \bm{x}^i) ^2 \left[ \left( \bm{c}_{j,\text{old}}^T \bm{x}^i + \frac{\nu^i - f_j(\bm{c}_{j,\text{old}}^T \bm{x}^i)}{\dot{f}_j(\bm{c}_{j,\text{old}}^T \bm{x}^i)} \right)  - \bm{c}_{j}^T \bm{x}^i \right]^2 \,.
\end{equation}
The optimization problem \eqref{eq:opt_ols} can be solved by the weighted least squares regression with target $u =  \bm{c}_{j,\text{old}}^T \bm{x}^i + \frac{\nu^i - f_j(\bm{c}_{j,\text{old}}^T \bm{x}^i)}{\dot{f}_j(\bm{c}_{j,\text{old}}^T \bm{x}^i)}$, regressors $\bm{x}^i$, weight $\omega =  \dot{f}_j(\bm{c}_{j,\text{old}}^T \bm{x}^i) ^2$, and no bias term. Let $\bm{W} \in \mathbb{R}^{N \times N}$ denote the diagonal matrix with entries $\omega$, $\bm{u} \in \mathbb{R}^{N}$ be the vector with entries $u$, and use the full data matrix $\bm{X} \in \mathbb{R} ^{N \times d}$, then $\bm{c}_j$ is found by solving the following
\begin{equation}
	\label{eq:ols}
	\argmin_{\bm{c}_{j}} \lVert  \bm{u} - \bm{X} \bm{c}_{j} \lVert ^2 _{\bm{W}} = \left( \bm{X}^T \bm{WX} \right) ^{-1} \bm{X}^T \bm{Wu} \,.
\end{equation}
With an updated $\bm{c}_{j}$, we find new projections $\bm{c}_{j}^T \bm{x}^i$ and refit $f_j$ to minimize the residual. Then $\bm{c}_{j}$ can be updated again with the new $f_j$. We repeat the alternating procedure until $(f_j, \bm{c}_{j})$ pair is converged. 

Although there is no restriction to the smooth functions $\{f_j\}$, it is better if the functions are once differentiable. PCE is a good fit here where each smooth function $f_j$ is a univariate PCE with respect to the projected variable $\bm{c}_{j}^T \bm{x}^i$. In this case, the variables $\bm{x}$ that can be arbitrarily distributed (based on prior knowledge or from the data) need to be first transformed to standard Gaussian distributions via, for example, inverse transform sampling or Rosenblatt transformation \cite{rosenblatt1952remarks}. Then, with a given $\bm{c}_{j}$, finding the optimal $f_j$ is equivalent to finding the PCE coefficients of $f_j$. One of the benefits of using PCEs as smooth functions is that the order of the PCE is adjustable, and we can gradually increase the polynomial orders until convergence is reached. Based on the flexibility of PCE, $r$ can be arbitrarily large to approximate any function in $L^2(\Omega)$, and the computational cost lies in constructing $r$ univariate PCEs. 

Note that there is no restriction on the projection in the optimization problem \eqref{eq:ols}, and each projection is discovered based on the current residual. The initial guess of $\bm{c}_j$ and the data determine where the ``optimal'' value land. Therefore, there are chances that the weighted least squares find the local optima and stop the stage-wise procedure when the local optima suggest there is not much room for improvement. In addition, there is no restriction on the projections. No coupling effects are considered among the projections since the smooth function $f_j$ only considers information on the $j$th projected variable. One of the advantages of such construction is that $\{f_j\}$ are univariate functions and are easy to compute with given projections and a given number of samples.

\subsection{Integration of projection pursuit and basis adaptation}
In the basis adaptation method, we already discovered that many scalar QoI could be represented in a low-dimensional space (\cite{ghauch2019, thimmisetty2017, zeng2021accelerated, tsilifis2017}). So the question is how to find the optimal low-dimensional space and the PCE on that space. The PPR method is an excellent tool for achieving that in a data-driven setting.

However, the original construction considers univariate representations in PPR with PCE employed as smooth functions. Therefore, it does not take advantage of the good mean-squared convergence property and the ability to represent any variable in the parameter space of the multivariate PCE construction. Accordingly, we will explore using multivariate PCE as a smooth function in PPR. In addition, the projections in this effort are assumed to be orthonormal, which has control of the exploring space of the next projection based on the current projections. Therefore, in this case, the projections are equivalent to the rows of the rotation matrix, and the smooth function is equivalent to multivariate PCE in the basis adaptation method. Hence, we refer to this method as the projection pursuit adaptation (PPA). 

Let $\bm{C}_r = \begin{bmatrix} \bm{c}_1 & \ldots & \bm{c}_r \end{bmatrix}^T$, with each row representing a projection and with the projections being mutually orthogonal. Also let $\bm{z}^i$ denote the projected variables such that
\begin{equation}
	\bm{z}^i = \begin{bmatrix} z_1^i & \ldots & z_r^i \end{bmatrix}^T =  \bm{C}_r \bm{x}^i =  \begin{bmatrix} \bm{c}_1 & \ldots & \bm{c}_r \end{bmatrix}^T \bm{x}^i \,, \qquad \bm{x}^i \in \mathbb{R}^d, \quad \bm{z}^i \in \mathbb{R}^r
\end{equation} 
The PPA method seeks the following model
\begin{equation}
	y^i \approx g_r(\bm{z}^i)\,,
\end{equation}
where $g_r(\cdot)$ is $r$-dimensional PCE. We want to find $(g_r, \bm{C}_r)$ such that the approximation fits the data well. Unlike PPR, where the regression model is an additive of many univariate functions, the PPA seeks an $r$-dimensional multivariate model. Same as in PPR, the following optimization problem is solved
\begin{equation}
	\label{eq:opt_ppa}
	\argmin_{(g_r, \bm{C}_r)} \sum_{i=1}^{N} \left[y^i - g_r (\bm{z}^i )\right]^2 \,.
\end{equation}
We will use a stage-wise greedy strategy to solve the optimization problem. In the $r$th stage, for example, the previously computed projections $\bm{C}_{r-1}$ are unchanged, and we want to get the optimal $\bm{c}_r$ that is appended to $\bm{C}_{r-1}$. While for the smooth function, the optimal $r$-dimensional  $g_r$ is computed independently of the previous information in $g_{r-1}$. That is, the stage-wise greedy strategy finds optimal $(g_r, \bm{c}_r)$ at each stage. Like PPR, $g_r$ and $\bm{c}_r$ are computed by an alternating optimization procedure, in which the optimization switches between $g_r$ and $\bm{c}_r$ until they converge. We can discover the optimal $r$-dimensional PCE $g_r$ by regression if $\bm{C}_r$ is given. While for a given $g_r$, we calculate the optimal $\bm{c}_r$ that is appended to the previously computed $\bm{C}_{r-1}$. The reason for not updating the whole $\bm{C}_r$ matrix at each stage is to avoid excessive computation and propose a convenient Gauss-Newton search algorithm for optimization. With a fixed $g_r$, the following problem remains to be solved
\begin{equation}
	\label{eq:opt_multid}
	\argmin_{\bm{c}_r} \sum_{i=1}^{N} \left[y^i - g_r (\bm{z}^i)\right]^2 \,,
\end{equation}
which is similar to Eq~\eqref{eq:opt_ppr_1} with $f_j$ fixed. However, the function $g_r(\cdot)$ here is a multivariate function. In order to use the Gauss-Newton search, let $\bm{c}_{r,\text{old}}$ be the current estimate of $\bm{c}_{r}$, the following multi-dimensional Taylor expansion is applied
\begin{equation}
	\label{eq:taylor_multid}
	\begin{split}
		&g_r(\bm{z}^i)  = g_r(\begin{bmatrix}
			z_1^i & \ldots & z_{r-1}^i & z_{r}^i
		\end{bmatrix}^T) = g_r(\begin{bmatrix}
		\bm{c}_{1} & \ldots & \bm{c}_{r-1} & \bm{c}_{r}
	\end{bmatrix}^T \bm{x}^i)\\
		\approx&
		g_r \left(\begin{bmatrix}
			\bm{c}_1 & \ldots & \bm{c}_{r-1} & \bm{c}_{r,\text{old}}
		\end{bmatrix}^T \bm{x}^i \right) 
		+
		\nabla g_r \left(\begin{bmatrix}
			\bm{c}_1 & \ldots & \bm{c}_{r-1} & \bm{c}_{r,\text{old}}
		\end{bmatrix}^T \bm{x}^i \right)  \\
		& \times \left( \begin{bmatrix}
		\bm{c}_1 & \ldots & \bm{c}_{r-1} & \bm{c}_{r}
		\end{bmatrix}^T - \begin{bmatrix}
		\bm{c}_1 & \ldots & \bm{c}_{r-1} & \bm{c}_{r,\text{old}}
		\end{bmatrix}^T \right) \bm{x}^i \\
		=&
		g_r \left(\begin{bmatrix}
			\bm{c}_1 & \ldots & \bm{c}_{r-1} & \bm{c}_{r,\text{old}}
		\end{bmatrix}^T \bm{x}^i \right) 
		+\frac{\partial g_r}{\partial z^i_{r, \text{old}} } \left(\begin{bmatrix}
			\bm{c}_1 & \ldots & \bm{c}_{r-1} & \bm{c}_{r,\text{old}}
		\end{bmatrix}^T \bm{x}^i \right) (\bm{c}_r - \bm{c}_{r,\text{old}})^T \bm{x}^i \\
		=&
		g_r \left(\bm{z}^{i}_{\text{old}} \right) + \frac{\partial g_r \left(\bm{z}^{i}_{\text{old}} \right)}{\partial z^i_{r, \text{old}} } (\bm{c}_r - \bm{c}_{r,\text{old}})^T \bm{x}^i\,,
	\end{split}
\end{equation}
where $\nabla g_r(\cdot)$ denotes the gradient of $g_r(\cdot)$; $\bm{z}^{i}_{\text{old}} = \begin{bmatrix}	\bm{c}_1 & \ldots & \bm{c}_{r-1} & \bm{c}_{r,\text{old}} \end{bmatrix}^T \bm{x}^i$; $z^i_{r, \text{old}} = \bm{c}_{r, \text{old}}^T \bm{x}^i$. In \eqref{eq:taylor_multid}, $\{\bm{c}_{j}\}$ for $j<r$ are the known previously obtained projections. Although the Taylor expansion is multi-dimensional, we avoid the gradient operations, thus replacing a multi-dimensional regression with a one-dimensional regression by updating only the $r$th projection. Substituting \eqref{eq:taylor_multid} to \eqref{eq:opt_multid} leads to the following optimization problem
\begin{equation}
	\label{eq:opt_ols_multid}
	\argmin_{\bm{c}_r}
	\sum_{i=1}^{N} \left( \frac{\partial g_r \left( \bm{z}^{i}_{\text{old}} \right)}{\partial z^i_{r, \text{old}} } \right) ^2 \left[ \left( \bm{c}_{r,\text{old}}^T\bm{x}^i + \frac{y^i - g_r \left(\bm{z}^{i}_{\text{old}} \right)} {\partial g_r \left(\bm{z}^{i}_{\text{old}} \right) / \partial z^i_{r, \text{old}}} \right)  - \bm{c}_r^T \bm{x}^i \right]^2 \,,
\end{equation}
Again, the optimization problem \eqref{eq:opt_ols_multid} can be solved by the weighted least squares regression with target $\hat{u} =  \frac{y^i - g_r \left(\bm{z}^{i}_{\text{old}} \right)} {\partial g_r \left(\bm{z}^{i}_{\text{old}} \right) / \partial z^i_{r, \text{old}}}$, regressors $\bm{x}^i$, weight $\hat{\omega} =  \left(\partial g_r \left( \bm{z}^{i}_{\text{old}} \right)/\partial z^i_{r, \text{old}} \right)^2$, and no bias term. Let $\hat{\bm{W}} \in \mathbb{R}^{N \times N}$ denote the diagonal matrix with entries $\hat{\omega}$, and $\hat{\bm{u}} \in \mathbb{R}^{N}$ be the vector with entries $\hat{u}$, then $\bm{c}_r$ can be found by a solution to the following optimization problem
\begin{equation}
\label{eq:ols_multid}
\argmin_{\bm{c}_r} \lVert  \hat{\bm{u}} - \bm{X} \bm{c}_r \lVert ^2 _{\hat{\bm{W}}} = \left( \bm{X}^T \hat{\bm{W}} \bm{X} \right) ^{-1} \bm{X}^T \hat{\bm{W}} \hat{\bm{u}} \,.
\end{equation}
Once $\bm{c}_r$ is computed, we apply a Gram-Schmidt procedure to all the projections to maintain orthogonality among them. Then, with the updated $\bm{c}_r$, we find new projections of $\bm{x}^i$ and refit $g_r$ to minimize the residual. Then $\bm{c}_r$ can be updated again with the new $g_r$. We repeat the alternating procedure until $(g_r, \bm{c}_r)$ pair is converged. 

The detailed PPA procedure is summarized in Algorithm \ref{algo:ppa}.
\begin{algorithm}
	\setstretch{1.25}
	\caption{Projection Pursuit Adaptation (PPA)}
	\label{algo:ppa}
	
	Let $r=1$;\\
	
	\While{$\{g_r(\bm{C}_r\bm{x}^i)\}$ is not converged}{
		Initial guess of $\bm{c}_r$; \\
		
		Project the original variables $\{\bm{x}^i\}$ to projected variables $\{\bm{z}^i\}$; \\
		
		Find the optimal $g_r$ (PCE) with respect to $\bm{z}^i$ by least squares; \\
		
		\While{$(g_r, \bm{c}_r)$ is not converged}{
			Compute corresponding $\hat{c}$, $\hat{\omega}$ in \eqref{eq:ols_multid}; \\
			
			Solve \eqref{eq:ols_multid} to obtain an updated $\bm{c}_r$; \\
			
			Perform Gram-Schmidt on $\bm{C}_r$ to maintain orthogonality;\\
			
			Update the projected variables $\{\bm{z}^i\}$; \\
			
			Find the updated $g_r$ (PCE) with respect to $\bm{z}^i$ by least squares; \\
		}
	$r = r+1$.
	}
\end{algorithm}

The PPA method always finds the optimal $r$-dimensional PCE adapted to the output or the QoI, starting from $r=1$. The optimal model is in the sense of optimal directions and optimal PCE. Same as in basis adaptation, the dimension $r$ is sequentially increased, resulting in PCE models that can better approximate the data. The procedure is stopped when adding a new projection cannot reduce the residual significantly.

Generally, fitting a multi-dimensional PCE model is much more expensive than fitting multiple univariate PCEs. However, the goal is to propose a representation of the actual model on a low-dimensional manifold. If the reduced dimension is small enough, the increase in computational cost is usually acceptable. In addition, the procedure here is entirely data-driven; no additional model evaluations are required to find the projections and the PCE coefficients.

\section{Applications}
\label{sec:application}
\subsection{Borehole model}
The PPA method is first applied on a borehole model (the same model as in \cite{zeng2021accelerated}) to make a comparison with the classical basis adaptation and the accelerated basis adaptation,
\begin{equation}
	\label{eq:borehole}
	f(\bm x) = \frac{2\pi T_u(H_u-H_l)}{\text{ln}(r/r_w)\left(1+\frac{2LT_u}{\text{ln}(r/r_w)r_w^2K_w}+\frac{T_u}{T_l}\right)}\,,
\end{equation}
which models the water flow through a borehole drilled from the ground surface through two aquifers. In this model we assume the radius of the borehole is $r_w = 0.1$ (m). The other parameters are randomized and listed in Table \ref{tab:borehole}.
\begin{table}[htb]
	\setstretch{1.25}
	\centering
	\caption{Random parameters and distributions of the borehole model}
	\label{tab:borehole}
	\begin{tabular}{l||ll}
		\hline \hline
		Variable & Distribution & Description \\
		\hline
		$r$ & Lognormal($\mu=7.71, \sigma=1.0056$) & Radius of influence (m) \\
		$T_u$ & Uniform(63070, 115600) & Transmissivity of upper aquifer (m$^2$/yr) \\
		$H_u$ & Uniform(990, 1110) & potentiometric head of upper aquifer (m) \\
		$T_l$ & Uniform(63.1, 116) & transmissivity of lower aquifer (m$^2$/yr) \\
		$H_l$ & Uniform(700, 820) & potentiometric head of lower aquifer (m)\\
		$L$ & Uniform(1120, 1680) & length of borehole (m) \\
		$K_w$ & Uniform(9855, 12045) & hydraulic conductivity of borehole (m/yr) \\
		\hline \hline
	\end{tabular}
\end{table}
The QoI $f(\bm{x})$ is the water flow rate of the borehole (m$^3$/yr).

A previous UQ exploration \cite{zeng2021accelerated} of the borehole model showed that an accurate probability density function (PDF) of the water flow rate could be obtained through (i) MCS of 3000 samples; (ii) classical basis adaptation of dimension 5 that requires 995 samples evaluated at specific Gauss quadrature points; (iii) accelerated basis adaptation of dimension 3 that required 308 samples evaluated on quadrature points. The maximum PCE order is $p=3$ in the last two methods, which is verified to be adequate.

\begin{figure}[htb]
	\centering
	\begin{minipage}{0.45\linewidth}
		\centering
		\includegraphics[width=\linewidth]{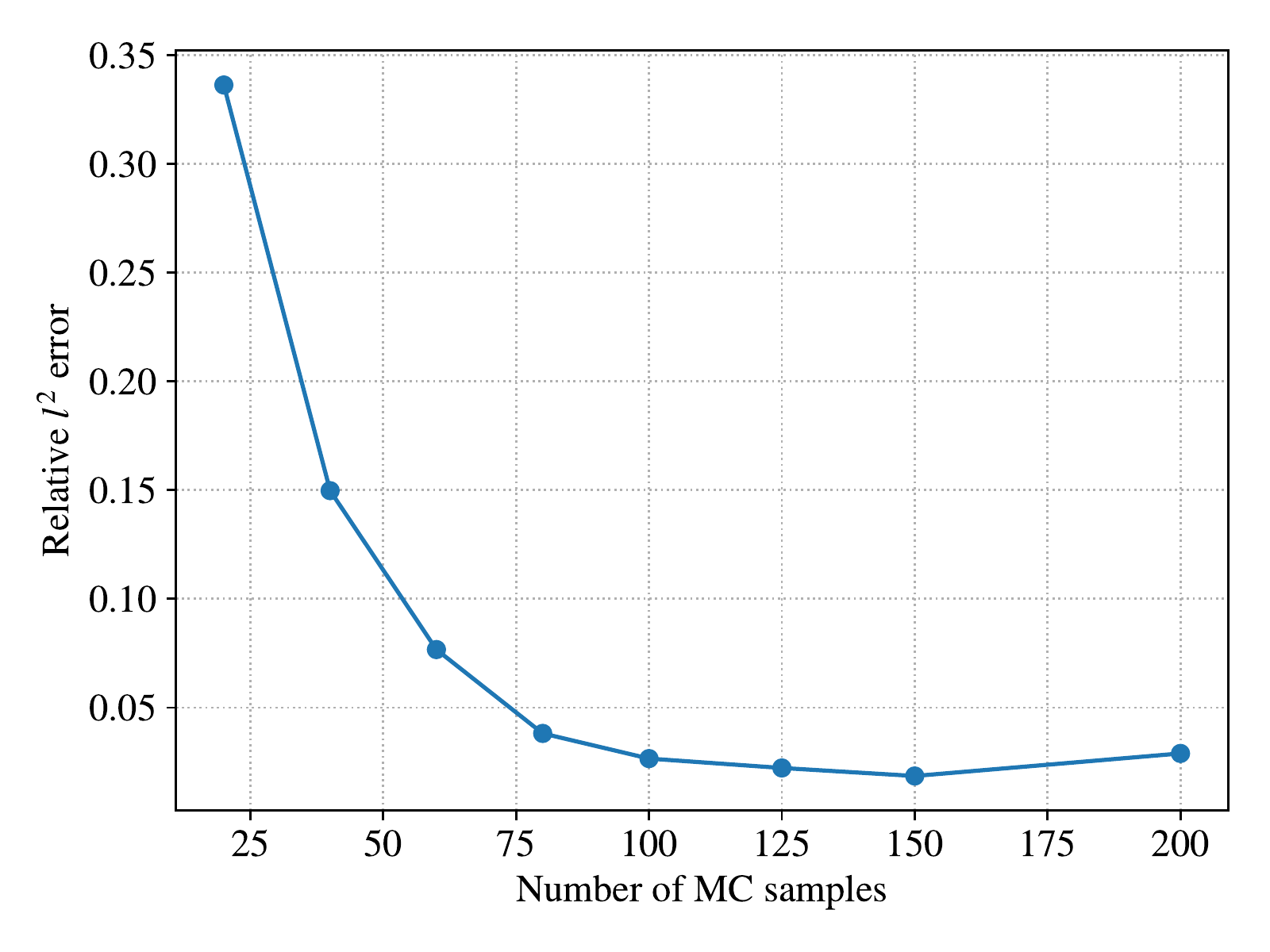}
		\caption{The convergence curve of the relative $l^2$ error of the PDF obtained from PPA with different number of MC samples for the borehole model.}
		\label{fig:rela_l2_err_bore}
	\end{minipage}
	\hspace{0.2cm}
	\begin{minipage}{0.45\linewidth}
		\centering
		\includegraphics[width=\linewidth]{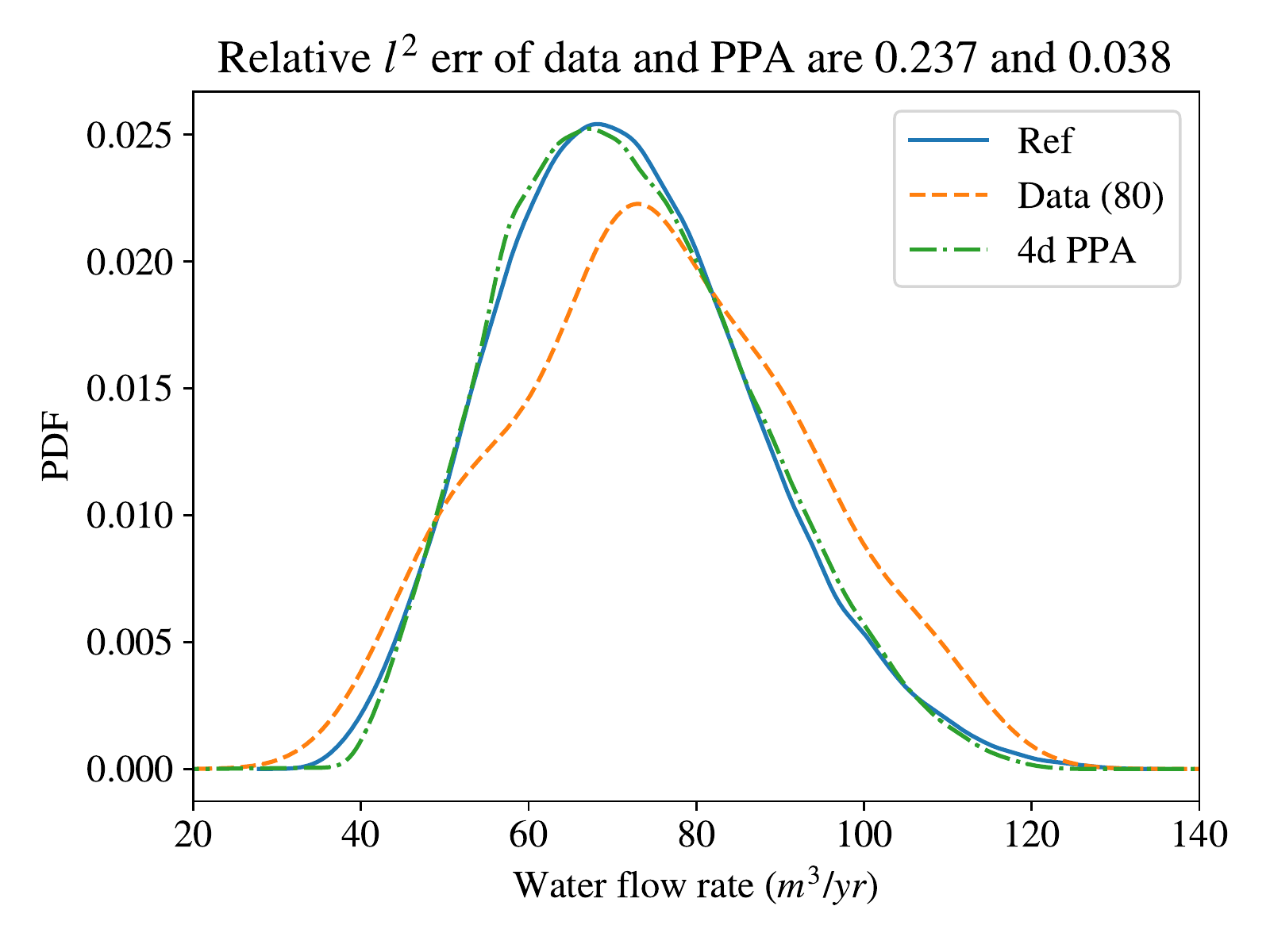}
		\caption{PDF comparison of the PPA method for the borehole model. Three PDFs are compared: the reference from 100,000 MC samples, the data, and the converged PPA model.}
		\label{fig:pdf_compare_bore_80}
	\end{minipage}
\end{figure}

We now apply the data-driven PPA method to learn the adapted projections and the PCE model with respect to the projected variables with a different number of MC samples $\{20, 40, 60, 80, 100, 125, 150, 200\}$. Using the same maximum PCE order $p=3$, we learn the converged optimally adapted projections and the corresponding optimal PCE models. Then we can generate many MC samples based on the resulting PCE models and estimate the PDF of the QoI by the KDE. A reference PDF is generated from 100,000 MC simulations of the borehole model by KDE for comparison. For any two vectors $\bm{u}$ and $\bm{v}$, we define the relative $l^2$ difference of $\bm{u}$ with respect to $\bm{v}$ as the following
\begin{equation}
	e = \frac{\lVert \bm{u} - \bm{v} \lVert_2}{\lVert \bm{v}\lVert_2}
\end{equation}
where $\lVert \cdot \lVert_2$ represents the $l^2$ norm. Then, compared to the reference PDF, we can compute the relative $l^2$ error of the PDFs generated from PPA with a different number of MC samples. The results are shown in Figure \ref{fig:rela_l2_err_bore}.
We see that the relative $l^2$ errors are smaller than 5\% when the MC samples are greater than 80 and converge with 150 MC samples. Therefore, we show the PDF comparison when 80 MC samples are used to learn the PPA model in Figure \ref{fig:pdf_compare_bore_80}. In the figure, three PDFs are compared, the reference PDF, the PDF generated by the data, and the PDF by the PPA model. Clearly, the PPA PDF is close to the reference, with a relative $l^2$ error of 3.8\%. In comparison, the PDF generated directly by the data has a large discrepancy from the reference, with an $l^2$ error of 23.7\%. 

While for evaluating lower order statistics, the accuracy of this PPA model using 80 MC samples is adequate, other evaluations, such as reliability and failure probabilities, could require more samples.  It is clear from Figure (\ref{fig:rela_l2_err_bore}) that 150 MC samples are sufficient for $l^2$ convergence, from which convergence in distribution, relevant for tail analysis, follows. 

\begin{figure}[htb]
	\centering
	\begin{minipage}{0.45\linewidth}
		\centering
		\includegraphics[width=\linewidth]{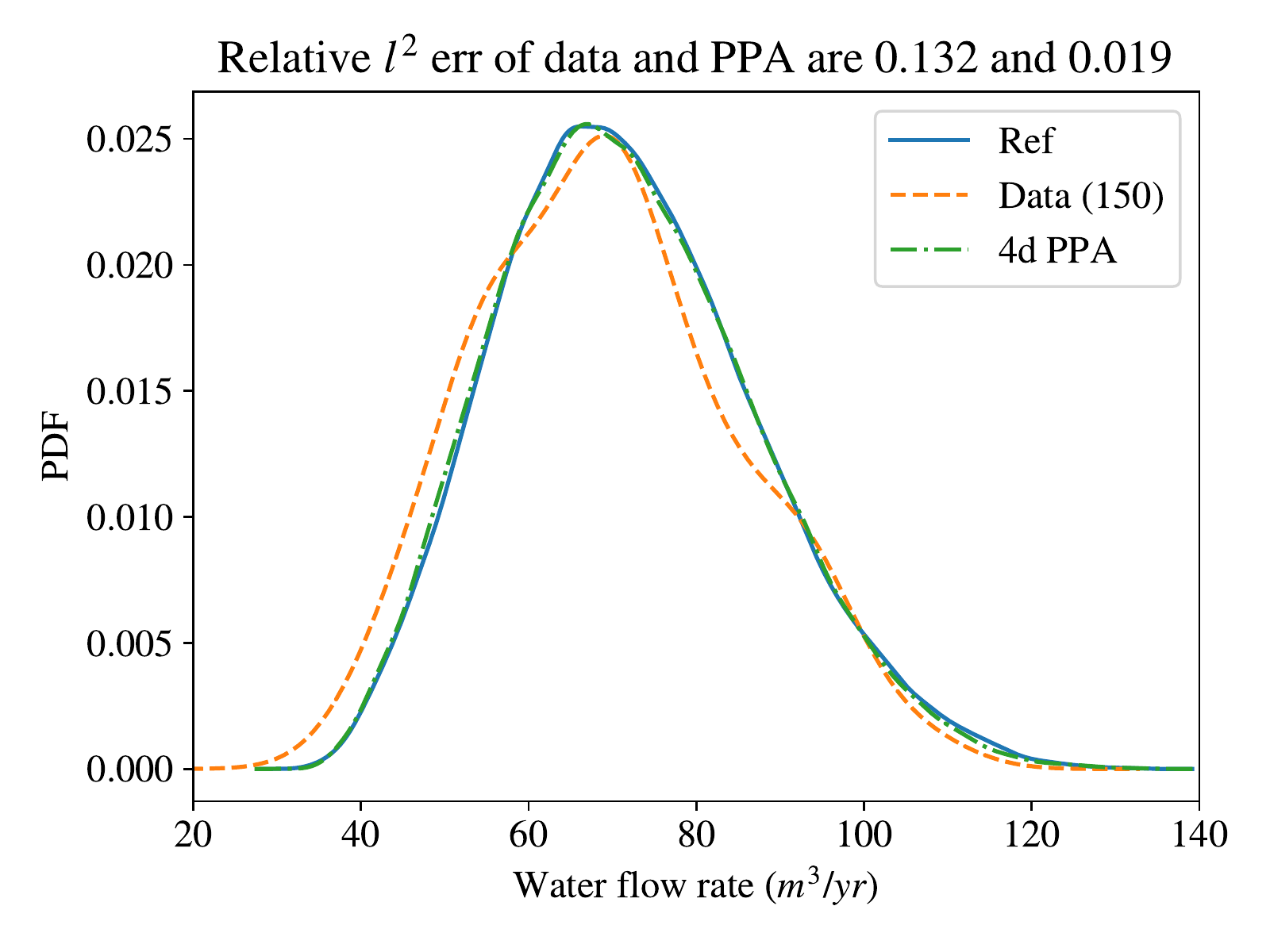}
		\caption{PDF comparison of the PPA method for the borehole model. Three PDFs are compared: the reference from 100,000 MC samples, the data, and the converged PPA model.}
		\label{fig:pdf_compare_bore}
	\end{minipage}
	\hspace{0.1cm}
	\begin{minipage}{0.45\linewidth}
		\centering
		\includegraphics[width=\linewidth]{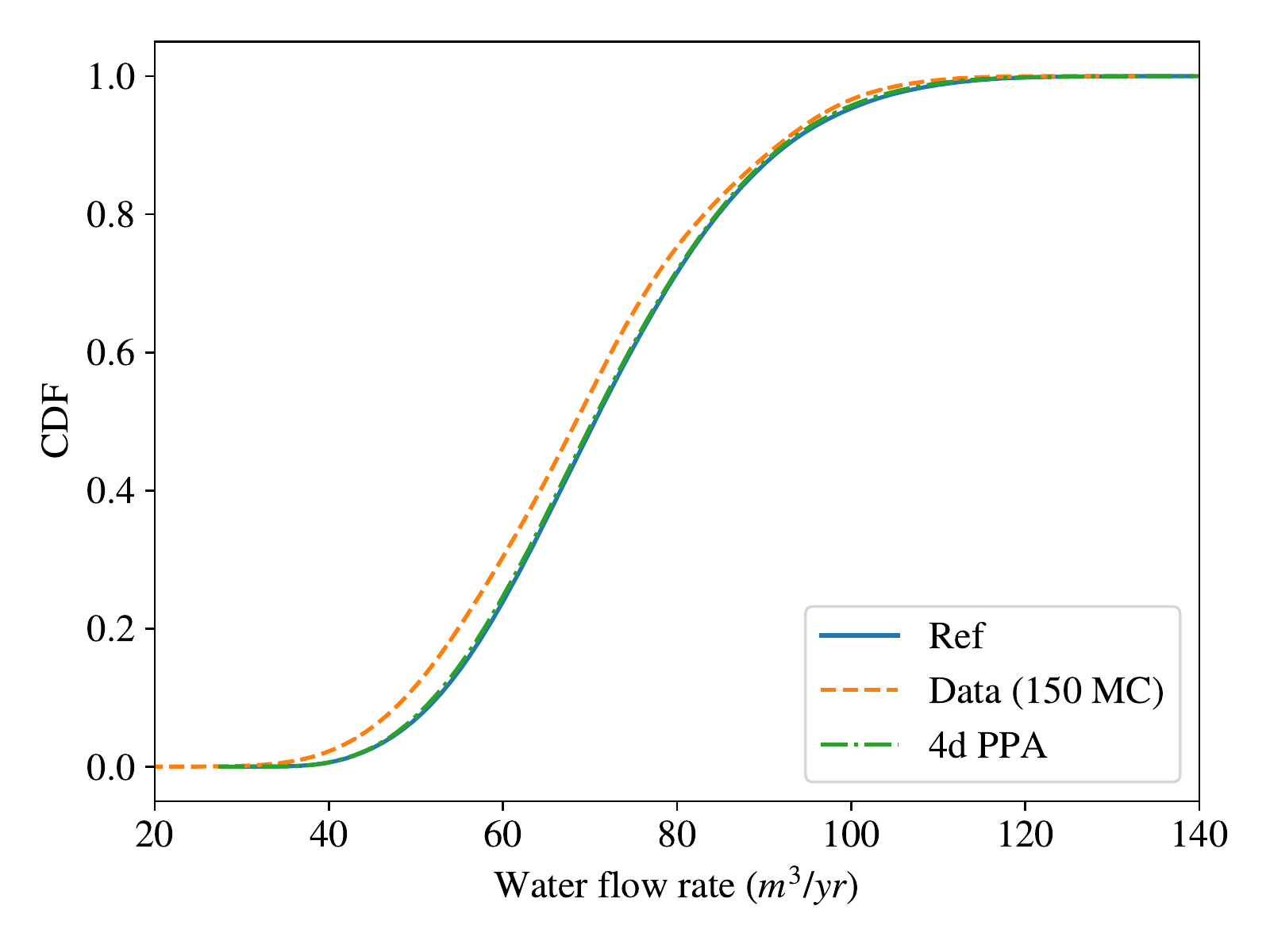}
		\caption{CDF comparison of the PPA method for the borehole model. Three PDFs are compared: the reference from 100,000 MC samples, the data, and the converged PPA model.}
		\label{fig:cdf_compare_bore}
	\end{minipage}
\end{figure}

To check the accuracy of the PPA model, Figures \ref{fig:pdf_compare_bore} and \ref{fig:cdf_compare_bore} present a comparison of PDFs and CDFs obtained from the reference model, the 150 data, and the PPA model generated from the 150 data. 
The first observation is that the resulting PDF of the PPA method is even closer to the reference compared to the case of using 80 MC samples. The relative $l^2$ difference between the PPA model and the reference decreases to 1.9\%. The accuracy improvement in the left tail region is noticeable. The second observation is that the PDF generated directly from the data differs from the reference, especially in the left tail region. From the CDF figure, the curve associated with the data is far from the reference, while the result from PPA is consistent with the reference.

\begin{figure}[htb]
	\centering
	\begin{minipage}{0.45\linewidth}
		\centering
		\includegraphics[width=\linewidth]{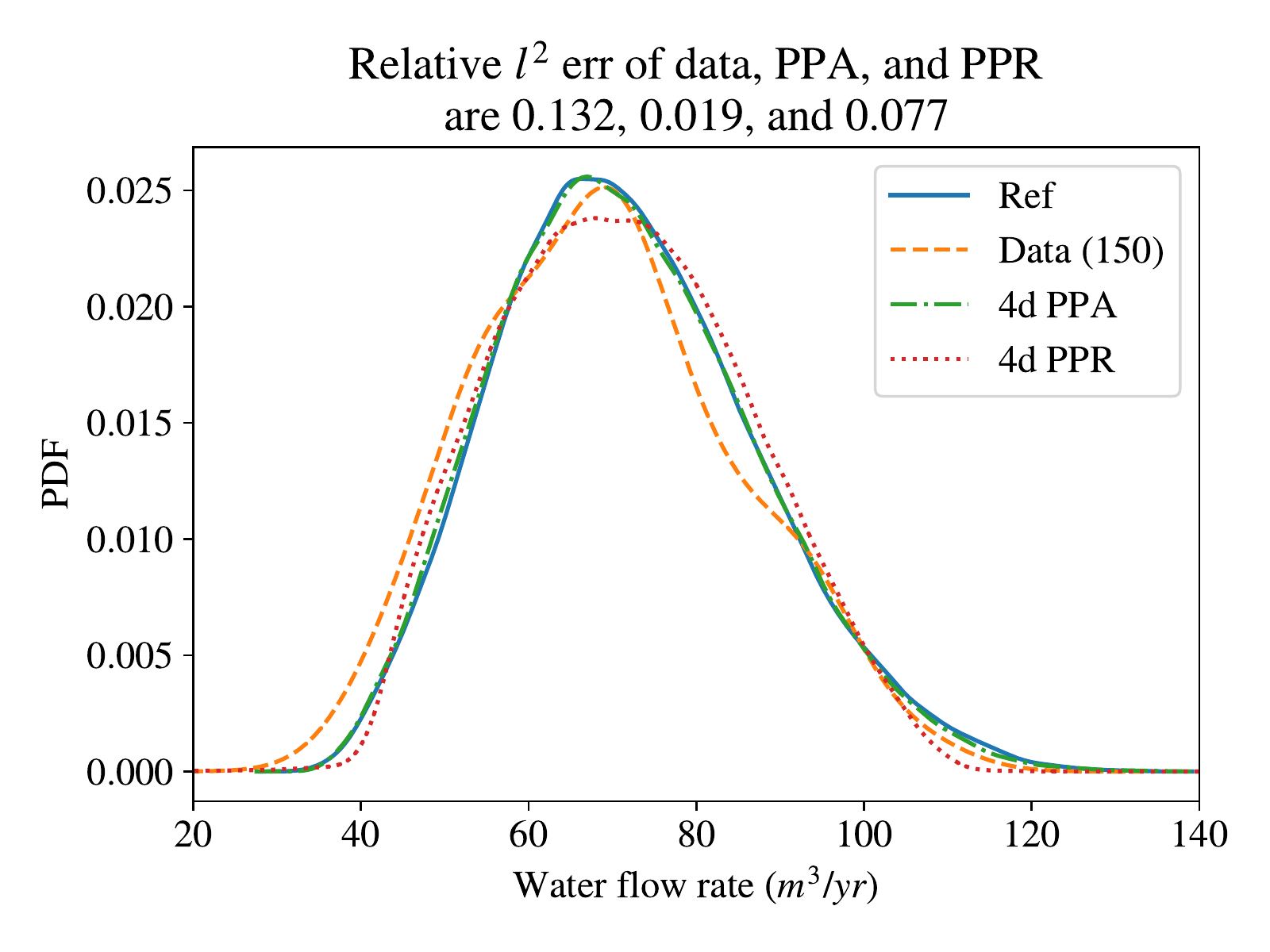}
		\caption{PDF comparison of the PPA method for the borehole model. Three PDFs are compared: the reference from 100,000 MC samples, the data, the converged PPA model, and the converged PPR model.}
		\label{fig:pdf_compare2_bore}
	\end{minipage}
	\hspace{0.1cm}
	\begin{minipage}{0.45\linewidth}
		\centering
		\includegraphics[width=\linewidth]{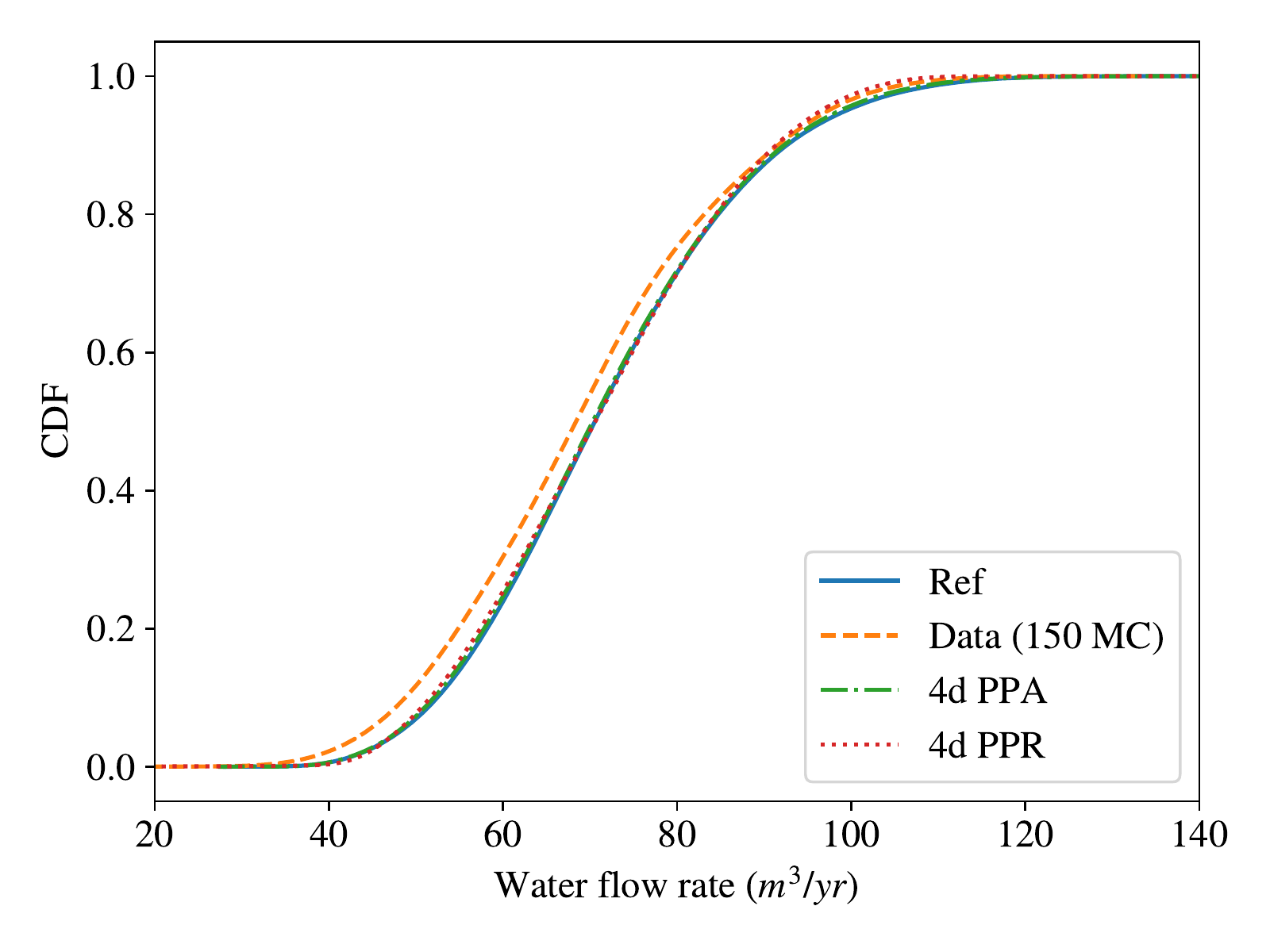}
		\caption{CDF comparison of the PPA method for the borehole model. Three CDFs are compared: the reference from 100,000 MC samples, the data, the converged PPA model, and the converged PPR model.}
		\label{fig:cdf_compare2_bore}
	\end{minipage}
\end{figure}

Since the PPA method is developed based on the PPR method, a comparison of these two techniques is conducted. PCEs are used as smooth functions in PPR, leading to an additive model of many univariate PCEs. We apply the PPR method on the 150 MC samples to obtain a converged model with the same stopping criterion as the PPA method. Then, the PPR model is compared to the PPA model. The resulting PDFs and CDFs are shown in Figures \ref{fig:pdf_compare2_bore} and \ref{fig:cdf_compare2_bore}.
We see that the PPR model also has improved over the data and is consistent with the reference overall. However, it misses the left and right tail regions, and the PDF peak at around 70 m$^3$/yr water flow rate. The relative $l^2$ error of the PPR model is 7.7\% compared to 1.9\% of the PPA model. The tail regions usually correspond to reliability or failure probability and thus are important in engineering. The CDF comparison presents no clear separation between PPA and PPR. However, we can still see a difference if we examine the right tail. For example, $P(f(\bm{x}) < 97) = 0.95$ for PPR, while  $P(f(\bm{x}) < 99) = 0.95$ for PPA, suggesting that these two models could make different conclusions regarding the reliability and failure probability.

Compared to MCS, where 3000 samples are required for a converged PDF, we use only 150 samples in PPA to construct a model with the same accuracy. The massive computation saving is because (i) the QoI can be represented on a low-dimensional manifold; (ii) the PPA method can discover the low-dimensional manifold and the function that maps the parameter space to the output space. In this application, the PPA reduced the dimension from 7 to 4, and the number of PCE terms of order 3 has reduced from 120 to 35. Moreover, the number of MC samples required to compute the PCE coefficients accurately will reduce significantly. However, even if we use the accelerated basis adaptation method where the dimension could be reduced to 3, the number of sparse quadrature points is 308, which is still two times greater than 150. In addition, the model is evaluated on specific quadrature points in accelerated basis adaptation to save cost, which in many applications might not be available.

\begin{figure}[htb]
	\centering
	\begin{minipage}{0.45\linewidth}
		\centering
		\includegraphics[width=\linewidth]{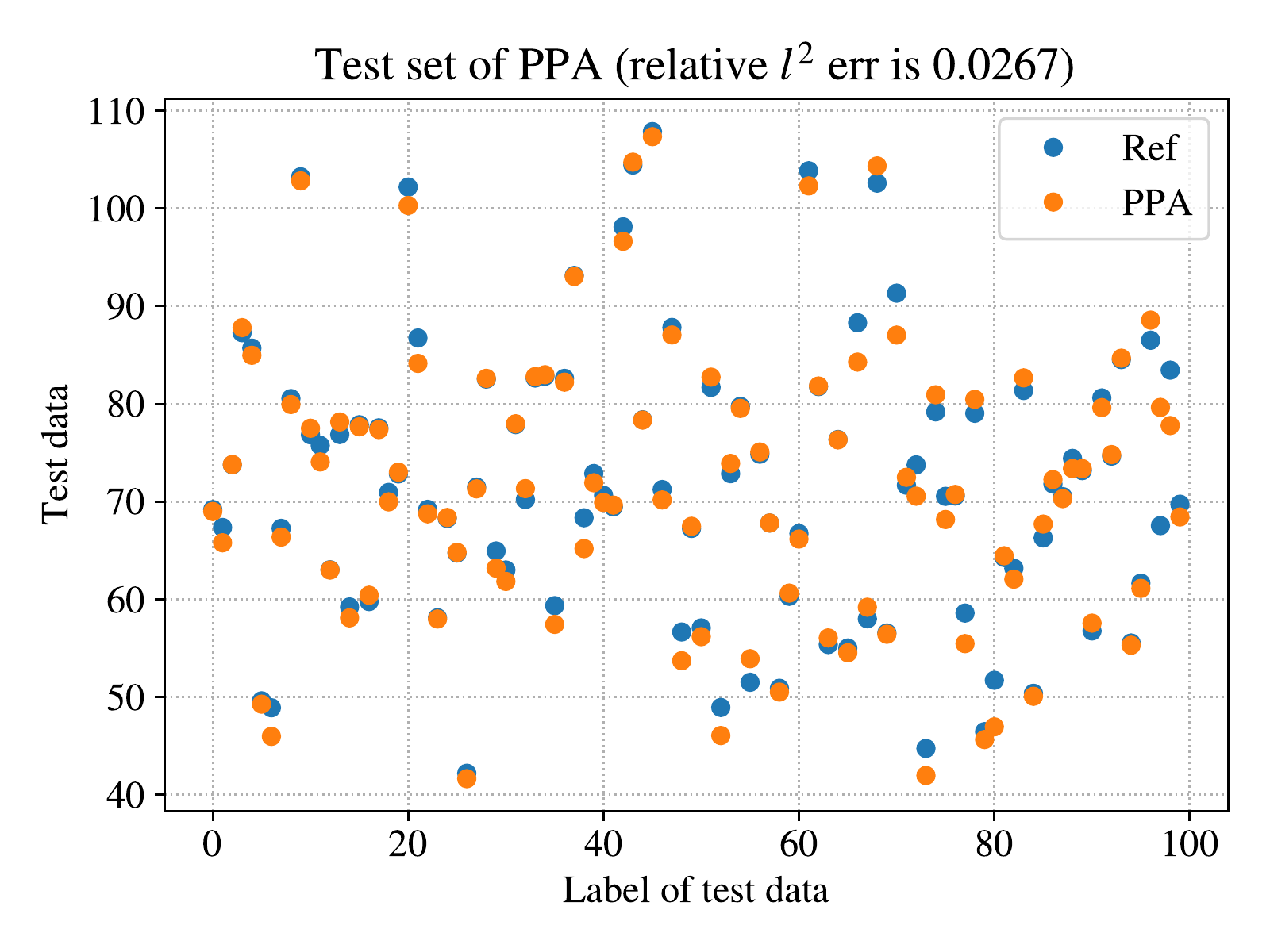}
	\end{minipage}
	\hspace{0.1cm}
	\begin{minipage}{0.45\linewidth}
		\centering
		\includegraphics[width=\linewidth]{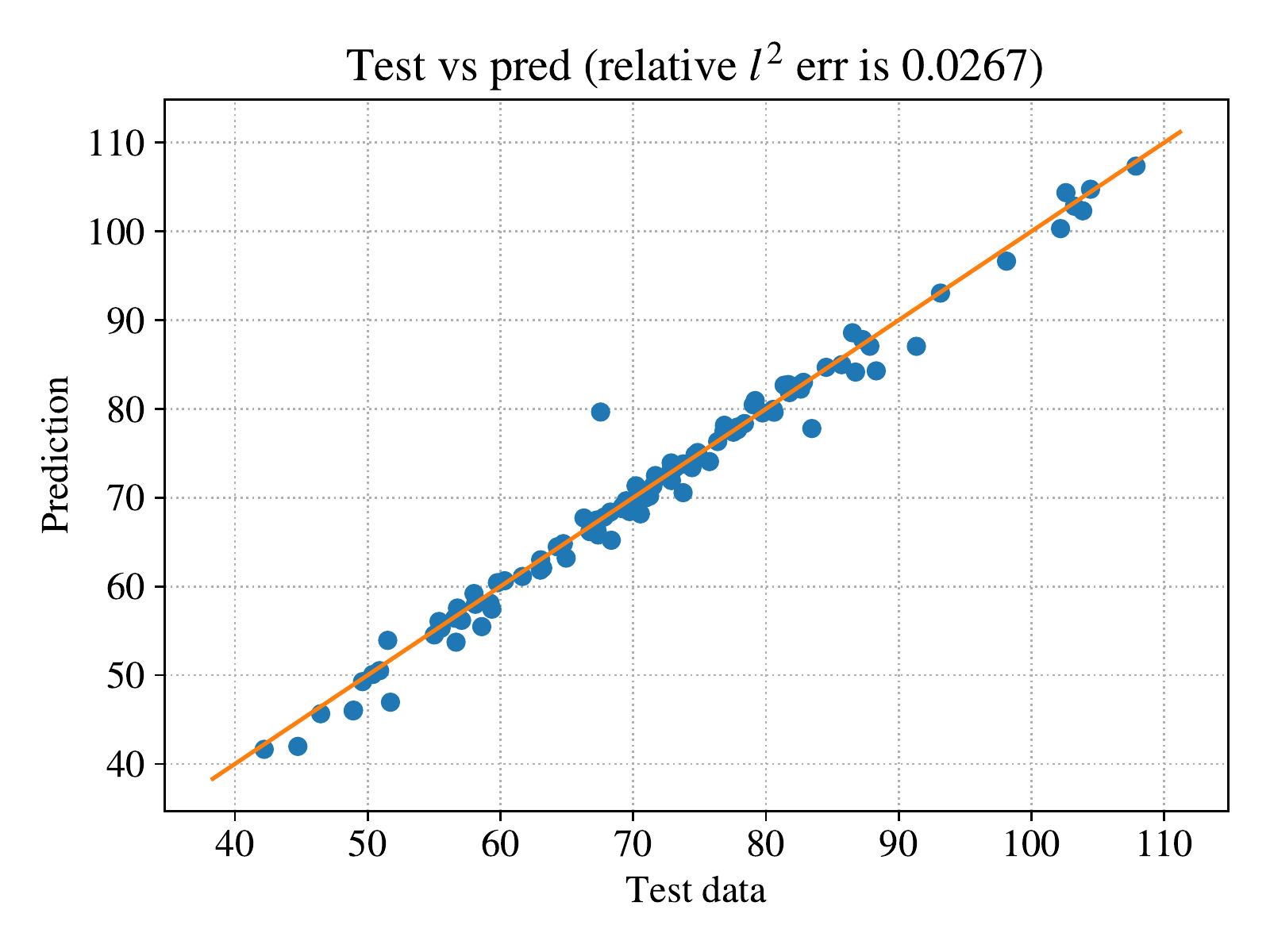}
	\end{minipage}
	\caption{Test data (100) of the PPA method for the borehole model. Figure on left shows the QoI value comparison of reference and prediction on various test data. Figure on the right plot the prediction against the test data where the solid line is of slope 1. }
	\label{fig:test_bore}
\end{figure}

The PPA model can not only be used for UQ but can also be used in model prediction. In the latter, the PPA model can serve as a surrogate model to predict QoI with different input parameters. The PPA model essentially constructed a nonlinear mapping from input space to output space based on the given data. For any given input parameter, the prediction of the associated output value can be obtained by querying from the PPA model. We call the data used to construct the PPA model training data, which in this case, is the 150 MC samples. To check the accuracy of the PPA model in prediction, we generate 100 MC samples that are different from the training data. We can predict the QoI on these 100 samples using the PPA model. The comparison of the prediction and the reference model on the test set is shown in Figure \ref{fig:test_bore}.
From the left figure, the PPA prediction (in orange) is close to the reference (in blue) on the majority of the test data. From the right figure, the scatter dots of prediction against the test set are closely aligned in the unit slope line, meaning that the predictions are close to the reference. The relative $l^2$ difference of the 100 test data is 2.0\%, which is close to the relative $l^2$ error of the PPA PDF compared to the reference PDF. The results suggest that the prediction is accurate.

\begin{figure}[htb]
	\centering
	\begin{minipage}{0.45\linewidth}
		\centering
		\includegraphics[width=\linewidth]{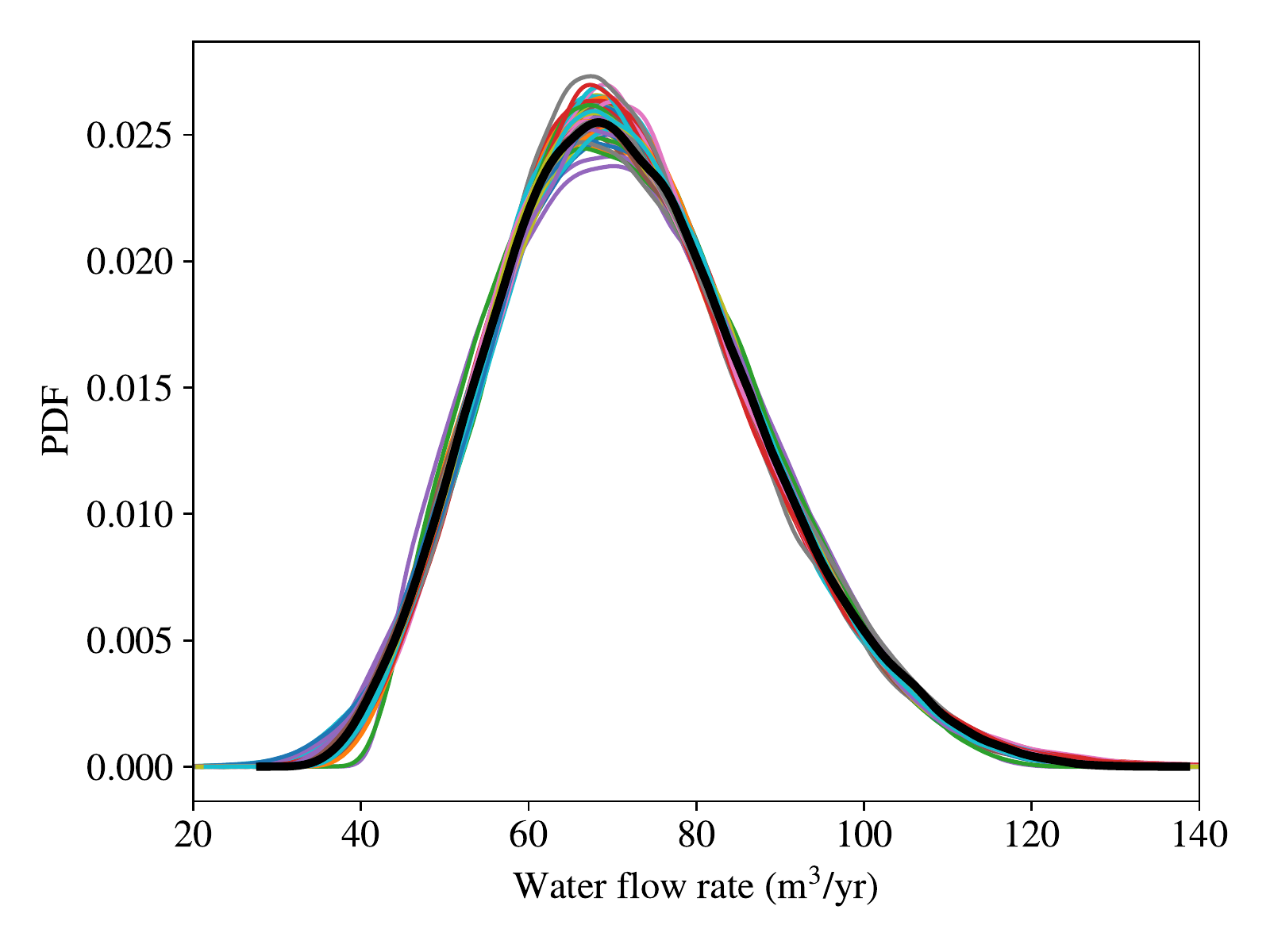}
		\caption{Plot of 100 PDF batch generated by PPA method with different 150 MC samples. The thick black curve represent the reference PDF.}
		\label{fig:bore_pdf_batch}
	\end{minipage}
	\hspace{0.1cm}
	\begin{minipage}{0.45\linewidth}
		\centering
		\includegraphics[width=\linewidth]{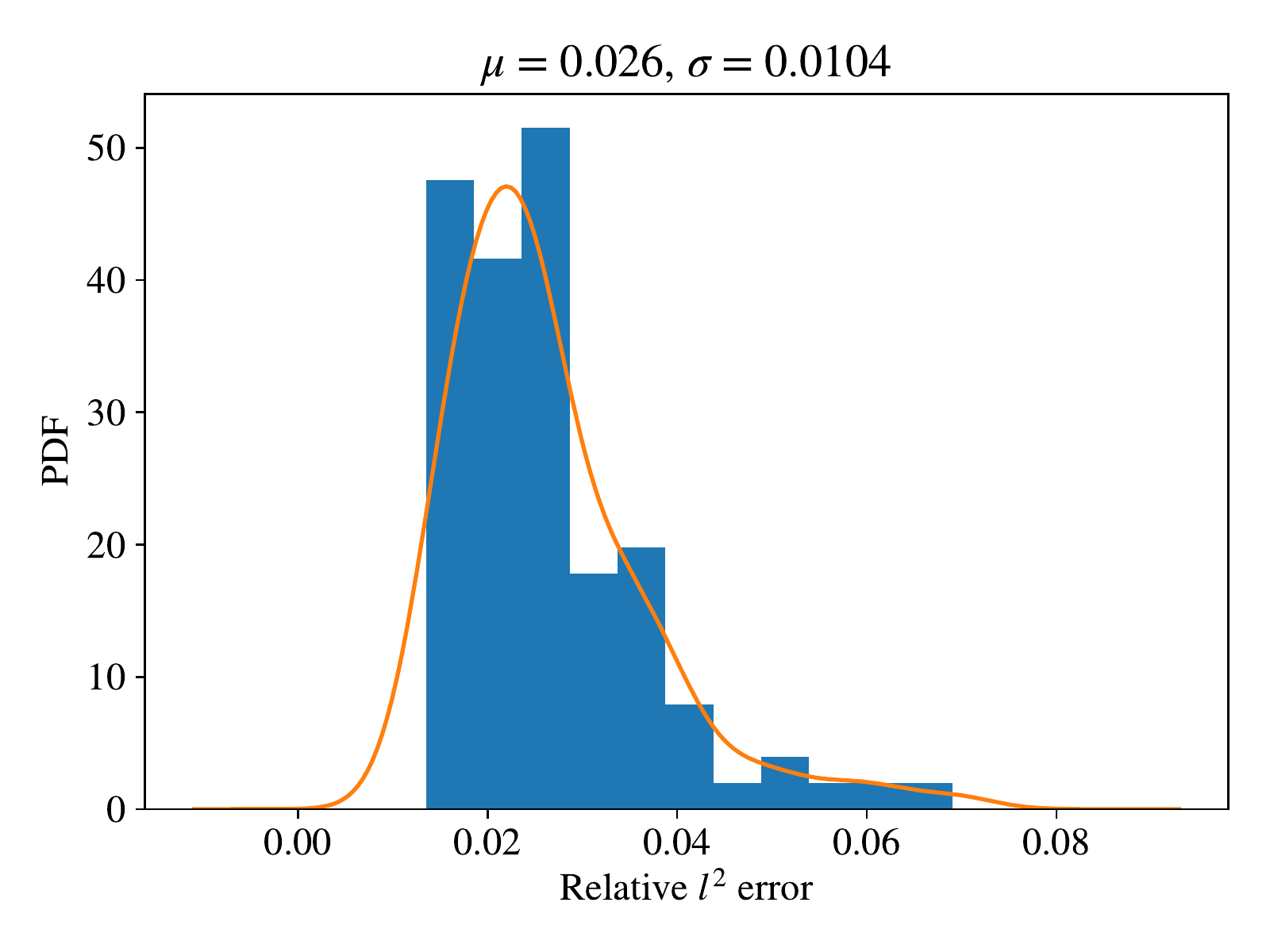}
		\caption{The histogram of relative $l^2$ errors of the 100 PDFs from Figure \ref{fig:bore_pdf_batch} compared to the reference PDF.}
		\label{fig:bore_l2_pdf}
	\end{minipage}
\end{figure}

To check the robustness of the PPA method, we generated 100 sets of MC samples, each containing 150 samples. Within each set, we apply the PPA method to find a converged PCE model and generate many MC samples and compute the associated KDE estimation of the PDF. The PDFs from these sets are plotted in Figure \ref{fig:bore_pdf_batch}, where the reference PDF in thick black is also presented.
Although generated with different MC samples, the PDFs are all close to the reference. To quantify the difference, we computed the relative $l^2$ errors of these 100 PDFs with respect to the reference. The histogram of the relative $l^2$ errors is shown in Figure \ref{fig:bore_l2_pdf}. The mean value and standard deviation of the relative $l^2$ errors are 0.026 and 0.0092, respectively. All the relative $l^2$ errors are less than 5\%, and most of the $l^2$ errors are less than 3\%. This means that the variation in the accuracy of the PDF estimation is small when different MC samples are used in the PPA method.

% \todo{Comparison of directions from adaptation and PPA}

\subsection{Structural Dynamics}
\begin{figure}[htb]
	\centering
	\begin{minipage}{0.25\linewidth}
		\centering
		\includegraphics[width=\linewidth]{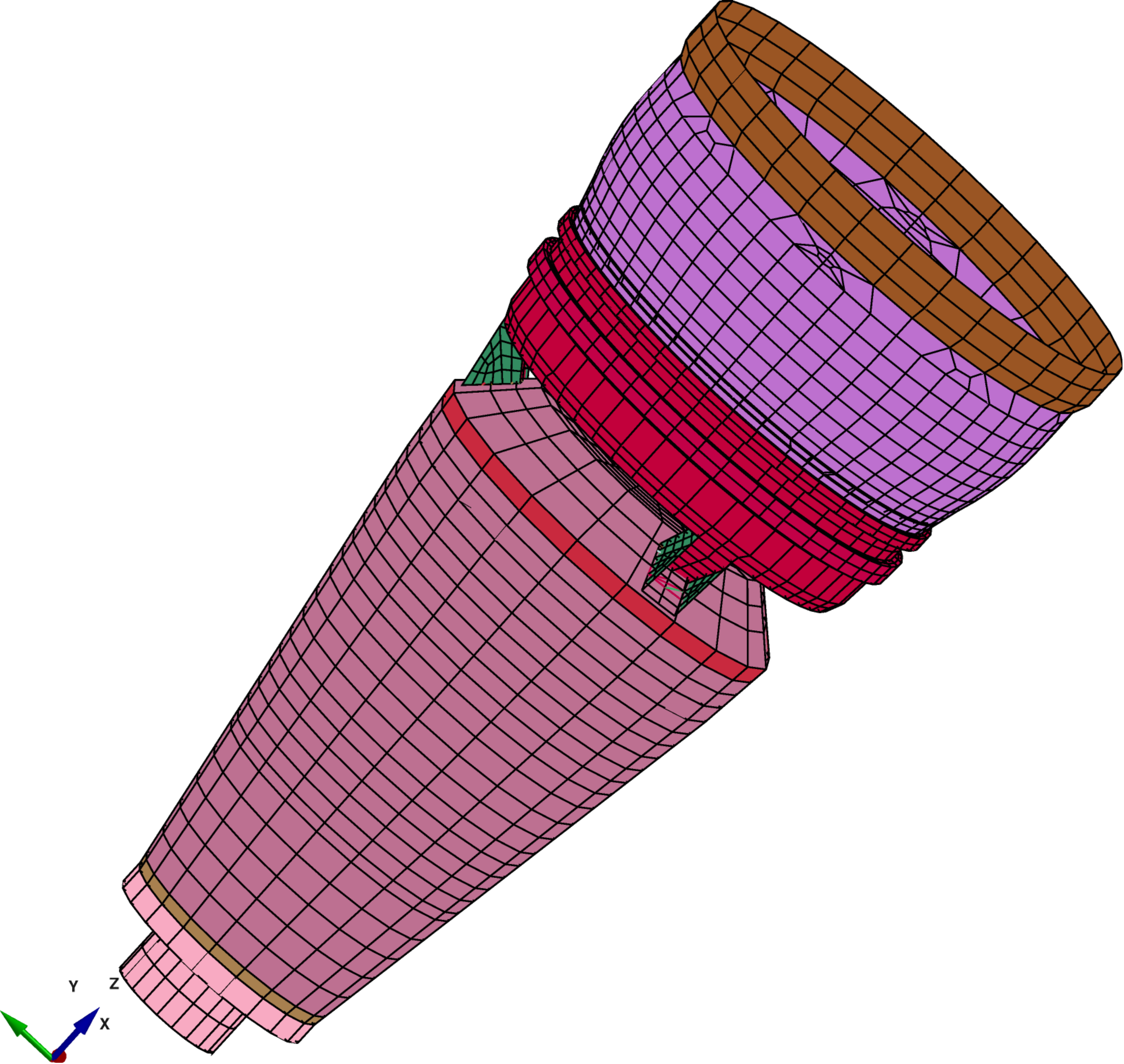}
	\end{minipage}
	\hspace{0.2cm}
	\begin{minipage}{0.25\linewidth}
		\centering
		\includegraphics[width=\linewidth]{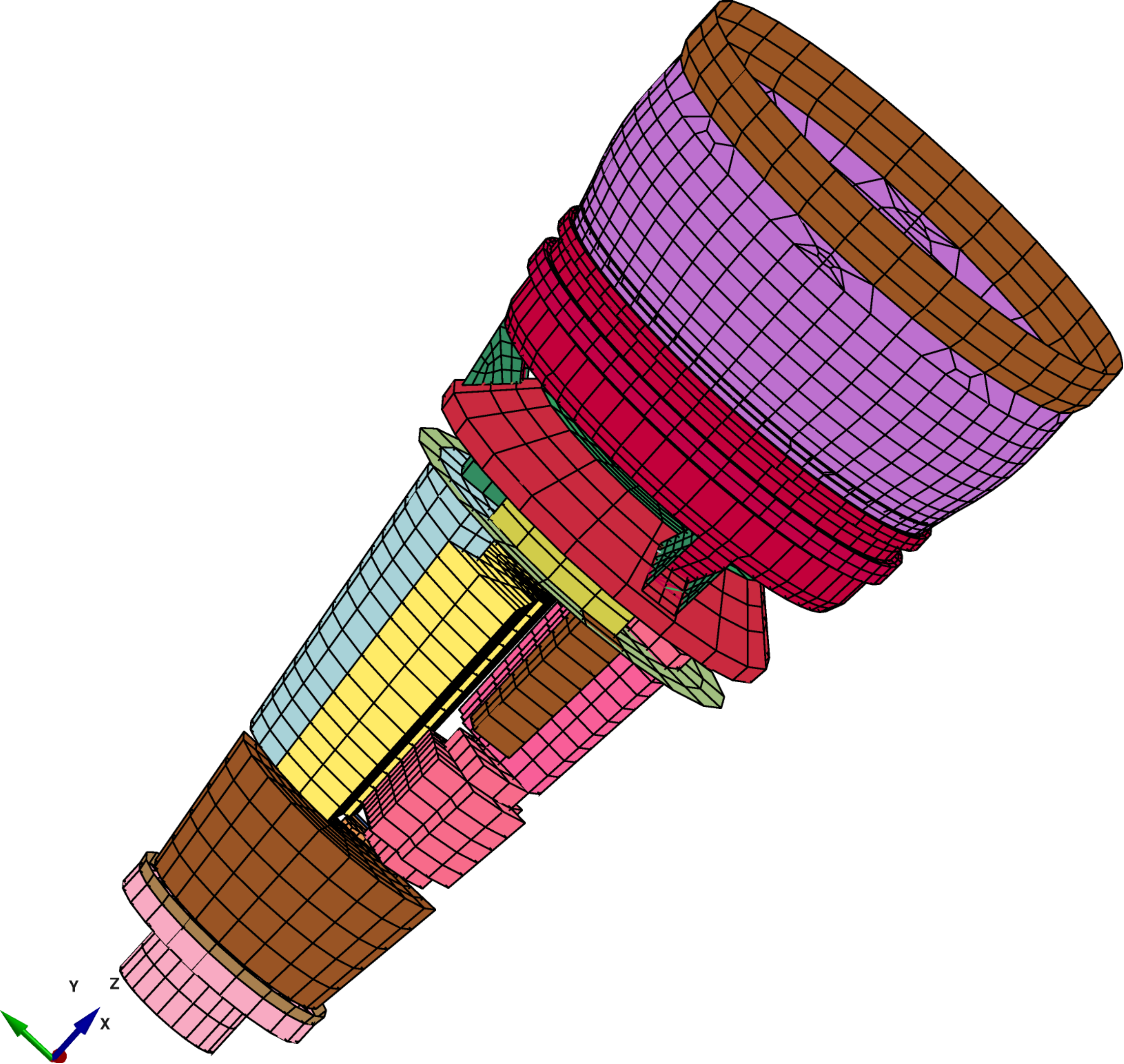}
	\end{minipage}
	\hspace{0.2cm}
	\begin{minipage}{0.4\linewidth}
		\centering
		\includegraphics[width=\linewidth]{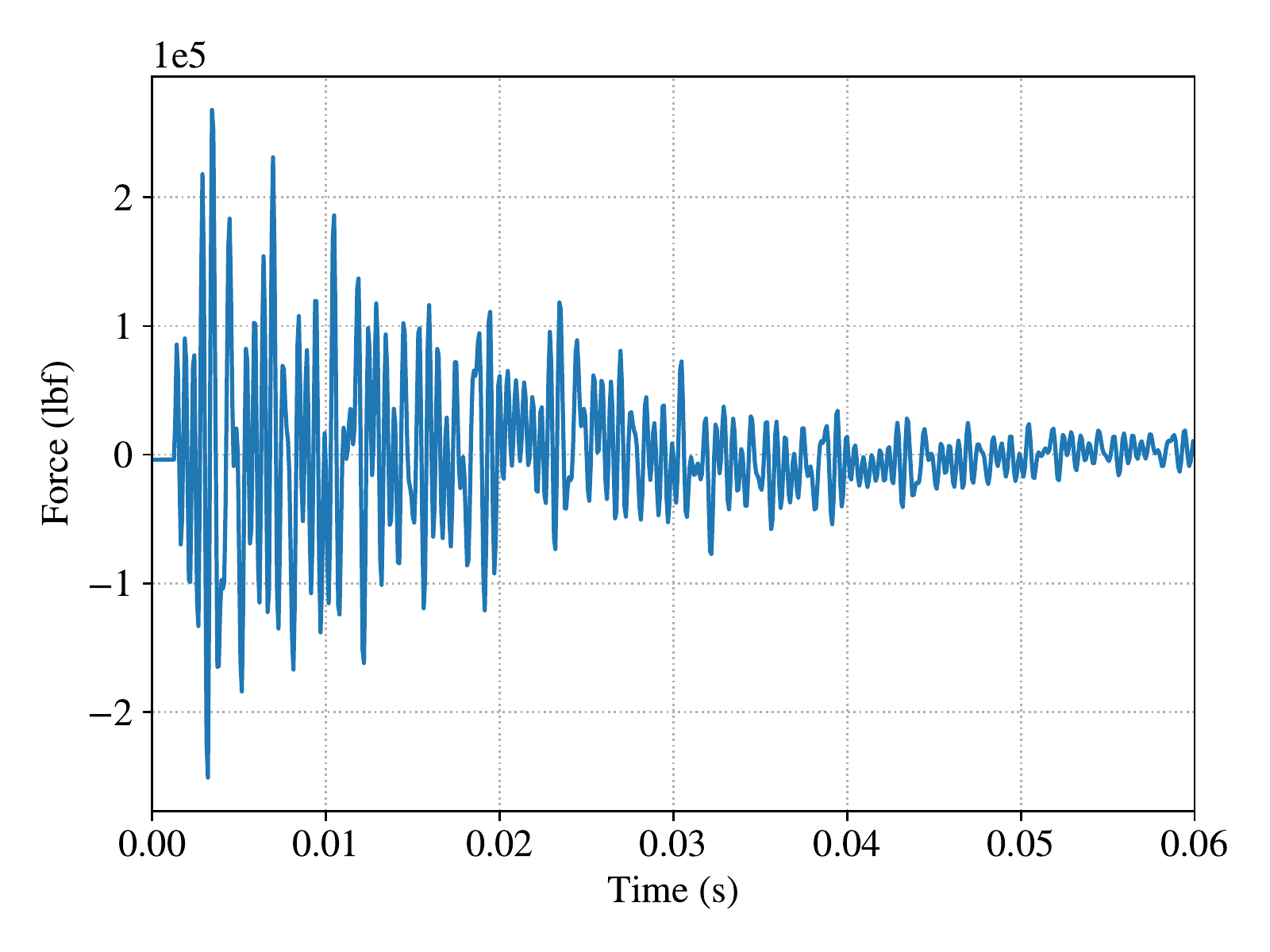}
	\end{minipage}
	\caption{The overall finite element model of the space structure (left); the finite element model of the space structure without the outer shell and shock absorption block (middle); time history of the impulse load (right).}
	\label{fig:space_model}
\end{figure}

\begin{figure}[htb]
	\centering
	\begin{minipage}{0.4\linewidth}
		\centering
		\includegraphics[width=\linewidth]{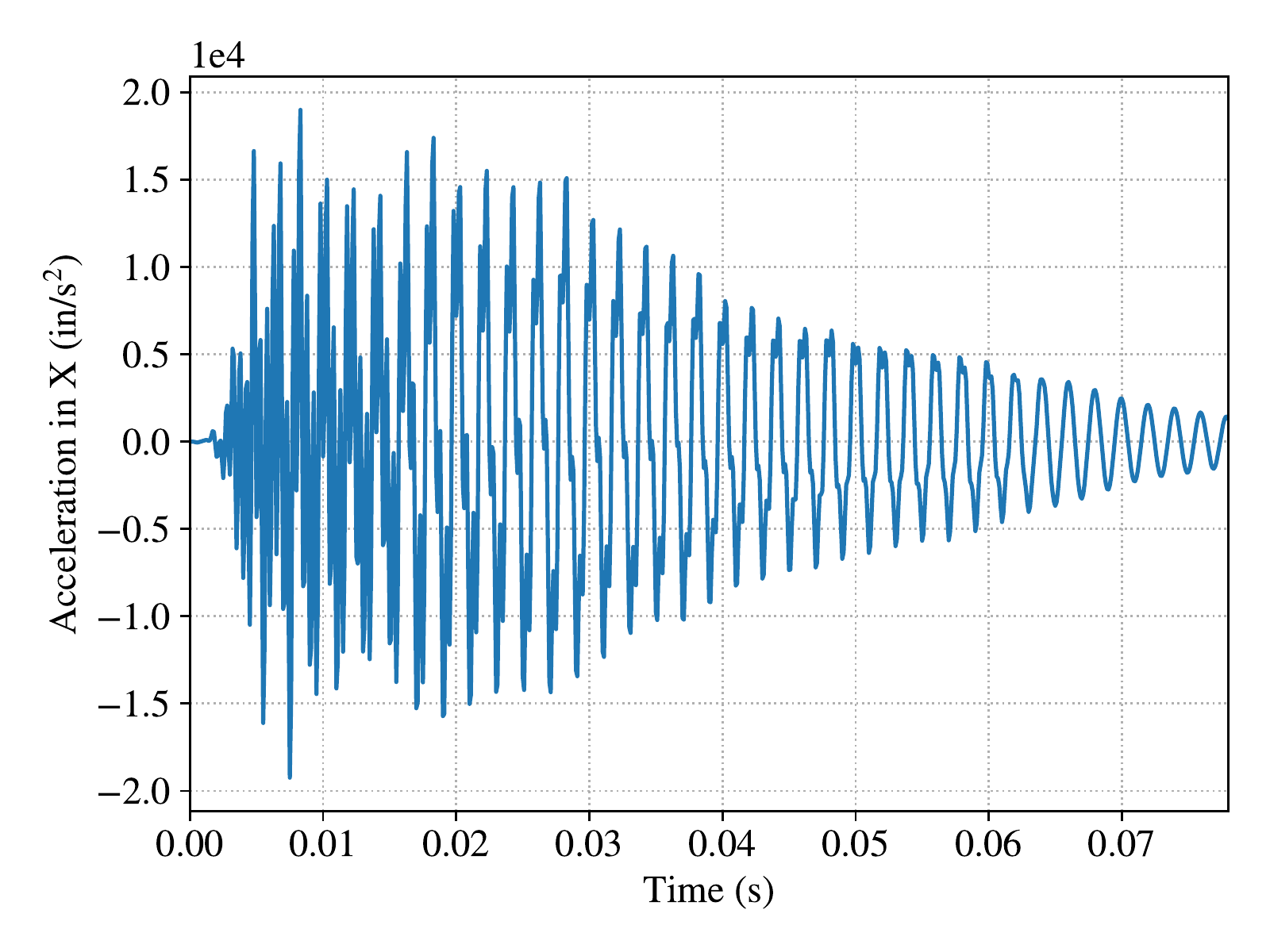}
	\end{minipage}
	\hspace{0.5cm}
	\begin{minipage}{0.4\linewidth}
		\centering
		\includegraphics[width=\linewidth]{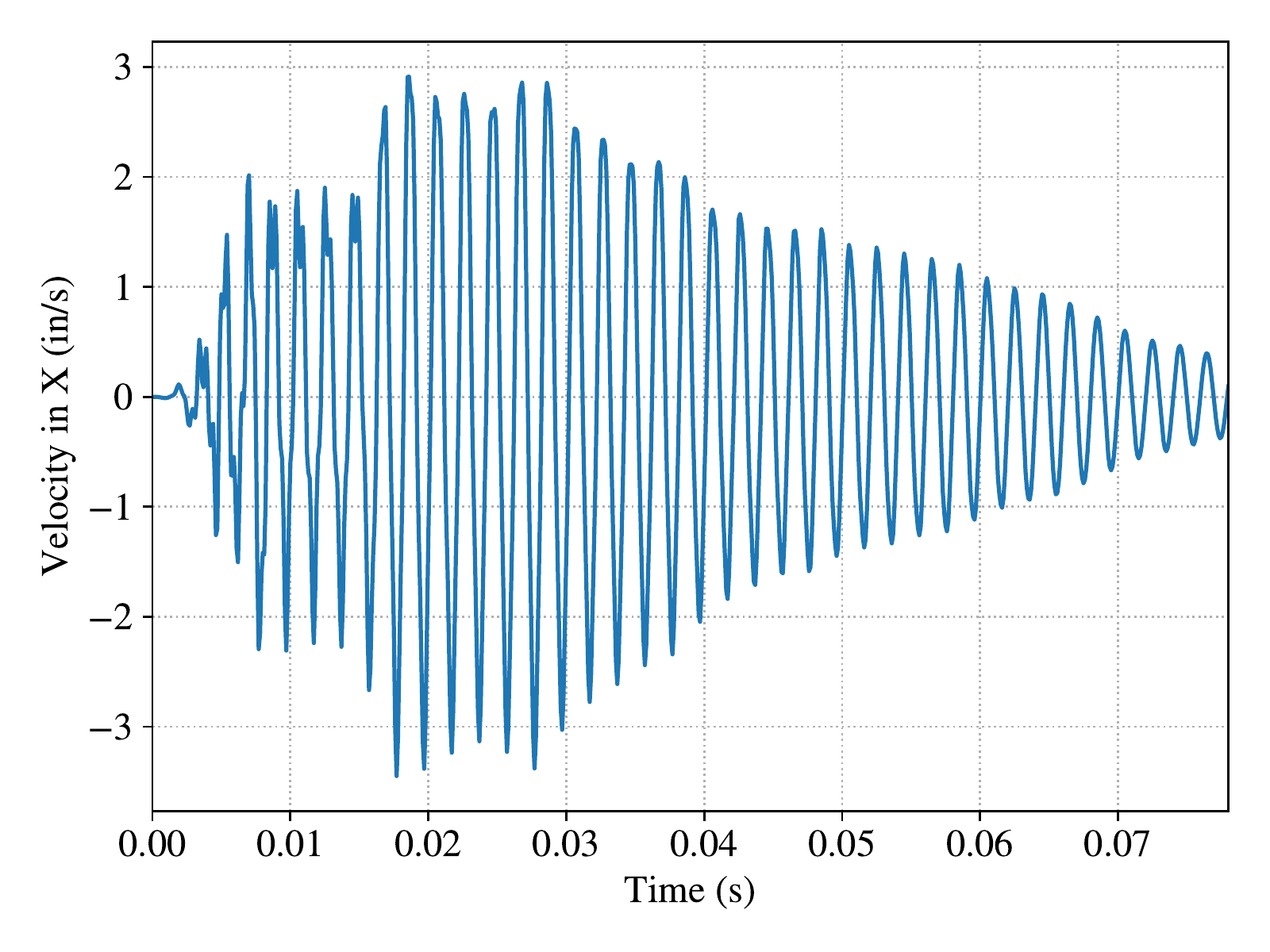}
	\end{minipage}
	\caption{Time history of acceleration (left) and velocity (right) of the essential part under the impulse loading shown in Figure \ref{fig:space_model}.}
	\label{fig:space_acc_vel_hist}
\end{figure}

The second application considers a space structure subjected to impulse load in the $Z$ direction. The computational model is shown in Figure \ref{fig:space_model}. 
The model can be separated into upper and lower parts connected by a three-point mounting pedestal. The upper part is open and is subjected to impulse load with the time history shown in Figure \ref{fig:space_model}. The lower part contains an outer shell, a solid shock absorption block, and essential components that can be seen in the middle figure of Figure \ref{fig:space_model}. The effects of the impulse loading are transferred from the upper part to the lower parts through the three-point mounting pedestals. Therefore, they can affect the functionality of the essential components in the lower part. In this application, we are particularly interested in one essential part maximum acceleration and maximum velocity responses. The time history of the acceleration and velocity of the essential part is shown in Figure \ref{fig:space_acc_vel_hist}.

\subsubsection{Uncertainty quantification and prediction}
The maximum acceleration and velocity can be directly read from the response histories in a deterministic setting. However, the material properties of the model are almost inevitably subjected to variations, especially for some special materials. The randomness of the materials can affect the results a lot. Therefore, in this application, some material properties are randomized, see Table \ref{tab2}.
\begin{table}[htb]
	\setstretch{1.25}
	\centering
	\caption{Random parameters and distributions of the space model}
	\label{tab2}
	\begin{tabularx}{\textwidth}{X||cXc}
		\hline \hline
		Variable & Distribution & Mean  & Coefficient of variation \\
		\hline
		Modulus of elasticity of 3 major components of the upper parts & Lognormal &  $1.6\times 10^7$ psi & 10\%\\
		\hline
		Stiffness of the 3 point mounting pedestals & Lognormal & $5.7\times 10^7$ lbs/in & 20\% \\
		\hline
		Modulus of elasticity of 18 components in the lower part & Lognormal & {14 with 1.0×107 psi; 2 with 9.76×107 psi; 1 with 2.93×107 psi; 1 with 8829.0 psi} & 20\% \\
		\hline \hline
	\end{tabularx}
\end{table}
Note that the randomization here differs from the model in \cite{zeng2021accelerated}. In this space model, many structural components are modeled in two parts, a solid element part and a shell element part that wraps around the solid part. In this paper, we randomized the shell and solid elements, while in \cite{zeng2021accelerated}, only shell elements are randomized. However, a shared random germ is used for a single structure component, leading to a random dimension of 24, the same as in \cite{zeng2021accelerated}. 

\begin{figure}[htb]
	\centering
	\begin{minipage}{0.45\linewidth}
		\centering
		\includegraphics[width=\linewidth]{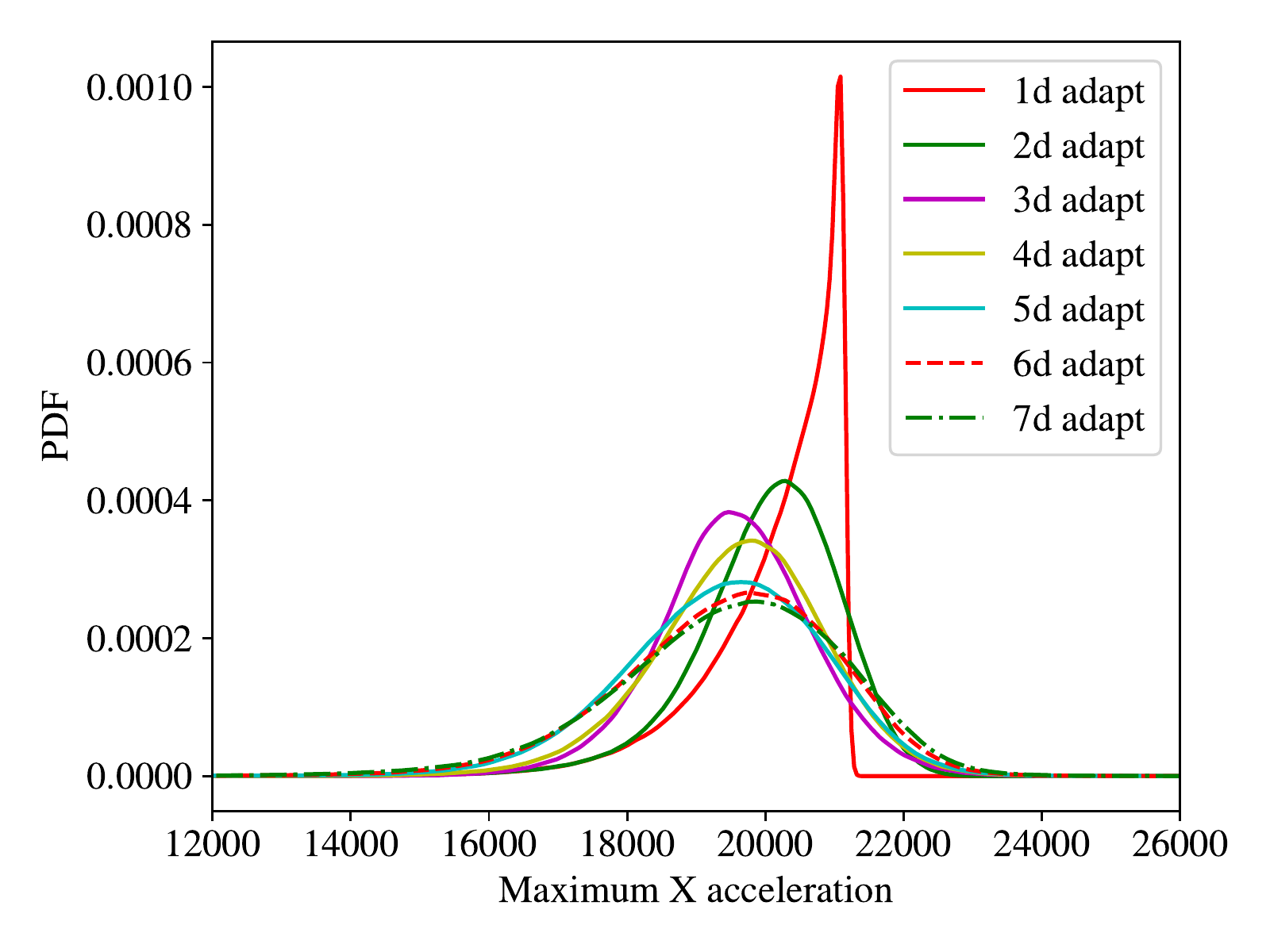}
		\caption{Classical basis adaptation of the space structure.}
		\label{fig:adapt_space}
	\end{minipage}
	\hspace{0.1cm}
	\begin{minipage}{0.45\linewidth}
		\centering
		\includegraphics[width=\linewidth]{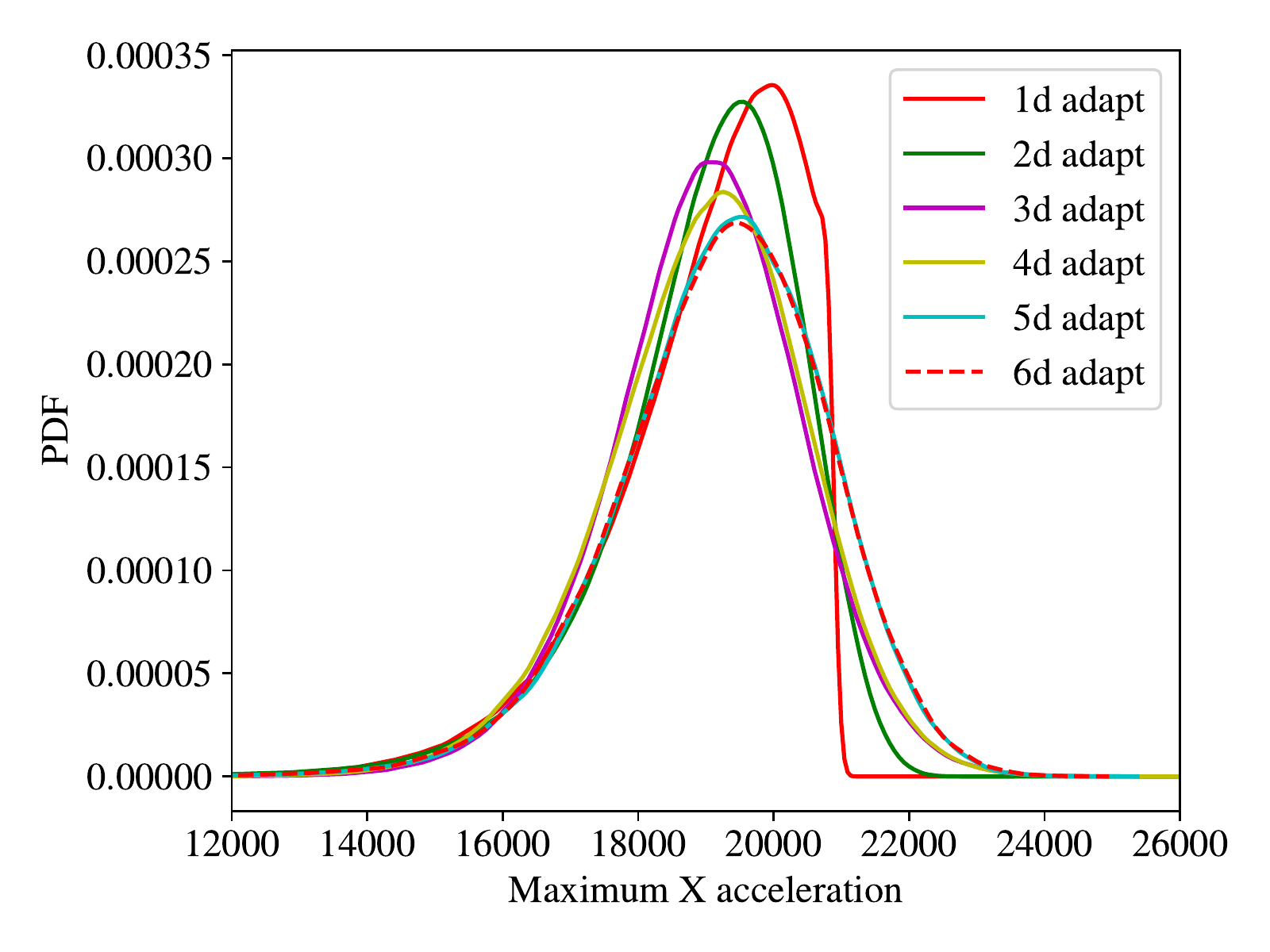}
		\caption{Accelerated basis adaptation method of the space structure.}
		\label{fig:acc_adapt_space}
	\end{minipage}
\end{figure}

In the first calculation, the QoI is the maximum $X$ acceleration along the time. Several methods can be applied to solve the forward UQ problem. Since the dimension is relatively high, a MCS requires hundreds of thousands of samples to reach a converged result. If a 3rd order PCE model is to be built, we can use a level 3 sparse Smolyak quadrature to reduce the cost; however, even in this case, the required quadrature points are $24,449$, which is still too expensive to compute. Therefore, we considered using basis adaptation and accelerated basis adaptation to solve the forward UQ problem. A first-order pilot PCE is first built from 49 level 1 quadrature points. Then, the rotation matrix adapted to the QoI is built for the classical adaptation and can be utilized to construct adaptations with sequentially increased dimensions. We update the rotation matrix for the accelerated algorithm whenever new information is fed into the model. The generated PDFs of the classical basis adaptation are shown in Figure \ref{fig:adapt_space}, where we see that it converges at dimension 7. The required quadrature points are 2731.
In contrast, the generated PDFs of the accelerated basis adaptation are shown in Figure \ref{fig:acc_adapt_space}, where we see that it converges at dimension 6. The required quadrature points are 1708. Then, for these two methods, the total model evaluation on the specific quadrature points is 2780 and 1757, respectively. The required samples are much fewer than PCE and MC methods. 

When model evaluation on specific quadrature points is unavailable, the above method will require much more MC samples to get accurate PCE models. However, the number of samples in practice is usually limited, and no more samples can be generated. Therefore, we can apply the PPA method introduced in this paper to use the available data to solve the forward UQ problem.

\begin{figure}[htb]
	\centering
	\begin{minipage}{0.45\linewidth}
		\centering
		\includegraphics[width=\linewidth]{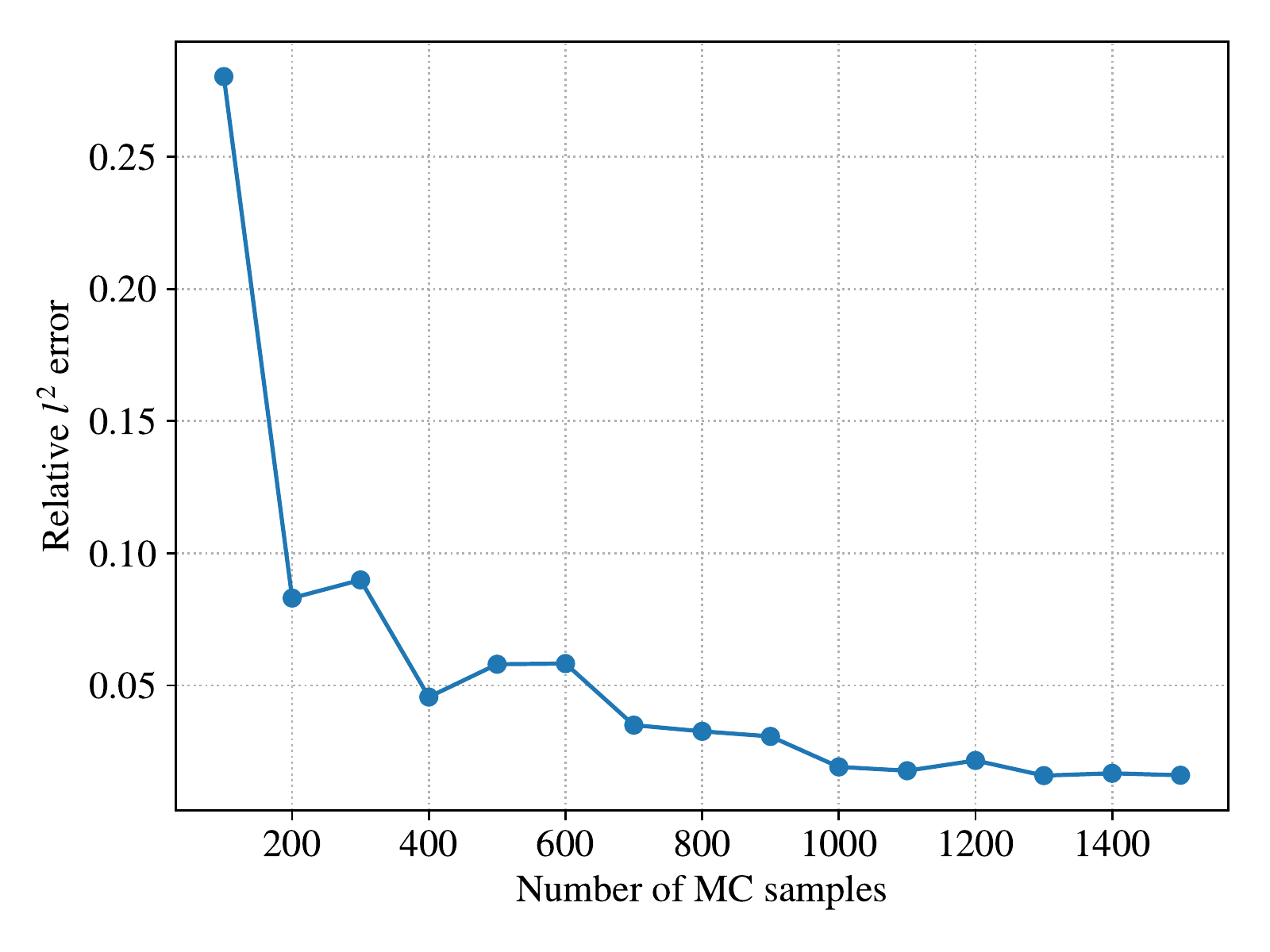}
		\caption{The convergence curve of the relative $l^2$ error of the PDF obtained from PPA with different MC samples for the space structure.}
		\label{fig:rela_l2_err_space}
	\end{minipage}
	\hspace{0.2cm}
	\begin{minipage}{0.45\linewidth}
		\centering
		\includegraphics[width=\linewidth]{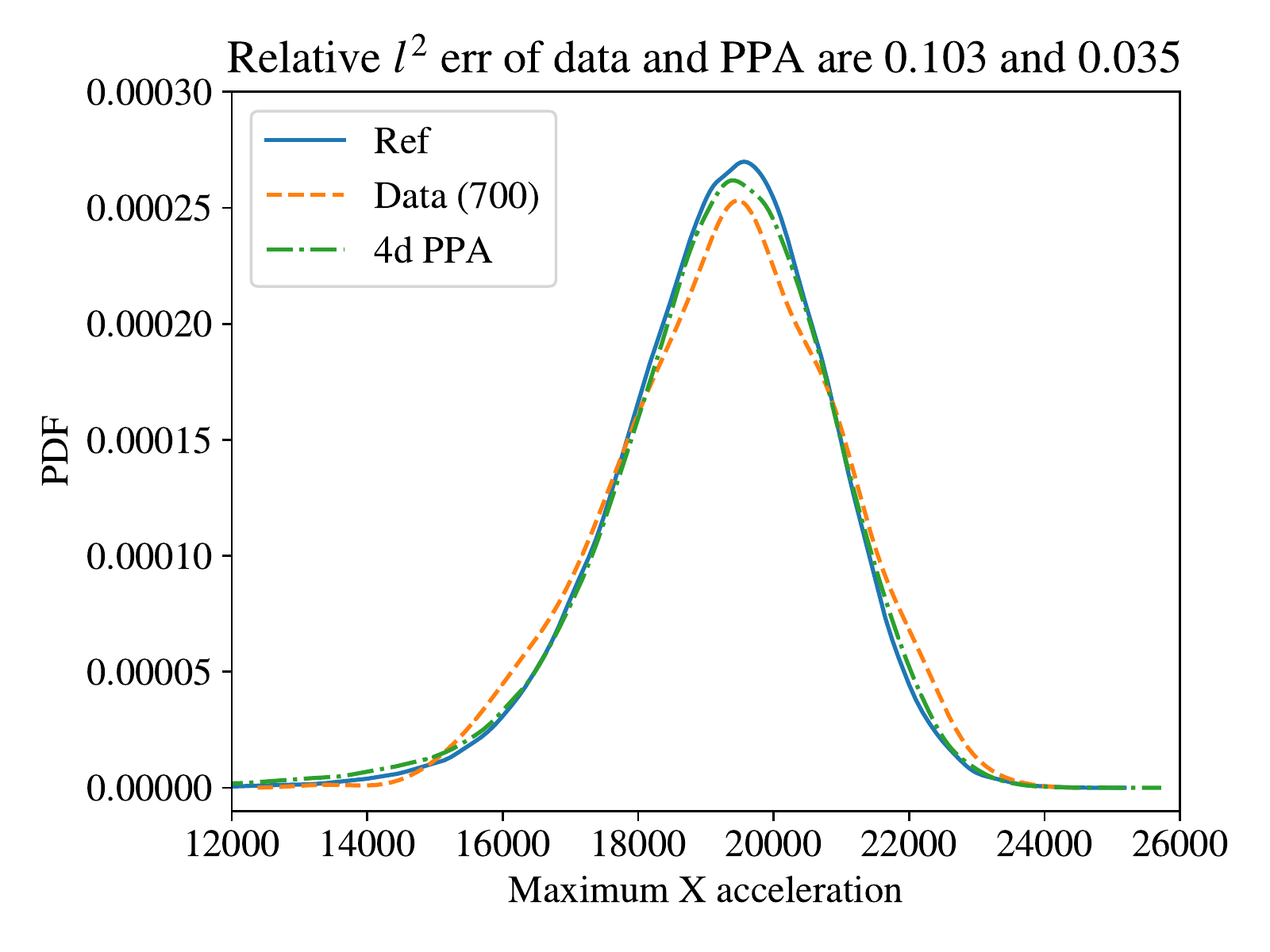}
		\caption{PDF comparison of the PPA method for the space structure. Three PDFs are compared: the reference, the data, and the converged PPA model.}
		\label{fig:pdf_compare_space_700}
	\end{minipage}
\end{figure}

As in the borehole model, we first use increasing MC samples in the data set to train converged PPA models. Then, once the PPA models were obtained, we generated many MC samples based on the model and compared the resulting PDF with the reference PDF obtained from the accelerated basis adaptation method. The relative $l^2$ errors of the PPA model obtained from different numbers of MC samples are shown in Figure \ref{fig:rela_l2_err_space}.
In this model, the relative $l^2$ errors are smaller than 5\% when the MC samples are greater than 700 and converge with 1100 MC samples. The number of MC samples required to obtain a converged PPA model is much greater than in the borehole model. Note that in practice, the reference PDF might not be available. In that case, one could compare the PDFs of the PPA model with different MC samples to check if they converge.
For most engineering applications, the accuracy of the PPA model with 700 MC samples is good enough. The PDF comparison, in this case, is shown in Figure \ref{fig:pdf_compare_space_700}. We see that the PPA PDF is close to the reference with a relative $l^2$ error of 3.5\%, while the PDF generated directly by the data has apparent discrepancies from the reference, the $l^2$ error of which is 10.3\%. 

\begin{figure}[htb]
	\centering
	\begin{minipage}{0.45\linewidth}
		\centering
		\includegraphics[width=\linewidth]{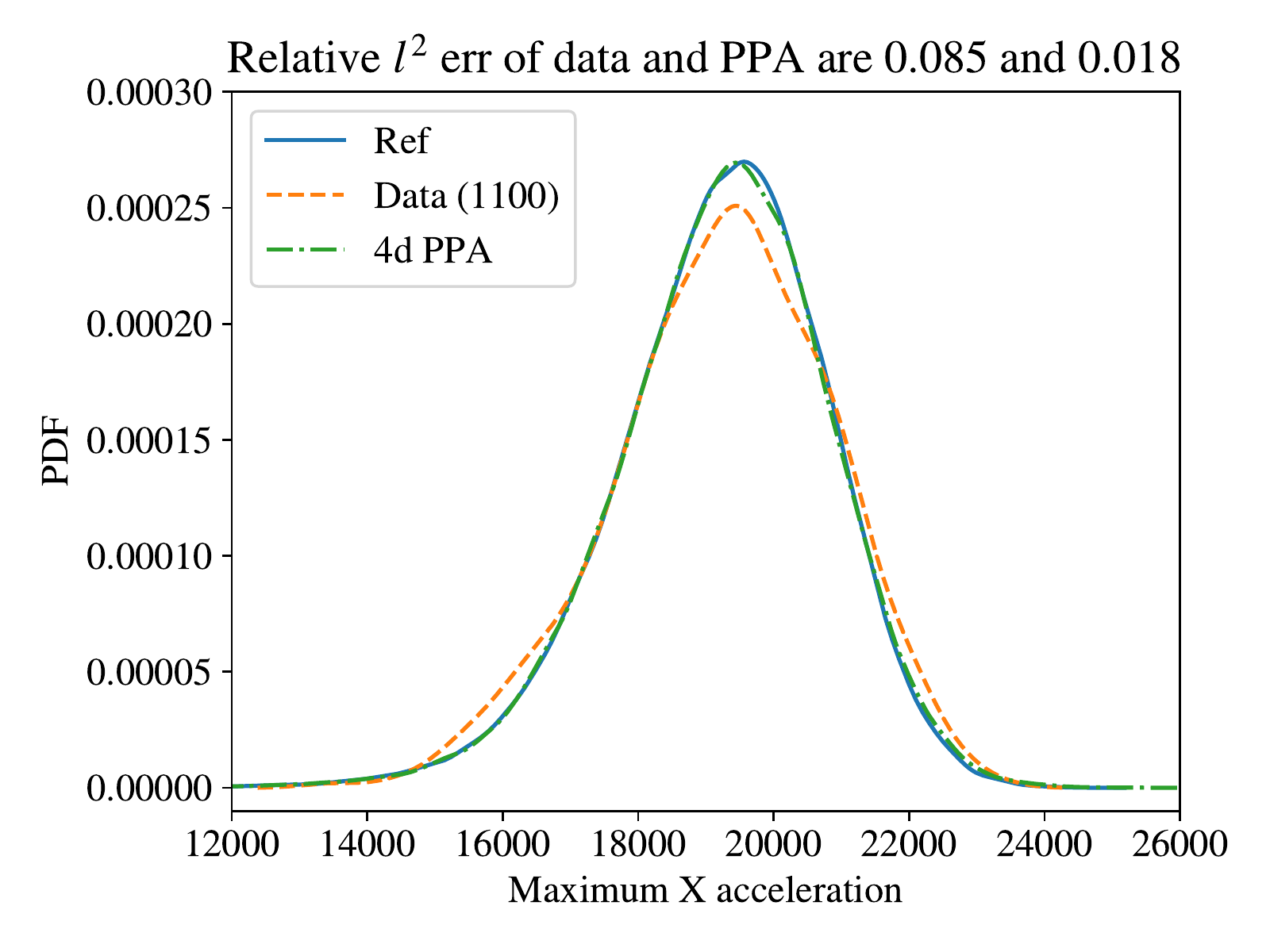}
		\caption{PDF comparison of the PPA method for the space structure. Three PDFs are compared: the reference, the data, and the converged PPA model.}
		\label{fig:pdf_compare_space}
	\end{minipage}
	\hspace{0.1cm}
	\begin{minipage}{0.45\linewidth}
		\centering
		\includegraphics[width=\linewidth]{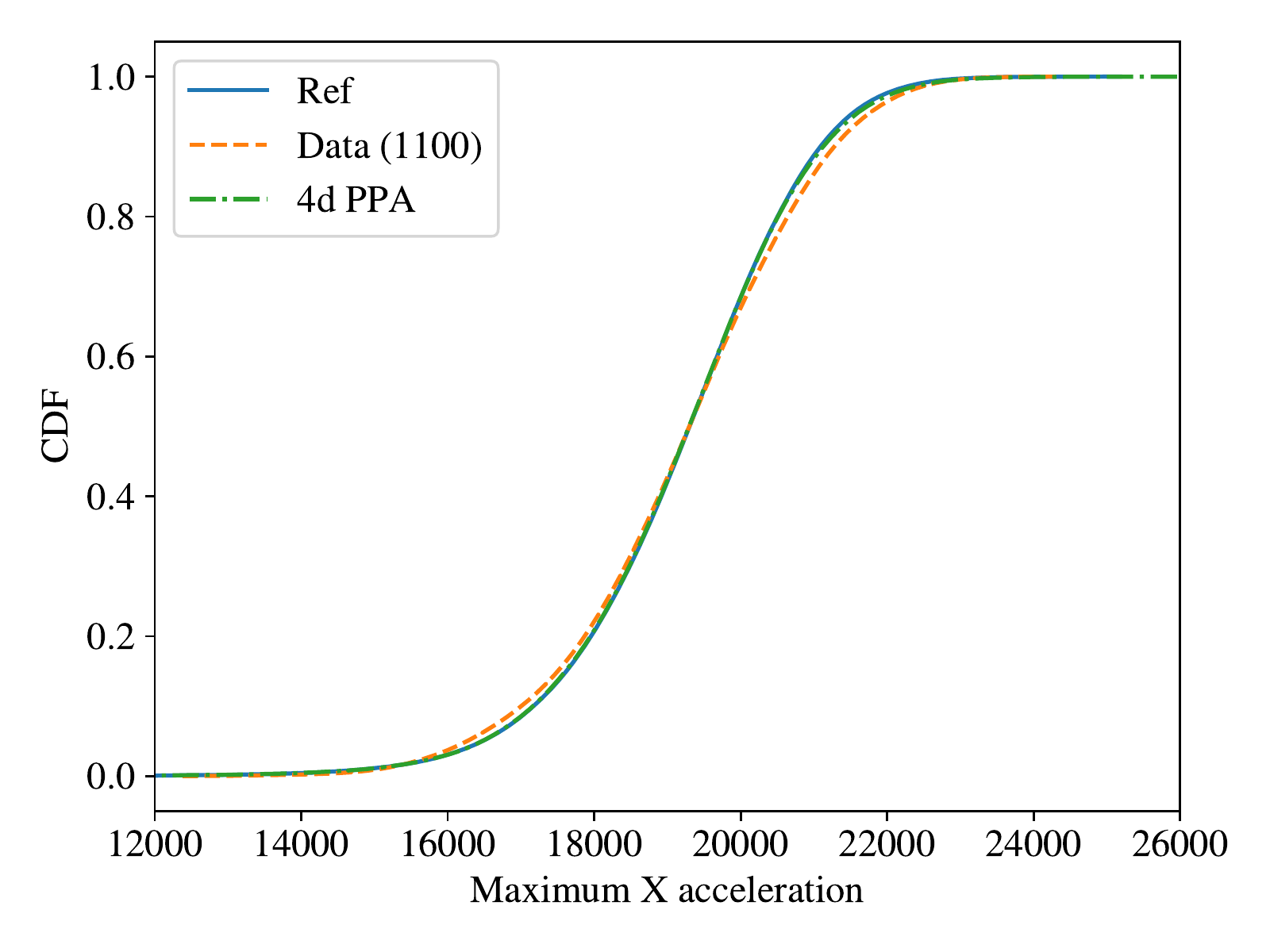}
		\caption{CDF comparison of the PPA method for the space structure. Three CDFs are compared: the reference, the data, and the converged PPA model.}
		\label{fig:cdf_compare_space}
	\end{minipage}
\end{figure}

To perform precise analyses, we will choose 1100 MC samples for further analysis. Figures \ref{fig:pdf_compare_space} and \ref{fig:cdf_compare_space} compare PDFs and CDFs obtained from different models.
We see that the resulting PDF of the PPA method is even closer to the reference than the case of using 700 MC samples. The relative $l^2$ difference between the PPA PDF and the reference PDF has been reduced to 1.9\%. The accuracy is improved in all support of the QoI. In addition, the PDF generated directly from the data is different from the reference, especially in the left tail region. Also, comparing the PDF generated from 1100 data in Figure \ref{fig:pdf_compare_space} and the PDF generated from 700 data in Figure \ref{fig:pdf_compare_space_700}, the data increase only slightly improves the PDF, implying a slow convergence of the MCS. From the CDF figure, the curve associated with the data slightly deviates from the reference, while the result from PPA almost overlaps with the reference.

\begin{figure}[htb]
	\centering
	\begin{minipage}{0.45\linewidth}
		\centering
		\includegraphics[width=\linewidth]{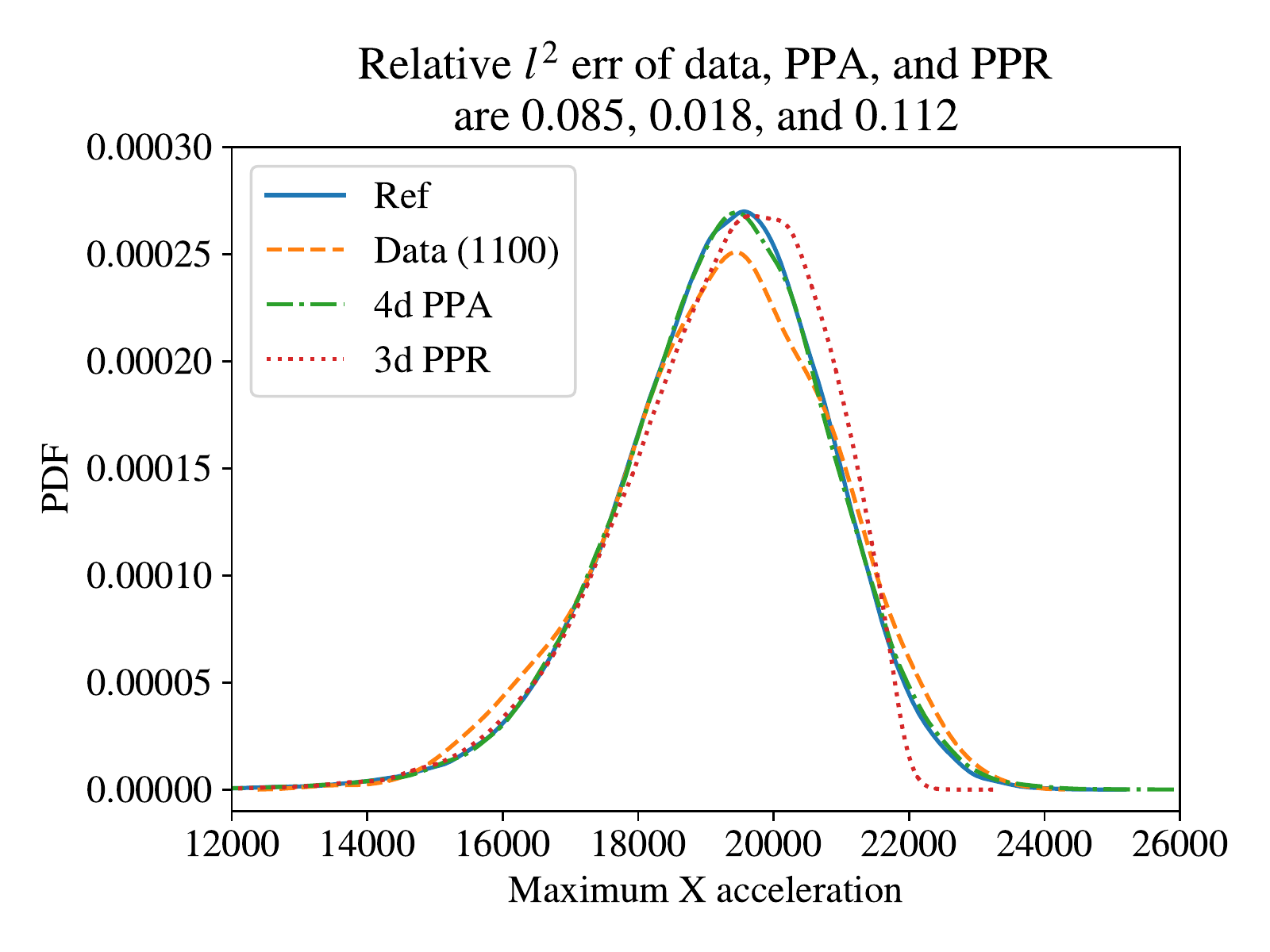}
		\caption{PDF comparison of the PPA method for the space structure. Three PDFs are compared: the reference, the data, the converged PPA model, and the converged PPR model.}
		\label{fig:pdf_compare2_space}
	\end{minipage}
	\hspace{0.1cm}
	\begin{minipage}{0.45\linewidth}
		\centering
		\includegraphics[width=\linewidth]{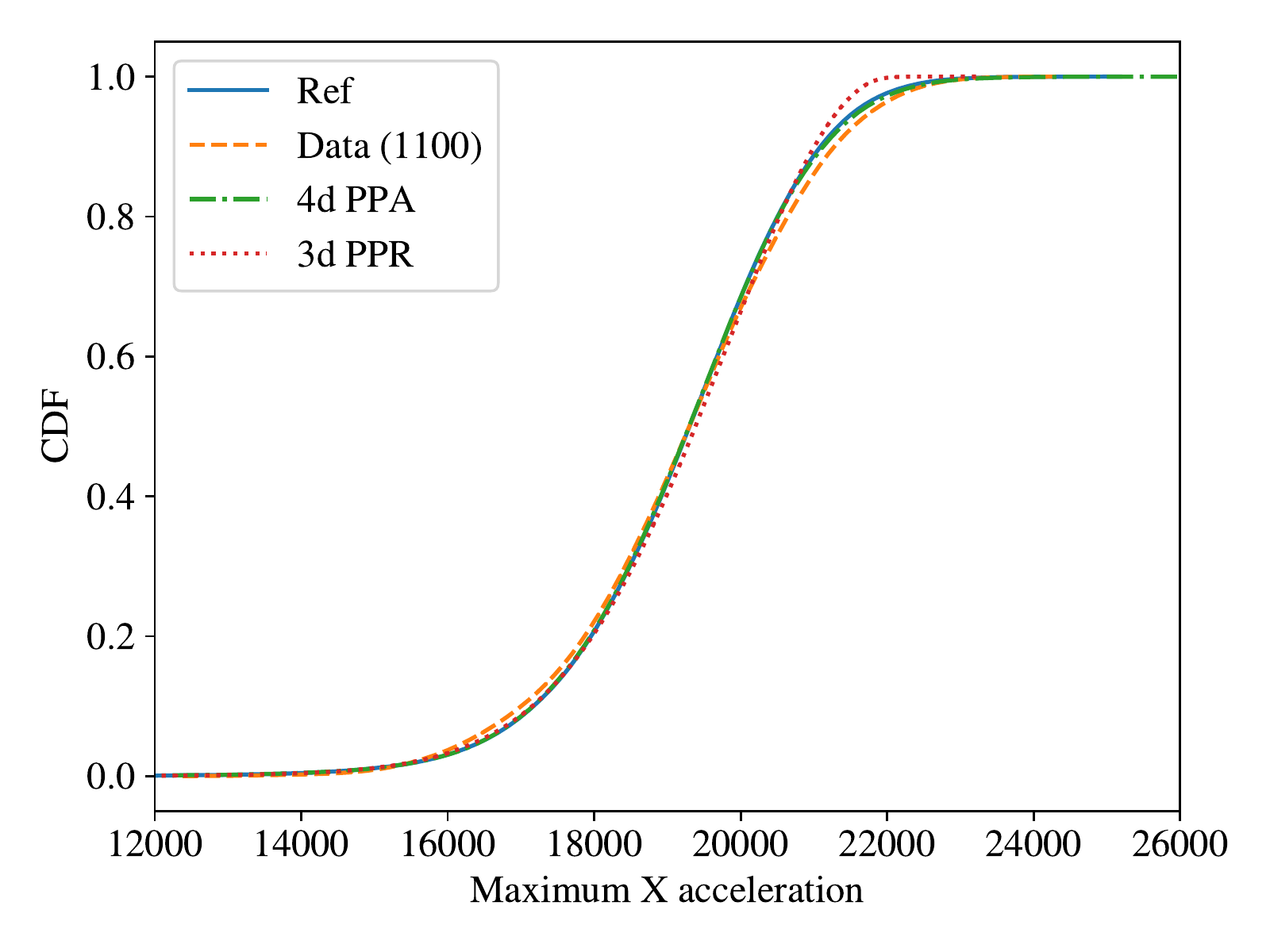}
		\caption{CDF comparison of the PPA method for the space structure. Three CDFs are compared: the reference, the data, the converged PPA model, and the converged PPR model.}
		\label{fig:cdf_compare2_space}
	\end{minipage}
\end{figure}

For this model, we also present the comparison to the PPR model, the additive of univariate PCE models. The same 1100 MC samples are used for PPA and PPR. The results of PDFs and CDFs are shown in Figures \ref{fig:pdf_compare2_space} and \ref{fig:cdf_compare2_space}.
We see that the PPR model has improved in the left half of the PDF compared to the data. However, it gives unsatisfactory predictions in the right tail regions. The relative $l^2$ error of the PPR model is 11.2\% compared to 1.8\% of the PPA model. From the CDF comparison, we can also clearly see that the PPR model gives inaccurate reliability information on the right tail. 

\begin{figure}[htb]
	\centering
	\begin{minipage}{0.45\linewidth}
		\centering
		\includegraphics[width=\linewidth]{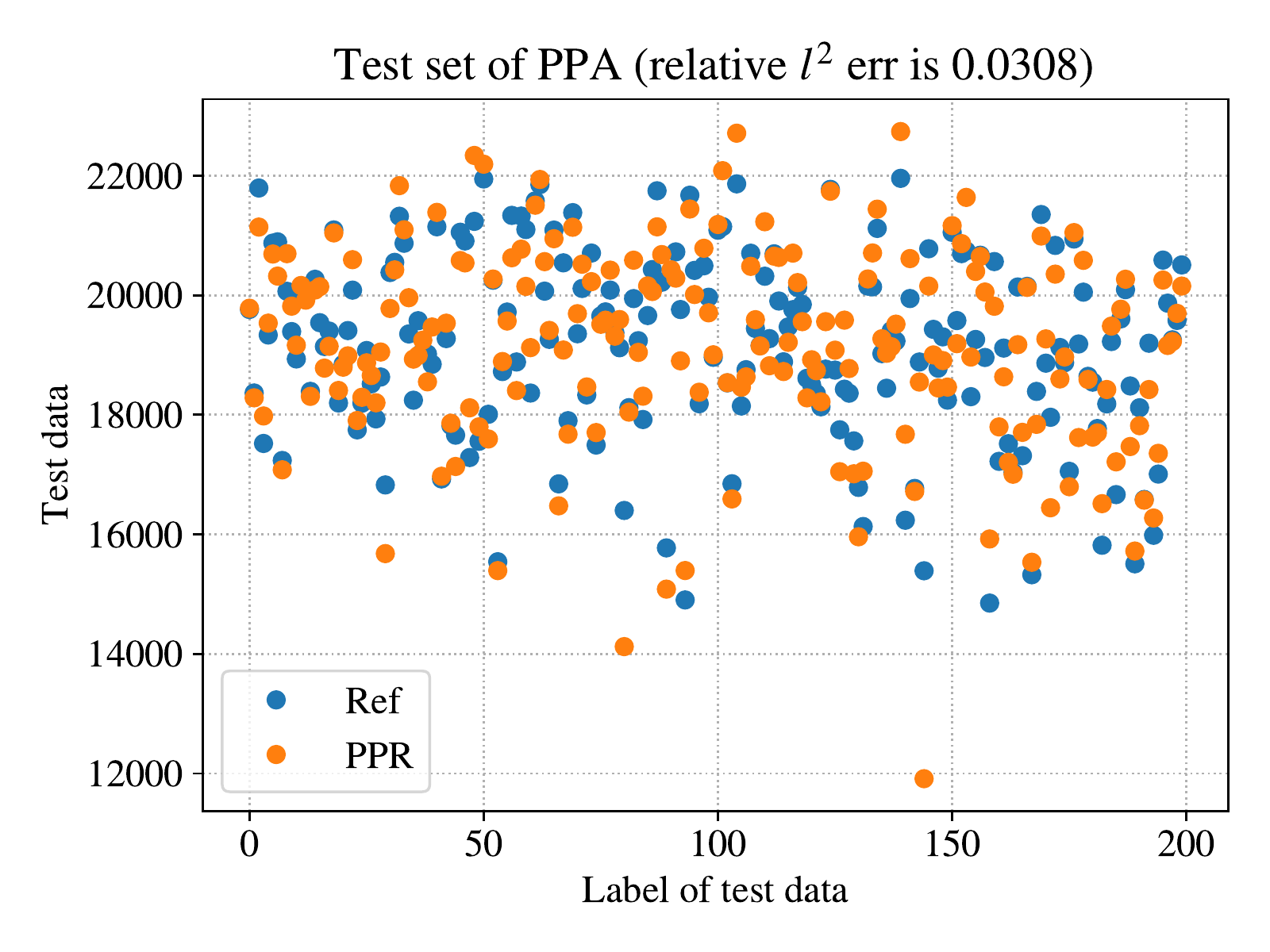}
	\end{minipage}
	\hspace{0.1cm}
	\begin{minipage}{0.45\linewidth}
		\centering
		\includegraphics[width=\linewidth]{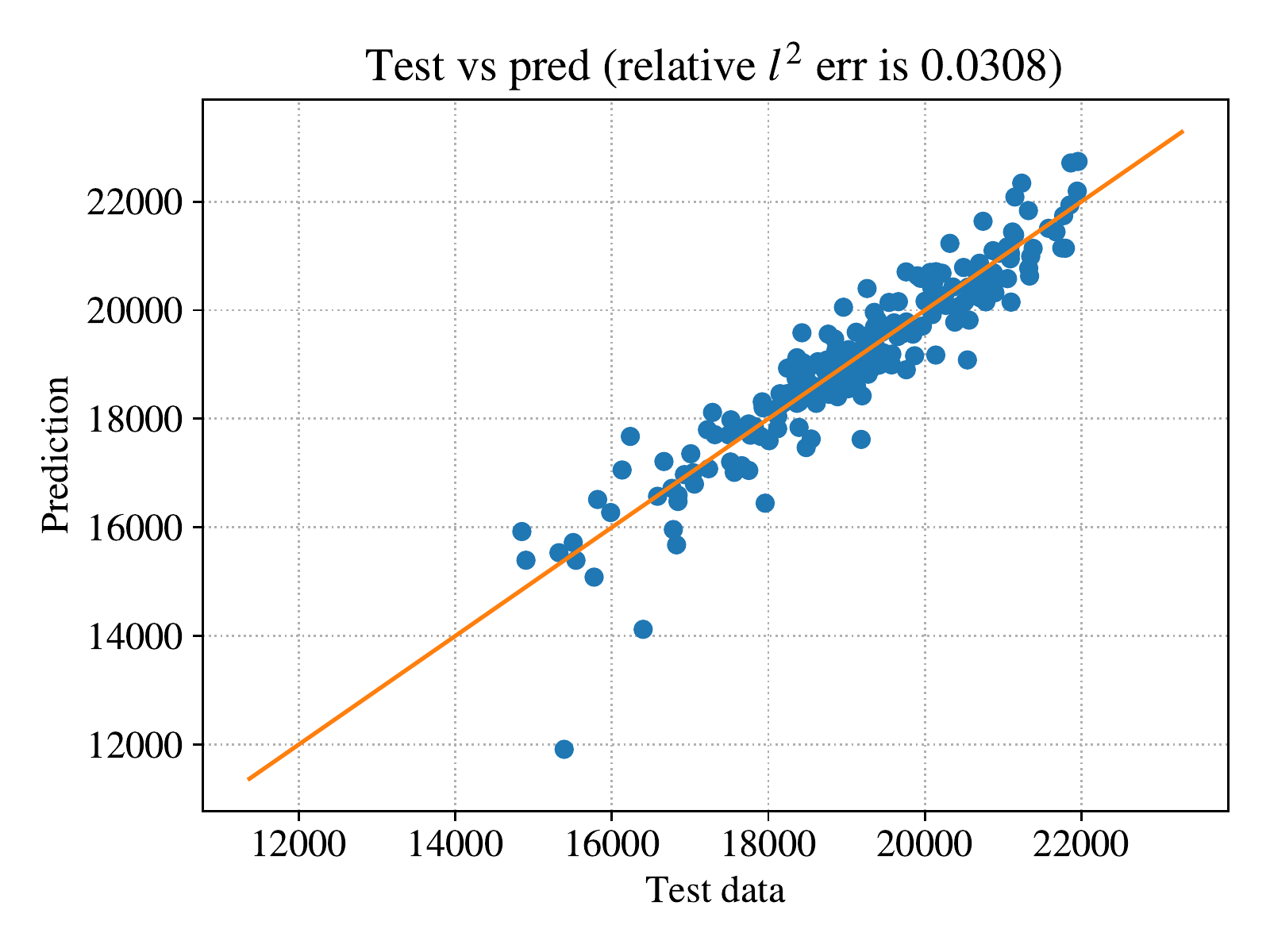}
	\end{minipage}
	\caption{Test data (200) of the PPA method for the space structure. Figure on left shows the QoI value comparison of reference and prediction on various test data. Figure on the right plot the prediction against the test data where the solid line is of slope 1. }
	\label{fig:test_space}
\end{figure}

In the case where the PPA model serves as a surrogate model to predict QoI with different input parameters, we generate 200 MC samples that are different from the training data. We can predict the QoI on these 200 samples using the PPA model. The comparison of the prediction and the reference model on the test set is shown in Figure \ref{fig:test_space}.
From the left figure, the PPA prediction is close to the reference on most of the test data. From the right figure, although several data have relatively worse predictions, the scatter dots of prediction against the test set are closely aligned in the unit slope line, meaning that the predictions are close to the reference. The relative $l^2$ difference of the 200 test data is 3.08\%, which is close to the relative $l^2$ error of the PPA PDF compared to the reference PDF. The results suggest that the prediction is accurate.

We see that the PPA model has reduced the dimension from 24 to 4 with only 1100 MC samples. The accuracy of the PPA model is comparable to the models proposed by basis adaptation and accelerated basis adaptation methods. The latter two methods adopted the sparse quadrature rule to compute the PCE coefficients and reduce computational costs. Even in that case, the required quadrature points are 2780 and 1757, respectively, which is still greater than 1100. Another observation is that the reduced dimension discovered by PPA is less than that of the basis adaptation methods, which suggests that the low-dimensional manifold found by PPA is better at revealing the probabilistic information of the QoI.

\subsubsection{Multi-QoI capability}
Another valuable feature of the data-driven essence of PPA is that the data can be reused for different QoIs. One of the limitations of the basis adaptation methods is that each rotation matrix is associated with a scalar QoI. If the QoI is changed, the whole procedure must be repeated to discover a low-dimensional space for the new QoI. This means that new model evaluations are required to compute the PCE coefficients of the newly adapted model. In the PPA, however, the data we are starting from can be reused to train a new model associated with any new QoI. The only additional computation lies in the training process, which is usually much easier, especially when a low-dimensional manifold is discovered. To test its capability, we use the same data and try to train a model for the maximum velocity of the space structure.

In this case, we are not showing the results from basis adaptation and accelerated basis adaptation. However, the accelerated basis adaptation method provides a reference model of the maximum velocity from a 6-dimensional PCE.

\begin{figure}[htb]
	\centering
	\begin{minipage}{0.45\linewidth}
		\centering
		\includegraphics[width=\linewidth]{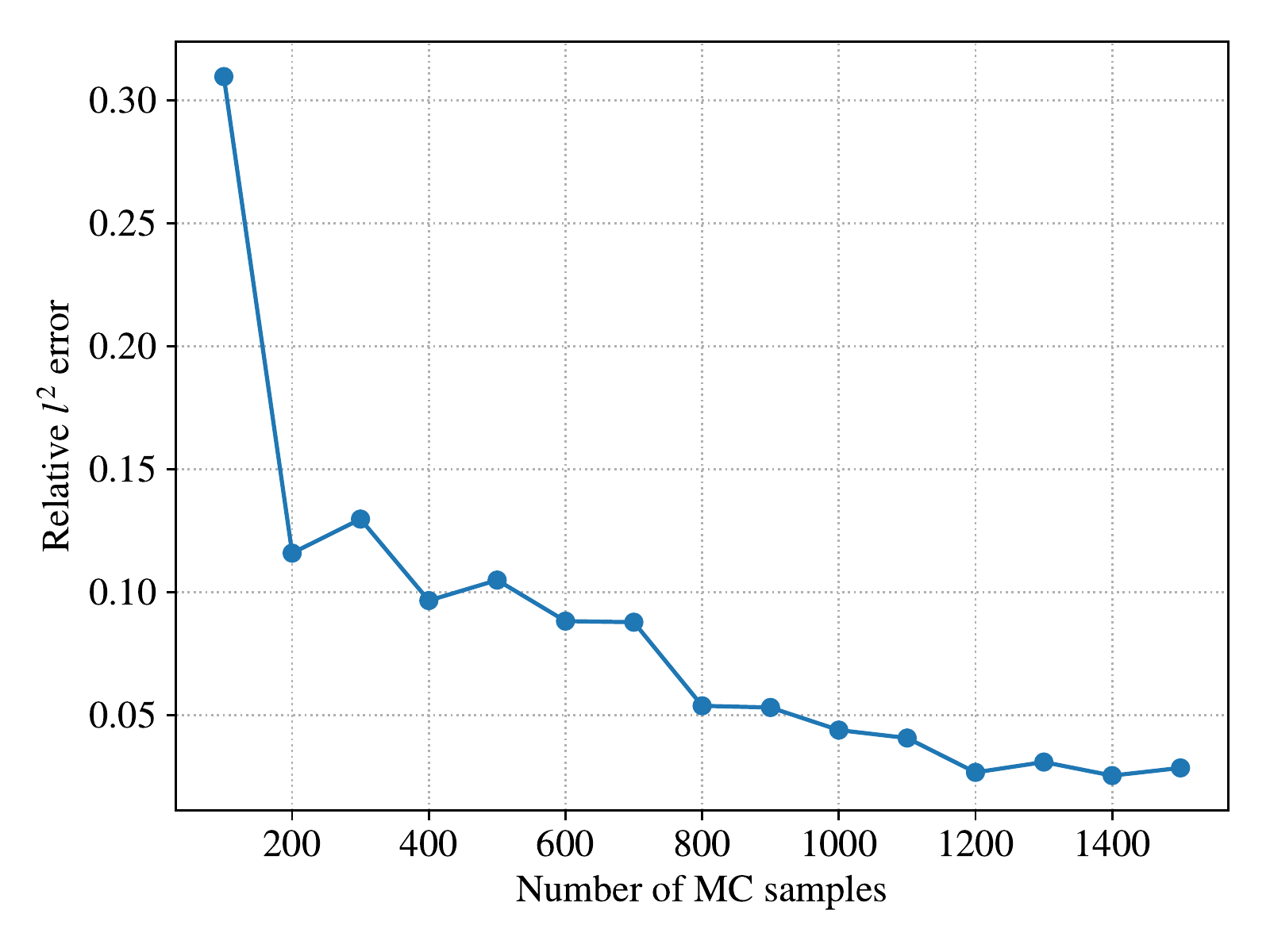}
		\caption{The convergence curve of the relative $l^2$ error of the PDF obtained from PPA with different MC samples for the space structure when QoI is the maximum velocity.}
		\label{fig:rela_l2_err_vel_space}
	\end{minipage}
	\hspace{0.2cm}
	\begin{minipage}{0.45\linewidth}
		\centering
		\includegraphics[width=\linewidth]{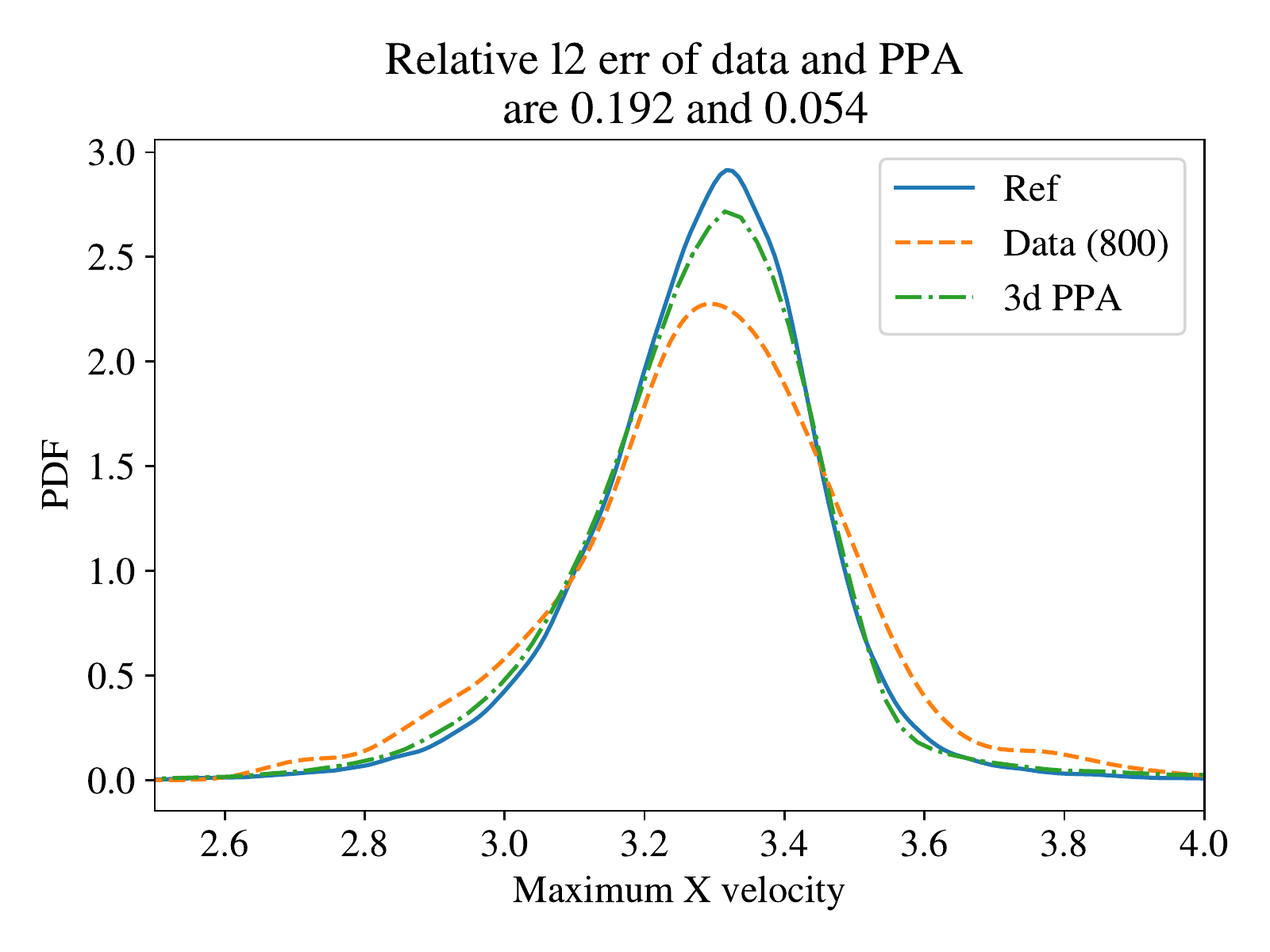}
		\caption{PDF comparison of the PPA method for the space structure when QoI is the maximum velocity. Three PDFs are compared: the reference, the data, and the converged PPA model.}
		\label{fig:pdf_compare_space_vel_800}
	\end{minipage}
\end{figure}

Starting from 100 MC samples in the data set, we gradually increase more MC samples. First, a converged PPA model is trained for each data set expressed in PCE. Then, we generate many MC samples based on the PPA model and obtain its PDF by KDE. Finally, the PDFs of the PPA from the data with different samples are compared to the reference PDF. Then, a convergence curve of the relative $l^2$ error can be obtained; see Figure \ref{fig:rela_l2_err_vel_space}.
For the maximum velocity, the PPA model has an $l^2$ error close to 5\% when the number of samples in the training data set is 800. Again, in most engineering applications, the accuracy is already enough. The PDF comparison, in this case, is shown in Figure \ref{fig:pdf_compare_space_vel_800}. Clearly, the PDF generated by the PPA model is close to the reference with an $l^2$ error of 5.4\%, while the PDF generated from the data directly has an $l^2$ error of 19.2\%, a large discrepancy from the reference.

\begin{figure}[htb]
	\centering
	\begin{minipage}{0.45\linewidth}
		\centering
		\includegraphics[width=\linewidth]{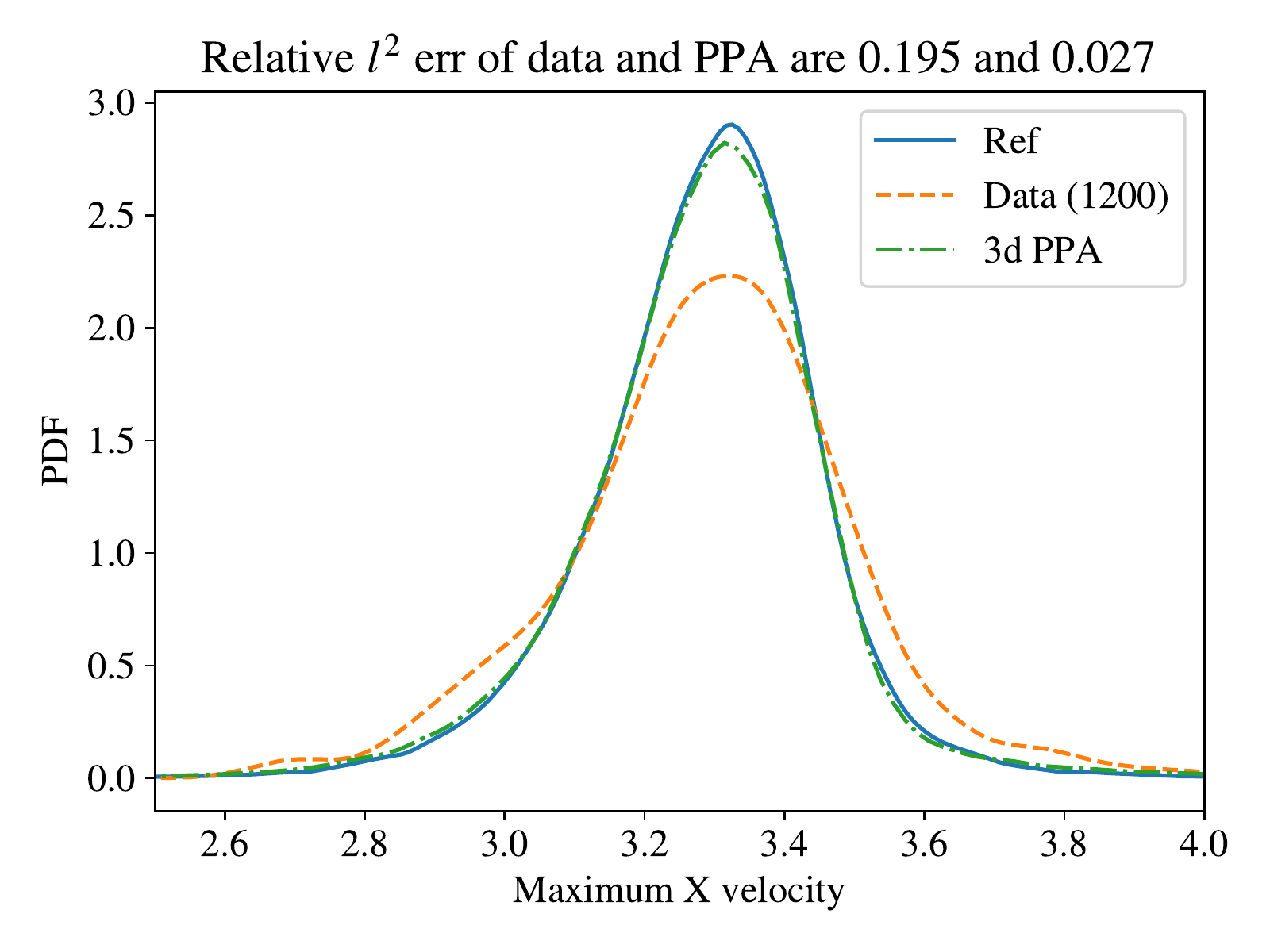}
		\caption{PDF comparison of the PPA method for the space structure when QoI is the maximum velocity. Three PDFs are compared: the reference, the data, and the converged PPA model.}
		\label{fig:pdf_compare_vel_space}
	\end{minipage}
	\hspace{0.1cm}
	\begin{minipage}{0.45\linewidth}
		\centering
		\includegraphics[width=\linewidth]{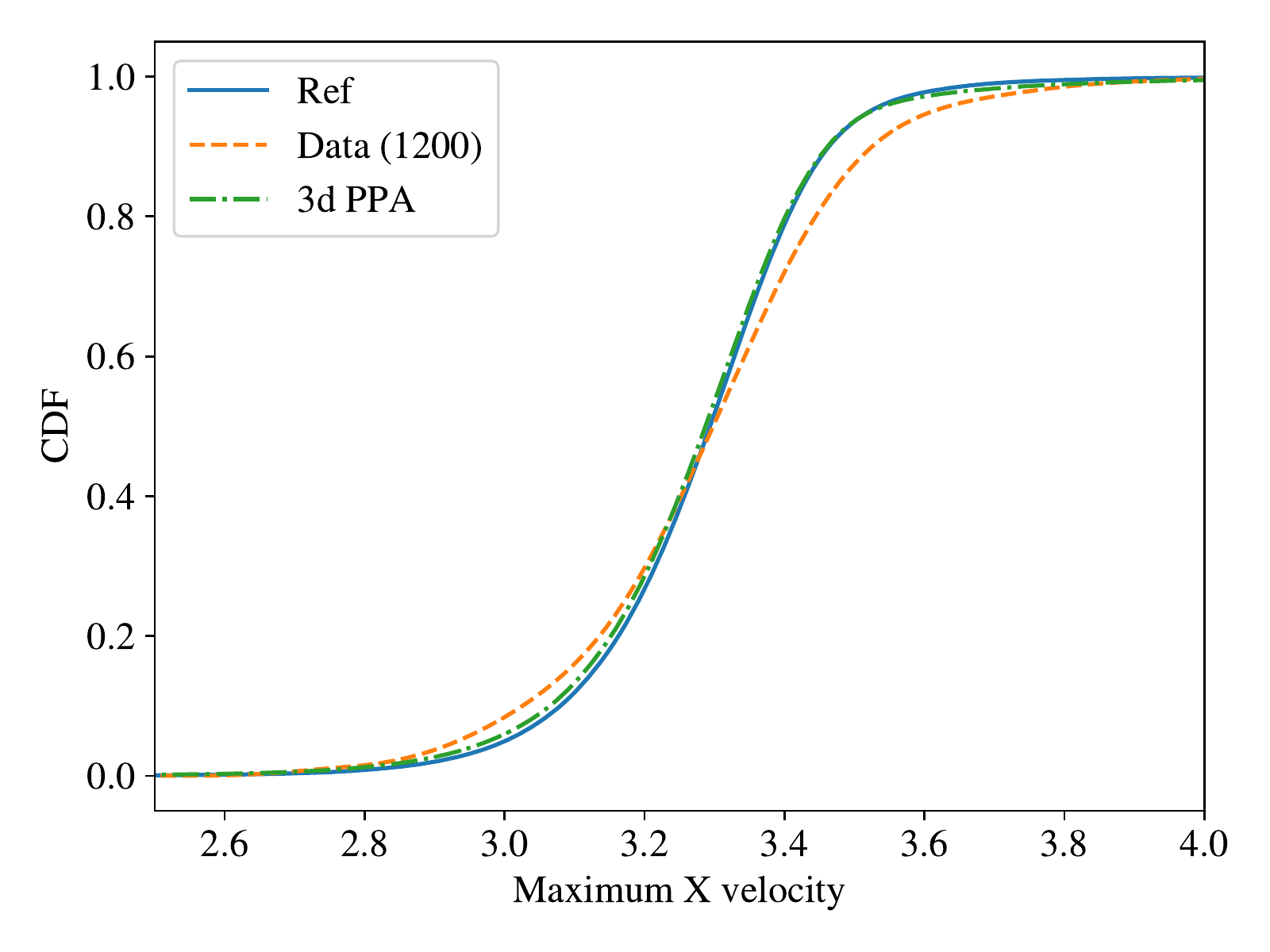}
		\caption{CDF comparison of the PPA method for the space structure when QoI is the maximum velocity. Three CDFs are compared: the reference, the data, and the converged PPA model.}
		\label{fig:cdf_compare_vel_space}
	\end{minipage}
\end{figure}

To perform precise analyses, we chose 1200 MC samples in the data set (where the PPA model is converged) for further analysis. The PDF and CDF comparisons of different models are shown in Figures \ref{fig:pdf_compare_vel_space} and \ref{fig:cdf_compare_vel_space}.
We see that the PDF obtained from the PPA model is even closer to the reference than when 800 data are used. The relative $l^2$ error has decreased from 5.4\% to 2.7\%. However, for the PDF obtained from the data, the $l^2$ error barely changes from the 800 data case to the 1200 data case. A significant discrepancy still exists even if the number of samples has increased. This signifies the slow convergence of the traditional MC methods. From Figure \ref{fig:cdf_compare_vel_space}, we also see that the CDF of the PPA model is consistent with the reference, while the CDF from the data is far from the reference.

\begin{figure}[htb]
	\centering
	\begin{minipage}{0.45\linewidth}
		\centering
		\includegraphics[width=\linewidth]{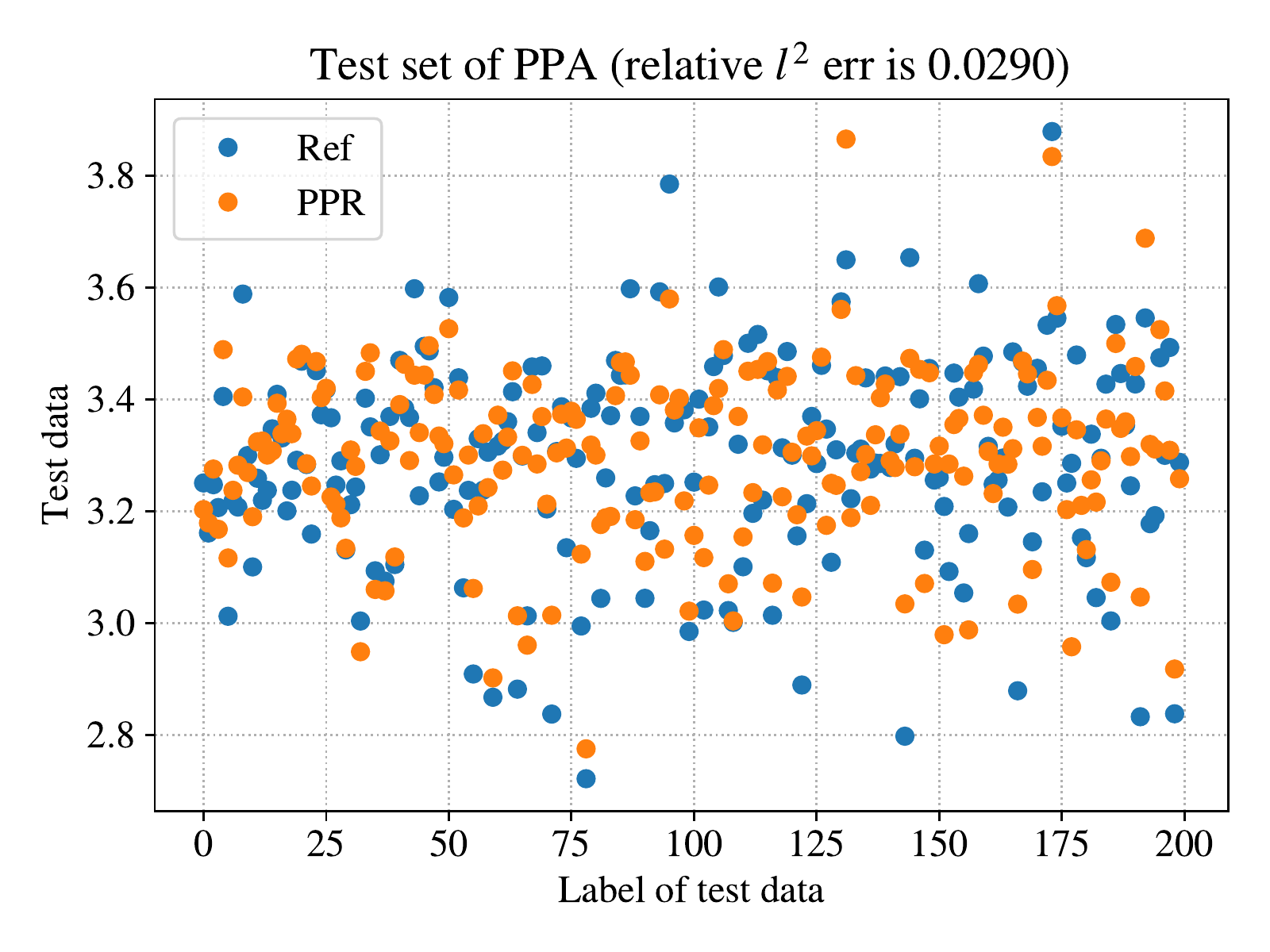}
	\end{minipage}
	\hspace{0.1cm}
	\begin{minipage}{0.45\linewidth}
		\centering
		\includegraphics[width=\linewidth]{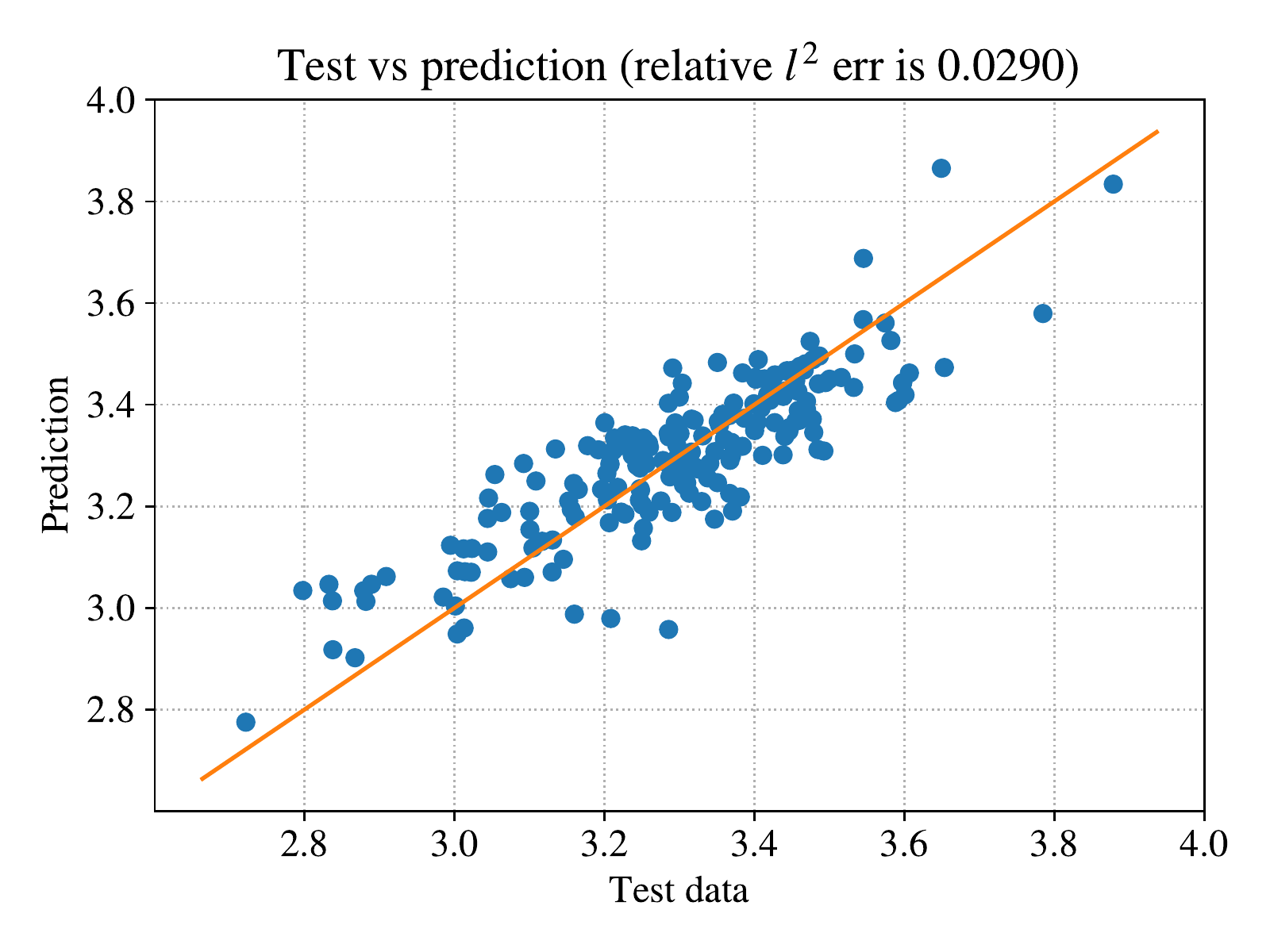}
	\end{minipage}
	\caption{Test data (200) of the PPA method for the space structure when QoI is the maximum velocity. Figure on left shows the QoI value comparison of reference and prediction on various test data. Figure on the right plot the prediction against the test data where the solid line is of slope 1. }
	\label{fig:test_vel_space}
\end{figure}

Again, we can use the PPA model as a surrogate to predict QoI with different input parameters. For example, we generate 200 MC samples that differ from the training data as the test set. We can predict the QoI on these 200 samples using the PPA model. The comparison of the prediction and the reference model on the test set is shown in Figure \ref{fig:test_vel_space}.
From the left figure, the PPA prediction is close to the reference on most of the test data. This can be seen more evident in the right figure, where the scatter dots of prediction against the test set are closely aligned in the unit slope line, meaning that the predictions are close to the reference. The relative $l^2$ difference of the 200 test data is 2.9\%, which is close to the relative $l^2$ error of the PPA PDF compared to the reference PDF. The results suggest that the prediction is accurate.

\section{Conclusions}
\label{sec:conclusions}
The paper proposed a novel, highly accurate, and efficient surrogate model for high-dimensional computational models and high-dimensional data. Specifically, we proposed a data-driven projection pursuit adaptation method (PPA) that can discover optimal low-dimensional space and polynomial chaos expansions (PCE) to represent the quantity of interest (QoI) accurately. The method has the following significances. First, the dimension reduction capability and optimal PCE representation are based on a-priori-specified data. No additional model evaluations are required. Second, the method can discover even a lower-dimensional localization than other techniques, using even less data. Third, the required data are independent samples, which have less restriction and are relatively easier to obtain than, for example, quadrature points. Lastly, the method has the multi-QoI capability as the same data set can be reused to learn different representations for different QoIs.

The method is developed for uncertainty quantification (UQ) and prediction. In UQ, the basis adaptation method in PCE and its accelerated algorithm present the potential to serve as a general dimension reduction technique. Although versatile, sparse quadrature is often used to construct the PCE models to reduce the computational cost, which requires evaluating the physical model on specific quadrature points. The projection pursuit regression (PPR) model approximates the response for high-dimensional data by a sum of smooth functions of some projected variables. The optimal projections and the smooth functions are found simultaneously from the given data. However, the convergence and capability to represent any variable in the admissible space in PPR depend on the choice of smooth functions. In this paper, we combine the advantages of these two methods and propose the novel PPA method for surrogate modeling of high-dimensional models or data. Specifically, we use multivariate PCE as a smooth function and approximate the response by a multivariate PCE on the projected variables. A stage-wise greedy strategy determines the number of projections, where we simultaneously compute the optimal projection and the optimal multivariate PCE at each stage and stop the procedure if additional information does not improve the approximation significantly. The PPA model inherits the convergence property and the capability to represent any variable in the admissible space from PCE. The applications of a borehole model and a space structure model demonstrate that the PPA model can find low-dimensional spaces to represent the responses accurately. In addition, the number of samples required to obtain accurate results in PPA is even less than the accelerated basis adaptation, where the latter uses sparse quadrature points. The space structure example also shows that the method can be applied to multiple QoIs using the same data set.

\section*{Acknowledgments}
We acknowledge support through ONR Grant No. N00014-21-1-2475 with Dr. Eric Marineau as Program Manager. Financial support is also acknowledged by the United States Department of Energy through Framework, Algorithms, and Scalable Technologies for Mathematics (FASTMATH) SciDac Institute under Grant No. DE-SC0021307 with Dr. Ceren Susut ad Program Manager.

%The financial support for this project was provided by the US Department of Energy through the Frameworks, Algorithms and Scalable Technologies for Mathematics (FASTMATH) SciDac Institute and by the US Department of Defense through the Multidisciplinary University Research Initiative (MURI) program. The findings presented herein are those of the authors and do not necessarily reflect the views of the sponsor.

\bibliographystyle{model1-num-names}

\bibliography{refs.bib}

\begin{thebibliography}{69}
\expandafter\ifx\csname natexlab\endcsname\relax\def\natexlab#1{#1}\fi
\providecommand{\url}[1]{\texttt{#1}}
\providecommand{\href}[2]{#2}
\providecommand{\path}[1]{#1}
\providecommand{\DOIprefix}{doi:}
\providecommand{\ArXivprefix}{arXiv:}
\providecommand{\URLprefix}{URL: }
\providecommand{\Pubmedprefix}{pmid:}
\providecommand{\doi}[1]{\href{http://dx.doi.org/#1}{\path{#1}}}
\providecommand{\Pubmed}[1]{\href{pmid:#1}{\path{#1}}}
\providecommand{\bibinfo}[2]{#2}
\ifx\xfnm\relax \def\xfnm[#1]{\unskip,\space#1}\fi
%Type = Book
\bibitem[{Ghanem and Spanos(1991)}]{ghanem1991}
\bibinfo{author}{R.~Ghanem}, \bibinfo{author}{P.~D. Spanos},
  \bibinfo{title}{Stochastic Finite Elements: A Spectral Approach},
  \bibinfo{publisher}{Springer-Verlag}, \bibinfo{year}{1991}.
%Type = Article
\bibitem[{Ghanem and Red-Horse(1999)}]{ghanem1999propagation}
\bibinfo{author}{R.~Ghanem}, \bibinfo{author}{J.~Red-Horse},
\newblock \bibinfo{title}{Propagation of probabilistic uncertainty in complex
  physical systems using a stochastic finite element approach},
\newblock \bibinfo{journal}{Physica D: Nonlinear Phenomena}
  \bibinfo{volume}{133} (\bibinfo{year}{1999}) \bibinfo{pages}{137--144}.
%Type = Article
\bibitem[{Sudret(2007)}]{sudret2007uncertainty}
\bibinfo{author}{B.~Sudret},
\newblock \bibinfo{title}{Uncertainty propagation and sensitivity analysis in
  mechanical models--contributions to structural reliability and stochastic
  spectral methods},
\newblock \bibinfo{journal}{Habilitationa diriger des recherches,
  Universit{\'e} Blaise Pascal, Clermont-Ferrand, France} \bibinfo{volume}{147}
  (\bibinfo{year}{2007}) \bibinfo{pages}{53}.
%Type = Article
\bibitem[{Najm(2009)}]{najm2009uncertainty}
\bibinfo{author}{H.~N. Najm},
\newblock \bibinfo{title}{Uncertainty quantification and polynomial chaos
  techniques in computational fluid dynamics},
\newblock \bibinfo{journal}{Annual review of fluid mechanics}
  \bibinfo{volume}{41} (\bibinfo{year}{2009}) \bibinfo{pages}{35--52}.
%Type = Book
\bibitem[{Soize(2017)}]{soize2017uncertainty}
\bibinfo{author}{C.~Soize}, \bibinfo{title}{Uncertainty quantification},
  \bibinfo{publisher}{Springer}, \bibinfo{year}{2017}.
%Type = Article
\bibitem[{Chen et~al.(2017)Chen, Zeng, and Peng}]{chen2017}
\bibinfo{author}{J.~Chen}, \bibinfo{author}{X.~Zeng},
  \bibinfo{author}{Y.~Peng},
\newblock \bibinfo{title}{Probabilistic analysis of wind-induced vibration
  mitigation of structures by fluid viscous dampers},
\newblock \bibinfo{journal}{Joural of Sound and Vibration}
  \bibinfo{volume}{409} (\bibinfo{year}{2017}) \bibinfo{pages}{287--305}.
%Type = Article
\bibitem[{Zeng et~al.(2017)Zeng, Peng, and Chen}]{zeng2017}
\bibinfo{author}{X.~Zeng}, \bibinfo{author}{Y.~Peng},
  \bibinfo{author}{J.~Chen},
\newblock \bibinfo{title}{Serviceability-based damping optimization of randomly
  wind-excited high-rise buildings},
\newblock \bibinfo{journal}{The structural design of tall and special
  buildings} \bibinfo{volume}{26} (\bibinfo{year}{2017})
  \bibinfo{pages}{e1371}.
%Type = Article
\bibitem[{Marzouk et~al.(2007)Marzouk, Najm, and Rahn}]{marzouk2007stochastic}
\bibinfo{author}{Y.~M. Marzouk}, \bibinfo{author}{H.~N. Najm},
  \bibinfo{author}{L.~A. Rahn},
\newblock \bibinfo{title}{Stochastic spectral methods for efficient {Bayesian}
  solution of inverse problems},
\newblock \bibinfo{journal}{Journal of Computational Physics}
  \bibinfo{volume}{224} (\bibinfo{year}{2007}) \bibinfo{pages}{560--586}.
%Type = Article
\bibitem[{Marzouk and Najm(2009)}]{marzouk2009dimensionality}
\bibinfo{author}{Y.~M. Marzouk}, \bibinfo{author}{H.~N. Najm},
\newblock \bibinfo{title}{Dimensionality reduction and polynomial chaos
  acceleration of {Bayesian} inference in inverse problems},
\newblock \bibinfo{journal}{Journal of Computational Physics}
  \bibinfo{volume}{228} (\bibinfo{year}{2009}) \bibinfo{pages}{1862--1902}.
%Type = Article
\bibitem[{Marzouk and Xiu(2009)}]{marzouk2009stochastic}
\bibinfo{author}{Y.~Marzouk}, \bibinfo{author}{D.~Xiu},
\newblock \bibinfo{title}{A stochastic collocation approach to {Bayesian}
  inference in inverse problems}  (\bibinfo{year}{2009}).
%Type = Inproceedings
\bibitem[{Sudret et~al.(2017)Sudret, Marelli, and Wiart}]{sudret2017surrogate}
\bibinfo{author}{B.~Sudret}, \bibinfo{author}{S.~Marelli},
  \bibinfo{author}{J.~Wiart},
\newblock \bibinfo{title}{Surrogate models for uncertainty quantification: An
  overview},
\newblock in: \bibinfo{booktitle}{2017 11th European conference on antennas and
  propagation (EUCAP)}, \bibinfo{organization}{IEEE}, \bibinfo{year}{2017}, pp.
  \bibinfo{pages}{793--797}.
%Type = Article
\bibitem[{Ghanem(1999)}]{ghanem1999}
\bibinfo{author}{R.~Ghanem},
\newblock \bibinfo{title}{Ingredients for a general purpose stochastic finite
  elements implementation},
\newblock \bibinfo{journal}{Computer Methods in Applied Mechanics and
  Engineering} \bibinfo{volume}{168} (\bibinfo{year}{1999})
  \bibinfo{pages}{19--34}.
%Type = Article
\bibitem[{Xiu and Karniadakis(2002{\natexlab{a}})}]{xiu2002}
\bibinfo{author}{D.~Xiu}, \bibinfo{author}{G.~E. Karniadakis},
\newblock \bibinfo{title}{The {Wiener}--{Askey} polynomial chaos for stochastic
  differential equations},
\newblock \bibinfo{journal}{SIAM journal on scientific computing}
  \bibinfo{volume}{24} (\bibinfo{year}{2002}{\natexlab{a}})
  \bibinfo{pages}{619--644}.
%Type = Article
\bibitem[{Xiu and Karniadakis(2002{\natexlab{b}})}]{xiu2002-2}
\bibinfo{author}{D.~Xiu}, \bibinfo{author}{G.~E. Karniadakis},
\newblock \bibinfo{title}{Modeling uncertainty in steady state diffusion
  problems via generalized polynomial chaos},
\newblock \bibinfo{journal}{Computer methods in applied mechanics and
  engineering} \bibinfo{volume}{191} (\bibinfo{year}{2002}{\natexlab{b}})
  \bibinfo{pages}{4927--4948}.
%Type = Article
\bibitem[{MacKay et~al.(1998)}]{mackay1998introduction}
\bibinfo{author}{D.~J. MacKay}, et~al.,
\newblock \bibinfo{title}{Introduction to gaussian processes},
\newblock \bibinfo{journal}{NATO ASI series F computer and systems sciences}
  \bibinfo{volume}{168} (\bibinfo{year}{1998}) \bibinfo{pages}{133--166}.
%Type = Article
\bibitem[{Seeger(2004)}]{seeger2004gaussian}
\bibinfo{author}{M.~Seeger},
\newblock \bibinfo{title}{Gaussian processes for machine learning},
\newblock \bibinfo{journal}{International journal of neural systems}
  \bibinfo{volume}{14} (\bibinfo{year}{2004}) \bibinfo{pages}{69--106}.
%Type = Article
\bibitem[{Bilionis and Zabaras(2012)}]{bilionis2012multi}
\bibinfo{author}{I.~Bilionis}, \bibinfo{author}{N.~Zabaras},
\newblock \bibinfo{title}{Multi-output local gaussian process regression:
  Applications to uncertainty quantification},
\newblock \bibinfo{journal}{Journal of Computational Physics}
  \bibinfo{volume}{231} (\bibinfo{year}{2012}) \bibinfo{pages}{5718--5746}.
%Type = Book
\bibitem[{Li and Chen(2009)}]{li2009}
\bibinfo{author}{J.~Li}, \bibinfo{author}{J.~Chen}, \bibinfo{title}{Stochastic
  Dynamics of Structures}, \bibinfo{publisher}{John Wiley \& Sons},
  \bibinfo{year}{2009}.
%Type = Article
\bibitem[{Chen et~al.(2016)Chen, Yang, and Li}]{chen2016}
\bibinfo{author}{J.~Chen}, \bibinfo{author}{J.~Yang}, \bibinfo{author}{J.~Li},
\newblock \bibinfo{title}{A {GF}-discrepancy for point selection in stochastic
  seismic response analysis of structures with uncertain parameters},
\newblock \bibinfo{journal}{Structural Safety} \bibinfo{volume}{59}
  (\bibinfo{year}{2016}) \bibinfo{pages}{20--31}.
%Type = Article
\bibitem[{Soize and Ghanem(2021)}]{soize2021probabilistic}
\bibinfo{author}{C.~Soize}, \bibinfo{author}{R.~Ghanem},
\newblock \bibinfo{title}{Probabilistic learning on manifolds constrained by
  nonlinear partial differential equations for small datasets},
\newblock \bibinfo{journal}{Computer Methods in Applied Mechanics and
  Engineering} \bibinfo{volume}{380} (\bibinfo{year}{2021})
  \bibinfo{pages}{113777}.
%Type = Article
\bibitem[{Zhang et~al.(2021)Zhang, Guilleminot, and
  Gomez}]{zhang2021stochastic}
\bibinfo{author}{H.~Zhang}, \bibinfo{author}{J.~Guilleminot},
  \bibinfo{author}{L.~J. Gomez},
\newblock \bibinfo{title}{Stochastic modeling of geometrical uncertainties on
  complex domains, with application to additive manufacturing and brain
  interface geometries},
\newblock \bibinfo{journal}{Computer Methods in Applied Mechanics and
  Engineering} \bibinfo{volume}{385} (\bibinfo{year}{2021})
  \bibinfo{pages}{114014}.
%Type = Article
\bibitem[{Giovanis and Shields(2020)}]{giovanis2020data}
\bibinfo{author}{D.~G. Giovanis}, \bibinfo{author}{M.~D. Shields},
\newblock \bibinfo{title}{Data-driven surrogates for high dimensional models
  using gaussian process regression on the grassmann manifold},
\newblock \bibinfo{journal}{Computer Methods in Applied Mechanics and
  Engineering} \bibinfo{volume}{370} (\bibinfo{year}{2020})
  \bibinfo{pages}{113269}.
%Type = Article
\bibitem[{Kougioumtzoglou and Spanos(2012)}]{kougioumtzoglou2012analytical}
\bibinfo{author}{I.~Kougioumtzoglou}, \bibinfo{author}{P.~Spanos},
\newblock \bibinfo{title}{An analytical wiener path integral technique for
  non-stationary response determination of nonlinear oscillators},
\newblock \bibinfo{journal}{Probabilistic Engineering Mechanics}
  \bibinfo{volume}{28} (\bibinfo{year}{2012}) \bibinfo{pages}{125--131}.
%Type = Article
\bibitem[{Psaros et~al.(2019)Psaros, Petromichelakis, and
  Kougioumtzoglou}]{psaros2019wiener}
\bibinfo{author}{A.~F. Psaros}, \bibinfo{author}{I.~Petromichelakis},
  \bibinfo{author}{I.~A. Kougioumtzoglou},
\newblock \bibinfo{title}{Wiener path integrals and multi-dimensional global
  bases for non-stationary stochastic response determination of structural
  systems},
\newblock \bibinfo{journal}{Mechanical Systems and Signal Processing}
  \bibinfo{volume}{128} (\bibinfo{year}{2019}) \bibinfo{pages}{551--571}.
%Type = Article
\bibitem[{Babu{\v{s}}ka et~al.(2007)Babu{\v{s}}ka, Nobile, and
  Tempone}]{babuvska2007}
\bibinfo{author}{I.~Babu{\v{s}}ka}, \bibinfo{author}{F.~Nobile},
  \bibinfo{author}{R.~Tempone},
\newblock \bibinfo{title}{A stochastic collocation method for elliptic partial
  differential equations with random input data},
\newblock \bibinfo{journal}{SIAM Journal on Numerical Analysis}
  \bibinfo{volume}{45} (\bibinfo{year}{2007}) \bibinfo{pages}{1005--1034}.
%Type = Book
\bibitem[{Le~Maitre and Knio(2010)}]{maitre2010}
\bibinfo{author}{O.~Le~Maitre}, \bibinfo{author}{O.~M. Knio},
  \bibinfo{title}{Spectral methods for uncertainty quantification: with
  applications to computational fluid dynamics}, \bibinfo{publisher}{Springer
  Science \& Business Media}, \bibinfo{year}{2010}.
%Type = Article
\bibitem[{Smolyak(1963)}]{smolyak1963}
\bibinfo{author}{S.~A. Smolyak},
\newblock \bibinfo{title}{Quadrature and interpolation formulas for tensor
  products of certain classes of functions},
\newblock \bibinfo{journal}{Doklady Akademii Nauk} \bibinfo{volume}{4}
  (\bibinfo{year}{1963}) \bibinfo{pages}{240--243}.
%Type = Article
\bibitem[{Gerstner and Griebel(1998)}]{gerstner1998}
\bibinfo{author}{T.~Gerstner}, \bibinfo{author}{M.~Griebel},
\newblock \bibinfo{title}{Numerical integration using sparse grids},
\newblock \bibinfo{journal}{Numerical algorithms} \bibinfo{volume}{18}
  (\bibinfo{year}{1998}) \bibinfo{pages}{209--232}.
%Type = Article
\bibitem[{Novak and Ritter(1999)}]{novak1999}
\bibinfo{author}{E.~Novak}, \bibinfo{author}{K.~Ritter},
\newblock \bibinfo{title}{Simple cubature formulas with high polynomial
  exactness},
\newblock \bibinfo{journal}{Constructive approximation} \bibinfo{volume}{15}
  (\bibinfo{year}{1999}) \bibinfo{pages}{499--522}.
%Type = Article
\bibitem[{Blatman and Sudret(2010)}]{blatman2010}
\bibinfo{author}{G.~Blatman}, \bibinfo{author}{B.~Sudret},
\newblock \bibinfo{title}{An adaptive algorithm to build up sparse polynomial
  chaos expansions for stochastic finite element analysis},
\newblock \bibinfo{journal}{Probabilistic Engineering Mechanics}
  \bibinfo{volume}{25} (\bibinfo{year}{2010}) \bibinfo{pages}{183--197}.
%Type = Article
\bibitem[{Blatman and Sudret(2011)}]{blatman2011}
\bibinfo{author}{G.~Blatman}, \bibinfo{author}{B.~Sudret},
\newblock \bibinfo{title}{Adaptive sparse polynomial chaos expansion based on
  least angle regression},
\newblock \bibinfo{journal}{Journal of Computational Physics}
  \bibinfo{volume}{230} (\bibinfo{year}{2011}) \bibinfo{pages}{2345--2367}.
%Type = Article
\bibitem[{Doostan and Owhadi(2011)}]{doostan2011}
\bibinfo{author}{A.~Doostan}, \bibinfo{author}{H.~Owhadi},
\newblock \bibinfo{title}{A non-adapted sparse approximation of {PDEs} with
  stochastic inputs},
\newblock \bibinfo{journal}{Journal of Computational Physics}
  \bibinfo{volume}{230} (\bibinfo{year}{2011}) \bibinfo{pages}{3015--3034}.
%Type = Article
\bibitem[{Sargsyan et~al.(2014)Sargsyan, Safta, Najm, Debusschere, Ricciuto,
  and Thornton}]{sargsyan2014}
\bibinfo{author}{K.~Sargsyan}, \bibinfo{author}{C.~Safta},
  \bibinfo{author}{H.~N. Najm}, \bibinfo{author}{B.~J. Debusschere},
  \bibinfo{author}{D.~Ricciuto}, \bibinfo{author}{P.~Thornton},
\newblock \bibinfo{title}{Dimensionality reduction for complex models via
  {Bayesian} compressive sensing},
\newblock \bibinfo{journal}{International Journal for Uncertainty
  Quantification} \bibinfo{volume}{4} (\bibinfo{year}{2014})
  \bibinfo{pages}{63--93}.
%Type = Article
\bibitem[{Hampton and Doostan(2015)}]{hampton2015compressive}
\bibinfo{author}{J.~Hampton}, \bibinfo{author}{A.~Doostan},
\newblock \bibinfo{title}{Compressive sampling of polynomial chaos expansions:
  Convergence analysis and sampling strategies},
\newblock \bibinfo{journal}{Journal of Computational Physics}
  \bibinfo{volume}{280} (\bibinfo{year}{2015}) \bibinfo{pages}{363--386}.
%Type = Article
\bibitem[{Constantine et~al.(2014)Constantine, Dow, and Wang}]{constantine2014}
\bibinfo{author}{P.~G. Constantine}, \bibinfo{author}{E.~Dow},
  \bibinfo{author}{Q.~Wang},
\newblock \bibinfo{title}{Active subspace methods in theory and practice:
  applications to kriging surfaces},
\newblock \bibinfo{journal}{SIAM Journal on Scientific Computing}
  \bibinfo{volume}{36} (\bibinfo{year}{2014}) \bibinfo{pages}{A1500--A1524}.
%Type = Book
\bibitem[{Constantine(2015)}]{constantine2015}
\bibinfo{author}{P.~G. Constantine}, \bibinfo{title}{Active Subspaces: Emerging
  Ideas for Dimension Reduction in Parameter Studies},
  volume~\bibinfo{volume}{2}, \bibinfo{publisher}{SIAM}, \bibinfo{year}{2015}.
%Type = Article
\bibitem[{Tipireddy and Ghanem(2014)}]{tipireddy2014}
\bibinfo{author}{R.~Tipireddy}, \bibinfo{author}{R.~Ghanem},
\newblock \bibinfo{title}{Basis adaptation in homogeneous chaos spaces},
\newblock \bibinfo{journal}{Journal of Computational Physics}
  \bibinfo{volume}{259} (\bibinfo{year}{2014}) \bibinfo{pages}{304--317}.
%Type = Article
\bibitem[{Thimmisetty et~al.(2017)Thimmisetty, Tsilifis, and
  Ghanem}]{thimmisetty2017}
\bibinfo{author}{C.~Thimmisetty}, \bibinfo{author}{P.~Tsilifis},
  \bibinfo{author}{R.~Ghanem},
\newblock \bibinfo{title}{Homogeneous chaos basis adaptation for design
  optimization under uncertainty: {Application} to the oil well placement
  problem},
\newblock \bibinfo{journal}{Ai Edam} \bibinfo{volume}{31}
  (\bibinfo{year}{2017}) \bibinfo{pages}{265--276}.
%Type = Article
\bibitem[{Ghauch et~al.(2019)Ghauch, Aitharaju, Rodgers, Pasupuleti, Dereims,
  and Ghanem}]{ghauch2019}
\bibinfo{author}{Z.~G. Ghauch}, \bibinfo{author}{V.~Aitharaju},
  \bibinfo{author}{W.~R. Rodgers}, \bibinfo{author}{P.~Pasupuleti},
  \bibinfo{author}{A.~Dereims}, \bibinfo{author}{R.~G. Ghanem},
\newblock \bibinfo{title}{Integrated stochastic analysis of fiber composites
  manufacturing using adapted polynomial chaos expansions},
\newblock \bibinfo{journal}{Composites Part A: Applied Science and
  Manufacturing} \bibinfo{volume}{118} (\bibinfo{year}{2019})
  \bibinfo{pages}{179--193}.
%Type = Article
\bibitem[{Zeng et~al.(2021)Zeng, Red-Horse, and Ghanem}]{zeng2021accelerated}
\bibinfo{author}{X.~Zeng}, \bibinfo{author}{J.~Red-Horse},
  \bibinfo{author}{R.~Ghanem},
\newblock \bibinfo{title}{Accelerated basis adaptation in homogeneous chaos
  spaces},
\newblock \bibinfo{journal}{Computer Methods in Applied Mechanics and
  Engineering} \bibinfo{volume}{386} (\bibinfo{year}{2021})
  \bibinfo{pages}{114109}.
%Type = Article
\bibitem[{Cai et~al.(2018)Cai, Luo, Wang, and Yang}]{cai2018feature}
\bibinfo{author}{J.~Cai}, \bibinfo{author}{J.~Luo}, \bibinfo{author}{S.~Wang},
  \bibinfo{author}{S.~Yang},
\newblock \bibinfo{title}{Feature selection in machine learning: A new
  perspective},
\newblock \bibinfo{journal}{Neurocomputing} \bibinfo{volume}{300}
  (\bibinfo{year}{2018}) \bibinfo{pages}{70--79}.
%Type = Book
\bibitem[{Murphy(2012)}]{murphy2012machine}
\bibinfo{author}{K.~P. Murphy}, \bibinfo{title}{Machine learning: a
  probabilistic perspective}, \bibinfo{publisher}{MIT press},
  \bibinfo{year}{2012}.
%Type = Inproceedings
\bibitem[{Howley et~al.(2005)Howley, Madden, O’Connell, and
  Ryder}]{howley2005effect}
\bibinfo{author}{T.~Howley}, \bibinfo{author}{M.~G. Madden},
  \bibinfo{author}{M.-L. O’Connell}, \bibinfo{author}{A.~G. Ryder},
\newblock \bibinfo{title}{The effect of principal component analysis on machine
  learning accuracy with high dimensional spectral data},
\newblock in: \bibinfo{booktitle}{International Conference on Innovative
  Techniques and Applications of Artificial Intelligence},
  \bibinfo{organization}{Springer}, \bibinfo{year}{2005}, pp.
  \bibinfo{pages}{209--222}.
%Type = Article
\bibitem[{Reddy et~al.(2020)Reddy, Reddy, Lakshmanna, Kaluri, Rajput,
  Srivastava, and Baker}]{reddy2020analysis}
\bibinfo{author}{G.~T. Reddy}, \bibinfo{author}{M.~P.~K. Reddy},
  \bibinfo{author}{K.~Lakshmanna}, \bibinfo{author}{R.~Kaluri},
  \bibinfo{author}{D.~S. Rajput}, \bibinfo{author}{G.~Srivastava},
  \bibinfo{author}{T.~Baker},
\newblock \bibinfo{title}{Analysis of dimensionality reduction techniques on
  big data},
\newblock \bibinfo{journal}{IEEE Access} \bibinfo{volume}{8}
  (\bibinfo{year}{2020}) \bibinfo{pages}{54776--54788}.
%Type = Article
\bibitem[{McInnes et~al.(2018)McInnes, Healy, and Melville}]{mcinnes2018umap}
\bibinfo{author}{L.~McInnes}, \bibinfo{author}{J.~Healy},
  \bibinfo{author}{J.~Melville},
\newblock \bibinfo{title}{Umap: Uniform manifold approximation and projection
  for dimension reduction},
\newblock \bibinfo{journal}{arXiv preprint arXiv:1802.03426}
  (\bibinfo{year}{2018}).
%Type = Article
\bibitem[{Lin and Zha(2008)}]{lin2008riemannian}
\bibinfo{author}{T.~Lin}, \bibinfo{author}{H.~Zha},
\newblock \bibinfo{title}{Riemannian manifold learning},
\newblock \bibinfo{journal}{IEEE Transactions on Pattern Analysis and Machine
  Intelligence} \bibinfo{volume}{30} (\bibinfo{year}{2008})
  \bibinfo{pages}{796--809}.
%Type = Inproceedings
\bibitem[{Han et~al.(2018)Han, Wang, Zhang, Li, and Xu}]{han2018autoencoder}
\bibinfo{author}{K.~Han}, \bibinfo{author}{Y.~Wang},
  \bibinfo{author}{C.~Zhang}, \bibinfo{author}{C.~Li}, \bibinfo{author}{C.~Xu},
\newblock \bibinfo{title}{Autoencoder inspired unsupervised feature selection},
\newblock in: \bibinfo{booktitle}{2018 IEEE international conference on
  acoustics, speech and signal processing (ICASSP)},
  \bibinfo{organization}{IEEE}, \bibinfo{year}{2018}, pp.
  \bibinfo{pages}{2941--2945}.
%Type = Inproceedings
\bibitem[{Baldi(2012)}]{baldi2012autoencoders}
\bibinfo{author}{P.~Baldi},
\newblock \bibinfo{title}{Autoencoders, unsupervised learning, and deep
  architectures},
\newblock in: \bibinfo{booktitle}{Proceedings of ICML workshop on unsupervised
  and transfer learning}, \bibinfo{organization}{JMLR Workshop and Conference
  Proceedings}, \bibinfo{year}{2012}, pp. \bibinfo{pages}{37--49}.
%Type = Article
\bibitem[{Friedman and Tukey(1974)}]{friedman1974projection}
\bibinfo{author}{J.~H. Friedman}, \bibinfo{author}{J.~W. Tukey},
\newblock \bibinfo{title}{A projection pursuit algorithm for exploratory data
  analysis},
\newblock \bibinfo{journal}{IEEE Transactions on computers}
  \bibinfo{volume}{100} (\bibinfo{year}{1974}) \bibinfo{pages}{881--890}.
%Type = Article
\bibitem[{Huber(1985)}]{huber1985projection}
\bibinfo{author}{P.~J. Huber},
\newblock \bibinfo{title}{Projection pursuit},
\newblock \bibinfo{journal}{The annals of Statistics}  (\bibinfo{year}{1985})
  \bibinfo{pages}{435--475}.
%Type = Article
\bibitem[{Friedman(1987)}]{friedman1987exploratory}
\bibinfo{author}{J.~H. Friedman},
\newblock \bibinfo{title}{Exploratory projection pursuit},
\newblock \bibinfo{journal}{Journal of the American statistical association}
  \bibinfo{volume}{82} (\bibinfo{year}{1987}) \bibinfo{pages}{249--266}.
%Type = Article
\bibitem[{Lee et~al.(2005)Lee, Cook, Klinke, and Lumley}]{lee2005projection}
\bibinfo{author}{E.-K. Lee}, \bibinfo{author}{D.~Cook},
  \bibinfo{author}{S.~Klinke}, \bibinfo{author}{T.~Lumley},
\newblock \bibinfo{title}{Projection pursuit for exploratory supervised
  classification},
\newblock \bibinfo{journal}{Journal of Computational and graphical Statistics}
  \bibinfo{volume}{14} (\bibinfo{year}{2005}) \bibinfo{pages}{831--846}.
%Type = Inproceedings
\bibitem[{Grochowski and Duch(2008)}]{grochowski2008projection}
\bibinfo{author}{M.~Grochowski}, \bibinfo{author}{W.~Duch},
\newblock \bibinfo{title}{Projection pursuit constructive neural networks based
  on quality of projected clusters},
\newblock in: \bibinfo{booktitle}{International Conference on Artificial Neural
  Networks}, \bibinfo{organization}{Springer}, \bibinfo{year}{2008}, pp.
  \bibinfo{pages}{754--762}.
%Type = Article
\bibitem[{Barcaru(2019)}]{barcaru2019supervised}
\bibinfo{author}{A.~Barcaru},
\newblock \bibinfo{title}{Supervised projection pursuit--a dimensionality
  reduction technique optimized for probabilistic classification},
\newblock \bibinfo{journal}{Chemometrics and Intelligent Laboratory Systems}
  \bibinfo{volume}{194} (\bibinfo{year}{2019}) \bibinfo{pages}{103867}.
%Type = Article
\bibitem[{Grear et~al.(2021)Grear, Avery, Patterson, and
  Jacobs}]{grear2021molecular}
\bibinfo{author}{T.~Grear}, \bibinfo{author}{C.~Avery},
  \bibinfo{author}{J.~Patterson}, \bibinfo{author}{D.~J. Jacobs},
\newblock \bibinfo{title}{Molecular function recognition by supervised
  projection pursuit machine learning},
\newblock \bibinfo{journal}{Scientific reports} \bibinfo{volume}{11}
  (\bibinfo{year}{2021}) \bibinfo{pages}{1--15}.
%Type = Article
\bibitem[{Olivier et~al.(2021)Olivier, Shields, and
  Graham-Brady}]{olivier2021bayesian}
\bibinfo{author}{A.~Olivier}, \bibinfo{author}{M.~D. Shields},
  \bibinfo{author}{L.~Graham-Brady},
\newblock \bibinfo{title}{Bayesian neural networks for uncertainty
  quantification in data-driven materials modeling},
\newblock \bibinfo{journal}{Computer Methods in Applied Mechanics and
  Engineering} \bibinfo{volume}{386} (\bibinfo{year}{2021})
  \bibinfo{pages}{114079}.
%Type = Article
\bibitem[{Yang et~al.(2021)Yang, Meng, and Karniadakis}]{yang2021b}
\bibinfo{author}{L.~Yang}, \bibinfo{author}{X.~Meng}, \bibinfo{author}{G.~E.
  Karniadakis},
\newblock \bibinfo{title}{B-pinns: Bayesian physics-informed neural networks
  for forward and inverse pde problems with noisy data},
\newblock \bibinfo{journal}{Journal of Computational Physics}
  \bibinfo{volume}{425} (\bibinfo{year}{2021}) \bibinfo{pages}{109913}.
%Type = Article
\bibitem[{Leibig et~al.(2017)Leibig, Allken, Ayhan, Berens, and
  Wahl}]{leibig2017leveraging}
\bibinfo{author}{C.~Leibig}, \bibinfo{author}{V.~Allken},
  \bibinfo{author}{M.~S. Ayhan}, \bibinfo{author}{P.~Berens},
  \bibinfo{author}{S.~Wahl},
\newblock \bibinfo{title}{Leveraging uncertainty information from deep neural
  networks for disease detection},
\newblock \bibinfo{journal}{Scientific reports} \bibinfo{volume}{7}
  (\bibinfo{year}{2017}) \bibinfo{pages}{1--14}.
%Type = Article
\bibitem[{Friedman and Stuetzle(1981)}]{friedman1981projection}
\bibinfo{author}{J.~H. Friedman}, \bibinfo{author}{W.~Stuetzle},
\newblock \bibinfo{title}{Projection pursuit regression},
\newblock \bibinfo{journal}{Journal of the American statistical Association}
  \bibinfo{volume}{76} (\bibinfo{year}{1981}) \bibinfo{pages}{817--823}.
%Type = Book
\bibitem[{Hastie et~al.(2009)Hastie, Tibshirani, Friedman, and
  Friedman}]{hastie2009elements}
\bibinfo{author}{T.~Hastie}, \bibinfo{author}{R.~Tibshirani},
  \bibinfo{author}{J.~H. Friedman}, \bibinfo{author}{J.~H. Friedman},
  \bibinfo{title}{The elements of statistical learning: data mining, inference,
  and prediction}, volume~\bibinfo{volume}{2}, \bibinfo{publisher}{Springer},
  \bibinfo{year}{2009}.
%Type = Article
\bibitem[{Qianjian and Jianguo(2011)}]{qianjian2011application}
\bibinfo{author}{G.~Qianjian}, \bibinfo{author}{Y.~Jianguo},
\newblock \bibinfo{title}{Application of projection pursuit regression to
  thermal error modeling of a cnc machine tool},
\newblock \bibinfo{journal}{The International Journal of Advanced Manufacturing
  Technology} \bibinfo{volume}{55} (\bibinfo{year}{2011})
  \bibinfo{pages}{623--629}.
%Type = Article
\bibitem[{Ferraty et~al.(2013)Ferraty, Goia, Salinelli, and
  Vieu}]{ferraty2013functional}
\bibinfo{author}{F.~Ferraty}, \bibinfo{author}{A.~Goia},
  \bibinfo{author}{E.~Salinelli}, \bibinfo{author}{P.~Vieu},
\newblock \bibinfo{title}{Functional projection pursuit regression},
\newblock \bibinfo{journal}{Test} \bibinfo{volume}{22} (\bibinfo{year}{2013})
  \bibinfo{pages}{293--320}.
%Type = Article
\bibitem[{Durocher et~al.(2015)Durocher, Chebana, and
  Ouarda}]{durocher2015nonlinear}
\bibinfo{author}{M.~Durocher}, \bibinfo{author}{F.~Chebana},
  \bibinfo{author}{T.~B. Ouarda},
\newblock \bibinfo{title}{A nonlinear approach to regional flood frequency
  analysis using projection pursuit regression},
\newblock \bibinfo{journal}{Journal of Hydrometeorology} \bibinfo{volume}{16}
  (\bibinfo{year}{2015}) \bibinfo{pages}{1561--1574}.
%Type = Article
\bibitem[{Cui et~al.(2017)Cui, Zhao, Chen, Zhang, Wang, Lu, Jia, and
  Wei}]{cui2017assessment}
\bibinfo{author}{H.-Y. Cui}, \bibinfo{author}{Y.~Zhao}, \bibinfo{author}{Y.-N.
  Chen}, \bibinfo{author}{X.~Zhang}, \bibinfo{author}{X.-Q. Wang},
  \bibinfo{author}{Q.~Lu}, \bibinfo{author}{L.-M. Jia}, \bibinfo{author}{Z.-M.
  Wei},
\newblock \bibinfo{title}{Assessment of phytotoxicity grade during composting
  based on eem/parafac combined with projection pursuit regression},
\newblock \bibinfo{journal}{Journal of hazardous materials}
  \bibinfo{volume}{326} (\bibinfo{year}{2017}) \bibinfo{pages}{10--17}.
%Type = Book
\bibitem[{Janson et~al.(1997)}]{janson1997}
\bibinfo{author}{S.~Janson}, et~al., \bibinfo{title}{Gaussian Hilbert Spaces},
  \bibinfo{publisher}{Cambridge University Press}, \bibinfo{year}{1997}.
%Type = Article
\bibitem[{Cohen and Migliorati(2017)}]{cohen2017optimal}
\bibinfo{author}{A.~Cohen}, \bibinfo{author}{G.~Migliorati},
\newblock \bibinfo{title}{Optimal weighted least-squares methods},
\newblock \bibinfo{journal}{The SMAI journal of computational mathematics}
  \bibinfo{volume}{3} (\bibinfo{year}{2017}) \bibinfo{pages}{181--203}.
%Type = Article
\bibitem[{Schwab and Gittelson(2011)}]{schwab2011}
\bibinfo{author}{C.~Schwab}, \bibinfo{author}{C.~J. Gittelson},
\newblock \bibinfo{title}{Sparse tensor discretizations of high-dimensional
  parametric and stochastic {PDEs}},
\newblock \bibinfo{journal}{Acta Numerica} \bibinfo{volume}{20}
  (\bibinfo{year}{2011}) \bibinfo{pages}{291--467}.
%Type = Article
\bibitem[{Rosenblatt(1952)}]{rosenblatt1952remarks}
\bibinfo{author}{M.~Rosenblatt},
\newblock \bibinfo{title}{Remarks on a multivariate transformation},
\newblock \bibinfo{journal}{The annals of mathematical statistics}
  \bibinfo{volume}{23} (\bibinfo{year}{1952}) \bibinfo{pages}{470--472}.
%Type = Article
\bibitem[{Tsilifis and Ghanem(2017)}]{tsilifis2017}
\bibinfo{author}{P.~Tsilifis}, \bibinfo{author}{R.~G. Ghanem},
\newblock \bibinfo{title}{Reduced {Wiener} chaos representation of random
  fields via basis adaptation and projection},
\newblock \bibinfo{journal}{Journal of Computational Physics}
  \bibinfo{volume}{341} (\bibinfo{year}{2017}) \bibinfo{pages}{102--120}.

\end{thebibliography}

\end{document}